\documentclass{elsart}

\usepackage{algorithmic}
\usepackage{algorithm}

\usepackage{amsmath}
\usepackage{amssymb}
\usepackage{ifthen}

\newcommand{\WR}{semi-reachable}

\newcommand{\SPAN}{\mathrm{Span}}

\newcommand \IO[1][{(Q \times T)^{+}}]
{ F(PC(T,\mathcal{U}) \times #1, \mathcal{Y}) }

\newcommand{\SwitchSysLin}[1][{}]{
  (\mathcal{X}_{#1},\mathcal{U},\mathcal{Y},Q,\{(A_{q}^{#1},B_{q}^{#1},C_{q}^{#1}) \mid q \in Q \})
}

\newcommand{\IM}{\mathrm{Im}}
\newcommand{\Rank}{\mathrm{rank}\mbox{ } }
\newtheorem{Definition}{Definition}
\newtheorem{Theorem}{Theorem}
\newtheorem{Lemma}{Lemma}

\newtheorem{Corollary}{Corollary}

\theoremstyle{definition}
\newtheorem{Remark}{Remark}
\newtheorem{Notation}{Notation}

\newtheorem{Construction}{Construction}

\newtheorem{Problem}{Problem}

\newcommand{\sep}{{\ }}

\newcommand{\LSS}{LSS}
\newcommand{\SLSS}{LSSs}
\newcommand{\MORPH}{\mathcal{S}}
\newcommand{\QNUM}{D}
\newcommand{\RMORPH}{\mathcal{S}}

\begin{document} 
\newenvironment{proof}{\begin{pf}}{\end{pf}}
\newcommand{\COL}{\mathbf{C}}
\bibliographystyle{plain}
\begin{frontmatter}
\title{ Partial-realization theory and algorithms \\
        for linear
         switched systems\\
         A formal power series approach }
\author{ Mih\'aly Petreczky}
\address{
Maastricht University
P.O. Box 616, 6200 MD Maastricht, The Netherlands 
 \texttt{M.Petreczky@maastrichtuniversity.nl}
  } 
\author{ Jan H. van Schuppen }
\address{Centrum voor Wiskunde en Informatica (CWI),\\
         P.O. Box 94079, 1090\ GB Amsterdam, The Netherlands \\
         \texttt{J.H.van.Schuppen@cwi.nl}
}
\begin{abstract}

The paper presents partial-realization theory  and
realization algorithms for linear switched systems. 
Linear switched systems are a particular subclass of hybrid systems.
We formulate a notion of a 
partial realization and we present conditions for existence of a minimal partial realization.
We propose
two partial-realization algorithms and we show that under certain conditions they
yield a complete realization.
Our main tool is the theory of rational formal power series. 
\end{abstract}
\end{frontmatter}
\section{Introduction}
\label{sect:intro}

\textbf{The main objective of the paper} \\
The immediate objective of the paper is to present partial-realization theory of
linear switched systems. However, the broader goal is to demonstrate that
\begin{itemize}
\item[(a)]
    it is possible to develop partial-realization theory for hybrid systems,
\item[(b)]
   the theory of rational formal power series can be used to obtain partial-realization
   theory of hybrid systems.
\end{itemize}
 The results of the paper just serve as evidence for the above claims.
 In fact, partial-realization theory can be developed along the same line for other classes of
 hybrid systems, \cite{MP:Phd}.
This paper is intended as the first one in a series of papers dealing with
partial-realization theory of hybrid systems.

\textbf{The class of linear switched systems } \\
  Linear switched systems have been studied for almost two decades, 
  see \cite{D:Lib,Sun:Book} for a survey.  Their practical and theoretical relevance is
  widely recognized. 
 \emph{Linear switched systems (abbreviated by \LSS)} are hybrid systems with external
  switching only, whose continuous dynamics in each discrete mode is determined
  by a linear continuous-time system, and whose discrete dynamics is trivial, i.e.
  any discrete state transition is allowed. The continuous subsystems are defined on
  the same state-, input- and output-spaces. 
  A discrete state transition occurs if the environment enforces one, i.e.
  the switching sequence itself functions as an input. 
  The time evolution of a linear switched system can be pictured as follows. A discrete mode
  is selected and a continuous input is fed into the corresponding linear subsystem.
  The state and output of the linear switched system is now determined by the time
  evolution of the chosen linear subsystem. At some point in time, a new discrete mode
  is selected. Then the linear subsystem corresponding to the new discrete mode is
  launched, using as initial state the state of the previous linear subsystem at the
  moment when a transition to the new mode occurred. From this point on, the state and
  output evolve according to the newly selected linear subsystem.

\textbf{Motivation} \\
  We believe that similarly to linear systems,  partial-realization 
  theory for \SLSS\  will be  useful for systems identification and model reduction of
  \SLSS\ and even of more general hybrid systems. We will
  elaborate on the possible applications in Section \ref{sect:concl}.


\textbf{Contribution of the paper} \\
The paper presents the following results.
\begin{itemize}
\item{\textbf{Partial-realization algorithm and theorem for \SLSS}} \\
   In Section \ref{sect:main_result} we propose a notion of partial-realization. 
   In Theorem \ref{part_real_lin:theo1} and \ref{part_real_lin:theo1.1} we formulate
   a sufficient condition for existence and minimality of a partial realization.
   We present two algorithms for 
  computing a partial-realization.
  In one of them, Algorithm \ref{alg0}, a partial-realization is constructed from the
  columns of a \emph{finite} sub-matrix Hankel-matrix, the other one, Algorithm \ref{PartReal}, 
  is based on finding a factorization of
  a \emph{finite} sub-matrix of the Hankel-matrix.
  The outcomes of both procedures are isomorphic.
  The former is potentially useful for theoretical purposes, while the latter one might serve as a basis for subspace identification-like methods.
  The factorization algorithm was implemented and the implementation
  is available from the the first-named author on request. 
\item{\textbf{Realization algorithm for \SLSS}}\\
 In Theorem \ref{part_real_lin:theo1} and \ref{part_real_lin:theo2} 
 we show that the above partial-realization algorithms return a complete
 minimal realization of the input-output behavior, provided the rank of the 
 finite  sub-matrix of the Hankel-matrix equals the rank of the full Hankel-matrix.
 We show that any finite sub-matrix whose size is above a certain lower bound has this
 property. 
\item{\textbf{Partial-realization theory for families
      of formal power series}} \\
  Both the (partial-) realization algorithm and the corresponding theoretical results
  for \SLSS\  follow quite easily from analogous statements for
  families of rational formal series. 
  In Appendix \ref{sect:pow:part} 
  we present what can be termed as partial-realization theory and algorithms
  for \emph{families of formal power series}.
  Although partial-realization theory of a 
  single rational formal power series has been known in various forms
  \cite{Isi:Tac,BilPart,Isi:Bilin,Son:Resp,Son:Real,BilSubMoor,MacChenBilSub},
  its extension
  to families of formal power series appears to be new. As it was already noted,
  this extension of the classical theory can be then used to develop partial-realization
  theory for a number of other (more general) classes of hybrid systems, \cite{MP:Phd}.
\end{itemize}

\textbf{Formal power series approach} \\
Our main tool is the theory of
rational formal power series.
Recall from \cite{MP:BigArticle,MP:RealForm,MP:RealBilin} that
there is a correspondence between 
\LSS\  realizations
and representations of certain rational families of former power series.
Recall that the  Hankel-matrix of a family of input-output maps
is in fact the Hankel-matrix of the corresponding family of
formal power series. Hence, by formulating a partial-realization theory and
algorithms for rational families of formal power series,
we immediately obtain partial-realization theory and algorithms for
\SLSS.

\textbf{Related work} \\
 To the best of our knowledge, the only results on partial-realization of linear
 switched systems is the first-named author's thesis \cite{MP:Phd}. In
 \cite{MP:RealBilin,MP:HybBilinReal} partial-realization theory of
 bilinear hybrid and switched systems are announced, but no proofs are provided.
 The thesis \cite{MP:Phd} contains most of
 the results of this paper. In addition, it also covers partial-realization theory of
 several other classes of hybrid systems.
 Reference \cite{MPRV:JumpMarkov} announces some results on partial-realization theory of
 stochastic jump-Markov systems.

Realization theory of rational formal power series is a classical topic, see
\cite{Son:Resp,Son:Real,Isi:Nonlin,Reut:Book,Salomaa:Book,MFliessHank,MFliessFormPow}.
It is known to be closely related to
 realization theory of bilinear and state-affine
 systems \cite{Isi:Tac,Isi:Bilin,Son:Resp,Son:Real}.
In turn, results on
partial-realization theory for discrete-time bilinear systems and state-affine systems
can be found in \cite{BilPart,Isi:Tac,Son:Resp,Son:Real,BilSubMoor,MacChenBilSub}. 
With respect to \cite{Isi:Tac,BilPart,Isi:Bilin,Son:Resp,Son:Real,BilSubMoor,MacChenBilSub}, the main
novelty of the paper is the following.
\begin{enumerate}
\item  In this paper (in Appendix \ref{sect:pow:part})
       partial-realization theory is stated directly for rational formal power series 
       representations, without reformulating it in terms of some particular system class
       such as bilinear or state-affine systems.
\item While the classical results \cite{Isi:Tac,BilPart,Isi:Bilin,Son:Resp,Son:Real,BilSubMoor,MacChenBilSub}
      can be thought of as a reformulation of partial-realization
      theory for a single rational formal power series, here we deal with \emph{families
      of formal power series}. 
\item We apply partial-realization theory of rational formal power series to
      linear switched systems.  This represents a novel result
      in the theory of hybrid systems.
      We view this as the main contribution of the paper.
\end{enumerate}
The statement of the main results on partial-realization theory for
families of formal power series was announced \cite{MP:HybPow,MPRV:JumpMarkov,MP:HybBilinReal,MP:RealBilin,MP:HybLinBilinTechReport,PetreczkyV07}, but no
proof was ever presented. The thesis \cite{MP:Phd} contains the statement and
the proof of the main results on formal power series.

\textbf{Outline of the paper} \\
The outline of the paper is the following. 
Section \ref{sect:info} presents an informal formulation of the partial-realization
problem for \SLSS. Section \ref{sect:num} discusses a numerical example. It is intended
as an accessible explanation of the main results of the paper by means of an example.
Section \ref{sect:switch} presents the definition and some elementary properties of linear 
switched systems together with the basic notation and terminology. 
Section \ref{sect:lin_sys:real_theo} provides a brief overview of the realization theory
of linear switched systems. This section is a prerequisite for
understanding the main results of the paper. 
Section \ref{sect:main_result} presents the formal statement of the main results of
the paper.  Section \ref{sect:lin_alg} presents the proof of the results stated in
Section \ref{sect:main_result}.
Finally, in Section \ref{sect:concl} we formulate the conclusions of the paper and
we discuss potential applications of the presented results.
As it was already remarked, partial-realization theory of \SLSS\ is based on 
partial-realization theory of rational formal power series. 
Therefore, in Appendix \ref{sect:pow}
we included a brief overview of 
the extension of the theory of formal power series to families of formal power series. 
More precisely, Appendix \ref{sect:pow}
reviews the relevant concepts and results on formal power series
from \cite{MP:BigArticle,MP:RealBilin,MP:RealForm}. In Appendix \ref{sect:pow:part}
we present partial-realization theory of families of formal power series.
and algorithms for computing a 
representation for a family of formal power series.
Appendix \ref{sect:pow:part} is in fact the main technical tool of the paper. It is a 
prerequisite for understanding the proof of the main results on partial-realization
theory of \SLSS.
\section{ Informal problem formulation }
\label{sect:info}

The goal of this section is to state the partial-realization problem for linear switched
systems in an informal way.
The next section, Section \ref{sect:num},
 provides a numerical example to illustrate the problem and the solution.
We defer the formal problem statement until Section \ref{sect:main_result}.


  \textbf{ Brief review of linear-partial realization theory} \\
  The partial-realization problem was originally formulated for linear  systems.
\cite{MR0255260,MR0245360,TetherPartLin} as follows.
  Assume that the first $N \in \mathbb{N}$ 
  Markov parameters $S=\{M_i\}_{i=1}^{N}$ of an input-output
  map $f$ are specified.
  \begin{itemize}
  \item[(a)]
    Find conditions for existence of a linear system
    whose first $N$ Markov-parameters coincide with $S=\{M_{i}\}_{i=1}^{N}$.
    Such a linear system is called a 
    \emph{partial-realization of $S$}. 
    Characterize minimal partial-realizations. Find an algorithm for 
    computing a (minimal) partial-realization.
  \item[(b)]
    Find conditions under which the thus obtained linear system is a minimal realization
    of the input-output map $f$.    
  \end{itemize}
   Notice that for discrete-time linear systems, the $k$th Markov parameter
   coincides with the output of the system at time $k$ for a particular input.
   Hence, the partial-realization problem can be seen as an identification problem.
   In fact, partial-realization theory of linear systems can be used
   for systems identification and model reduction.

  \textbf{Markov-parameters for \SLSS} \\
   Let $\Phi$ be a family of input-output map which map continuous-valued inputs and
   switching sequences to continuous outputs. That is, the elements of $\Phi$ are 
   of the same form as the input-output maps generated by \SLSS.
   Below  we will define \emph{generalized Markov-parameters}
   of a family $\Phi$ of input-output maps which could potentially be realized by a \LSS.  
   The generalized Markov-parameters of $\Phi$ are defined as certain high-order 
   derivatives with respect to the switching times of the elements of $\Phi$.
   More precisely, the Markov-parameters are indexed by triples consisting of the
   following components
   \begin{itemize}
   \item[(a)] elements of $\Phi$, or
               pairs $(q,j)$ where $q$ runs through the set of discrete modes and
              $j$ runs through the set $\{1,2,\ldots,m\}$, where $m$ denotes the
              number of continuous-valued input channels, and
   \item[(b)] sequences of discrete modes, and
   \item[(c)] discrete modes.
   \end{itemize}
   Without going into details, the
   intuition behind the presented indexing is the following. A Markov-parameter
   of $\Phi$ indexed by $f \in \Phi$, a discrete mode $q$, and
   a sequence of discrete modes $q_1,q_2,\ldots,q_k$, stands for the partial derivative
   of $f$ with respect to the switching times, evaluated at zero, for the 
   following switching scenario. 
   The system goes through the modes $q_1,q_2,\cdots, q_k$ and then jumps to $q$.  
   Here the continuous input is set to zero. 
   The intuitive meaning of 
   a Markov-parameter indexed by a pair $(q_0,j)$, a discrete mode $q$ and
   a sequence of discrete modes $q_1,q_2,\ldots,q_k$ is the following.
    Each input-output map $f \in \Phi$ can be written as a sum $f=a_{f}+y^{\Phi}$ of two maps, where $a_f$ is independent of the continuous input and $y^{\Phi}$
   is common for all the elements of $\Phi$ and it is linear in continuous inputs.
   Roughly speaking, $a_f$ accounts for the output from a certain initial state
   under zero input, and $y^{\Phi}$ represents the input-output map induced by the
   zero initial state. The Markov-parameter indexed by the pair $(q_0,j)$, discrete
   mode $q$ and sequence $q_1,q_2,\cdots q_k$ is the derivative of $y^{\Phi}$ 
   with respect
   to the switching times evaluated at zero, for the switching scenario where
   the system goes from mode $q_0$ to $q_1, q_2,\cdots q_k$ and ends in $q$, and
   all the continuous input channels are $0$ except the
   $j$th one which is $1$.

  \textbf{Partial-realization problem for \SLSS} \\
    We will refer to a Markov parameter indexed by a sequence of discrete modes of
    length $k$ as a Markov-parameter of \emph{order $k$}. Fix a natural number $N>0$
    and let $S$ be the collection of all Markov-parameters of $\Phi$
    of order at most $N$. That is, $S$ is simply a collection of high-order
    derivatives of the elements of $\Phi$, such that the degree of derivation is bounded by $N$.
    A \LSS\  $\Sigma$ is said 
   to be a \emph{$N$-partial realization} of $\Phi$, if certain products of 
   the matrices of $\Sigma$ are equal to the corresponding elements of $S$. 
   In other words, $\Sigma$ is a $N$-partial realization of $\Phi$ if
   the input-output maps of $\Sigma$ and those of $\Phi$ have the property that their
   derivatives corresponding to the Markov-parameters of order at most $N$ coincide.
   From realization theory of \SLSS\ \cite{MP:Phd,MP:RealForm,MP:BigArticle} it follows
   that a \LSS\ is a realization of $\Phi$ if and only if it is a $N$-partial realization of
   $\Phi$ for all $N \in \mathbb{N}$.
  The \emph{partial realization problem} for \SLSS\  can be now stated as follows.
  \begin{itemize}
     \item
            Find conditions for existence of a $N$-partial realization of $\Phi$ by \LSS.
            Characterize minimal dimensional $N$-partial realizations of $\Phi$.
            Find an algorithm for computing a minimal $N$-partial \LSS\  realizations of $\Phi$.
     \item
            Find conditions under which a minimal $N$-partial \LSS\ realization of $\Phi$
            is a complete realization of $\Phi$.
     \end{itemize}

   The motivation for studying the partial-realization problem for \SLSS\  is similar
   to that of for linear systems, i.e. we expect it to be
   useful for model reduction and systems identification. In Section \ref{sect:concl}
   we will present a more detailed description of the motivation and possible 
   applications. 

\section{Numerical example}
\label{sect:num}
The purpose of this section is to demonstrate the main results of the paper by means of a 
numerical example. In this section we will tacitly use the notation and terminology
of Section \ref{sect:switch} and Section \ref{sect:main_result}. 

Consider the linear switched system of the form
\begin{equation}
\label{lin_switch01}
\Sigma \left \{ \begin{split}
   \dot x(t)=& A_{q(t)}x(t)+B_{q(t)}u(t)  \\
   y(t)=&C_{q(t)}x(t)
\end{split}
\right.
\end{equation}
where 
$q(t) \in \{q_1,q_2\}$ is the discrete mode at time $t$,
$x(t) \in \mathbb{R}^{5}$ is the continuous-state at $t$,
$y(t) \in \mathbb{R}$ is the scalar output at $t$, and
$u(t) \in \mathbb{R}$ is the scalar input at $t$.
The system matrices $A_q,B_q,C_q$ describing the linear control system residing in a 
state (mode) $q \in \{q_1,q_2\}=Q$ are of the following form. 
\[ 
  \begin{split}
  & A_{q_1} =
   \begin{bmatrix} 
     0 &  0 &  0 &  0 &  0 \\
    0 &  0 &  1 &  0 &  0 \\
    0 &  0 &  0 &  1 &  0 \\
    0 & 0 &  0 &  0 &  1  \\
    0 & 0 &  0 &  0 &  0
  \end{bmatrix},
  B_{q_1}=\begin{bmatrix} 0 \\ 0 \\ 0 \\ 0 \\ 1 \end{bmatrix}, 
  C_{q_1}=\begin{bmatrix} 0 \\  1 \\  0 \\  0 \\  0 \end{bmatrix}^{T} \\
 & A_{q_2}=
 \begin{bmatrix} 0 &  0 &  0 &  0 &  0 \\
                 0 &  2 &  0 &  0 &  0 \\
                 0 &  0 &  0 &  0 &  0 \\
                 0 &  0 &  0 &  0 &  0 \\
                 0 &  0 &  0 &  0 &  3
  \end{bmatrix},
  B_{q_2}=\begin{bmatrix} 0 \\ 1 \\ 0 \\ 0 \\ 0 \end{bmatrix}, 
  C_{q_2}=\begin{bmatrix} 0 \\  0 \\  0 \\  0 \\  1 \end{bmatrix}^{T}
\end{split}
\]
Consider the initial states $x_1=0$ and 
$x_2=\begin{bmatrix} 0 & 0 & 0 & 0 & 1 \end{bmatrix}^{T}$.
Consider the set of input-output maps $\Phi=\{f_1,f_2\}$ such that
$f_i$ is realized by $\Sigma$ from the initial state $x_i$, $i=1,2$, i.e.
$f_i(u,w)=y_{\Sigma}(x_i,u,w)$, for $i=1,2$, for each continuous-valued input
$u \in PC(T,\mathcal{U})$ and finite switching sequence $w \in (Q \times T)^{+}$.

Notice that $\Phi$ can be realized by the following minimal linear switched
system
\begin{equation}
\label{lin_switch02}
\Sigma_m \left \{ \begin{split}
   \dot x(t)=& A^{m}_{q(t)}x(t)+B^{m}_{q(t)}u(t)  \\
   y(t)=&C^{m}_{q(t)}x(t)
\end{split}
\right.
\end{equation}
 where for each $q \in \{q_1,q_2\}$, the matrices
 $A_{q}^{m},B_{q}^{m},C_{q}^{m}$ are of the following form
 \[
  \begin{split}
 &  A^{m}_{q_1}=
  \begin{bmatrix} 0 &  0 &  0 &  0 \\
   0 &  0 &  0 &  1 \\
   1 &  0 &  0 &  0 \\
   0 &  0 &  1 &  0
 \end{bmatrix},
 B^{m}_{q_1}=\begin{bmatrix} 1 \\ 0 \\ 0 \\ 0 \end{bmatrix},
 C^{m}_{q_1}=\begin{bmatrix} 0 \\ 1 \\  0 \\  0 \end{bmatrix}^{T} \\
 & A^{m}_{q_2}=\begin{bmatrix}
   3 &   0 &  0 &   0 \\
   0 &  2  &  0 &   0 \\
   0 &  0  &  0 &   0 \\
   0 &  0  &  0 & 0 
  \end{bmatrix}, 
   B^{m}_{q_2}=\begin{bmatrix} 0 \\ 1 \\ 0 \\ 0 \end{bmatrix},
   C^{m}_{q_2}=\begin{bmatrix} 1 \\ 0 \\ 0 \\ 0 \end{bmatrix}^{T}
 \end{split}
\]
More precisely, $\Sigma_m$ realizes $f_i$ from the 
initial states $x^{m}_i$, for $i=1,2$. Here, $x_1^{m}=0$ and
$x_2^{m}=\begin{bmatrix} 1 & 0 & 0 & 0 \end{bmatrix}^{T}$.

Since $\Phi$ has a realization by a \LSS, it is clear that it has a generalized
kernel representation and hence its Markov-parameters can be defined.
In order to give a better intuition on Markov-parameters of $\Phi$, we have listed some
of them in Table \ref{table1}.
\begin{table}
\caption{\label{table1} Markov-parameters of $\Phi$}
\begin{tabular}{|l|c|c|r|} 
 \hline 
 $(q_0,j) \in (\{q_1,q_2\} \times \{1\}) $ & $q \in Q$ & $w \in Q^{*} $ & Markov-parameter  $S_{q,q_0,j}(w)$ \\ 
\hline 
$(q_2, 1)$ &  $q_2$ & $\epsilon$ &  $0$ \\
$(q_2,1)$  & $q_2$ &  $q_1$ & $0$ \\
\vdots & \vdots & \vdots & \vdots \\
$(q_2,1)$ & $q_2$ &  $q_2q_2q_2$ &  $0$ \\ 
$(q_2, 1)$ & $q_1$ & $\epsilon$ & $1$ \\
 $(q_2,1)$ & $q_1$ &  $q_1$ &  $0$ \\
 $(q_2,1)$ & $q_1$ & $q_2$ & $2$ \\
\vdots & \vdots & \vdots & \vdots \\
 $(q_2,1)$ & $q_1$ & $q_2q_2$ &  $4$ \\ 
 $(q_2,1)$ & $q_1$ & $q_2q_2q_2$ & $8$ \\
$(q_1, 1)$ &  $q_1$ & $\epsilon$ & $0$ \\
$(q_1,1)$ & $q_1$ &  $q_1$  & $0$ \\
\vdots & \vdots & \vdots & \vdots \\
$(q_1,1)$ & $q_1$ & $q_2q_2q_2$ & 0 \\
$(q_1, 1)$ & $q_2$ & $\epsilon$ &  1 \\
 $(q_1,1)$ & $q_2$ & $q_1$ &  $0$ \\
 $(q_1,1)$ & $q_2$ & $q_2$ & $3$ \\
 $(q_1,1)$ & $q_2$ & $q_2q_2$ & $9$ \\
$(q_1,1)$ & $q_2$ & $q_2q_2q_2$ & $27$ \\
 \hline 
 $j \in \Phi $ & $q \in Q$ & $w \in Q^{*} $ & Markov-parameter  $S_{j,q}(w)$ \\ 
\hline
$f_1$  & $q_1$ & $\epsilon$ & $0$ \\
$f_1$ & $q_1$ & $q_1$ & $0$ \\
$f_2$ & $q_1$ & $\epsilon$ & $0$ \\
$f_2$ & $q_1$ & $q_1$ & $0$ \\
$f_2$ & $q_2$ & $\epsilon$ & $1$ \\
$f_2$ & $q_2$ & $q_1$ & $0$ \\
$f_2$ & $q_2$ & $q_2$ & $3$ \\
$f_2$ & $q_2$ & $q_2q_2$ & $9$ \\
$f_2$ & $q_2$ & $q_2q_2q_2$ & $27$ \\
\hline
\end{tabular}
\end{table}
Consider the upper-left sub-matrix $H_{\Phi,K,L}$ of the Hankel-matrix $H_{\Phi}$ of
$\Phi$. Recall that $H_{\Phi,K,L}$ is formed by the intersection of
 the columns of $H_{\Phi}$ indexed
by a sequence of discrete modes of length at most $L$, and by the rows of $H_{\Phi}$ 
indexed by a sequence of discrete modes of length at most $K$.
From Theorem \ref{part_real_lin:theo1} it follows
that $\Rank H_{\Phi,K,L} \le \dim \Sigma_m=4$ for
all $K,L \ge 4$. 
In fact, it turns out that the Hankel-matrix
$H_{\Phi,N,N}$ for $N=2$ has already rank $4$.
By Theorem \ref{part_real_lin:theo1}
it means that 
we can already compute a minimal \LSS\  realization of $\Phi$
from the generalized Markov parameters of $\Phi$ of 
indexed by sequences of discrete modes 
of length at most $5=2+3$.
Applying Algorithm \ref{PartReal} to $H_{\Phi,3,2}$ yields
the following minimal \LSS\ realization $\Sigma_{f}$ of $\Phi$. 
\begin{equation}
\label{lin_switch020}
\Sigma_{f}: \left \{ \begin{split}
   \dot x(t)=& A^{f}_{q(t)}x(t)+B^{f}_{q(t)}u(t)  \\
   y(t)=&C^{f}_{q(t)}x(t)
\end{split}
\right.
\end{equation}
 where the system matrices are of the form
\[
  \begin{split}
  & A_{q_1}^{f}=
    \begin{bmatrix}  0 &    0 &   0 &    0 \\
                     0 &    0 &    0 &    0.79 \\
                 -0.49 &    0 &    0 &   0  \\ 
                     0 &    0 &   -2.55  & 0 
    \end{bmatrix},              
  B_{q_1}^{f} = \begin{bmatrix} 1.46 \\ 0 \\  0 \\ 0 \end{bmatrix},
  C_{q_1}^{f} = \begin{bmatrix} 0 \\    0.71 \\  0 \\    0 \end{bmatrix}^{T} \\ 
  & A_{q_2}^{f} = \begin{bmatrix}
                 3 &  0 & 0 & 0 \\
                 0 &  2 & 0  & 0 \\
                 0 &  0 & 0 &  0 \\
                 0 &  0 & 0 &  0
                 \end{bmatrix},
    B_{q_2}^{f}=\begin{bmatrix}  0 \\  1.42 \\ 0  \\ 0 \end{bmatrix},
    C_{q_2}^{f}=\begin{bmatrix}  0.69 \\  0 \\   0 \\ 0 \end{bmatrix}
  \end{split}
\]
The \LSS\  $\Sigma_{f}$ realizes the input-output map $f_1$ from the initial state
$x^{f}_1=(0,0,0,0)^{T}$, and the map $f_2$ from the initial state
$x^{f}_2=( 1.46,0,0,0)^{T}$.

However, it turns out that $\Rank H_{\Phi,1,0}=\Rank H_{\Phi,0,1}=\Rank H_{\Phi,0,0}$.
That is, the sequence of Markov parameters of $\Phi$ indexed by
sequences of discrete modes of length at most $1$ already satisfies 
the partial realization theorems Theorem \ref{part_real_lin:theo1} and 
Theorem \ref{part_real_lin:theo2}.
Hence, we can apply Algorithm \ref{PartReal} to $H_{\Phi,1,0}$ to obtain a 
minimal \LSS\ $1$-partial realization of $\Phi$. However, it can
be checked that the thus obtained \LSS\ \emph{is not a \LSS\ realization of $\Phi$}.
Indeed, by applying Algorithm \ref{PartReal} to $H_{\Phi,1,0}$ we obtain the following
\LSS\  realization of $\Phi$
\begin{equation}
\label{lin_switch03}
\Sigma_{part}: \left \{ \begin{split}
   \dot x(t)=& A^{part}_{q(t)}x(t)+B^{part}_{q(t)}u(t)  \\
   y(t)=&C^{part}_{q(t)}x(t)
\end{split}
\right.
\end{equation}
where the system matrices $A^{part}_{q},B^{part}_{q},C_{q}^{part}$,
 $q \in \{q_1,q_2\}=Q$ are as follows
\[
 \begin{split}
 & A^{part}_{q_1}=\begin{bmatrix}
                 0 &  0 \\
                 0 &  0
                \end{bmatrix},
 B_{q_1}^{part} = \begin{bmatrix}  -1.5 \\ 0 \end{bmatrix},
 C_{q_1}^{part} = \begin{bmatrix} 0 \\ 0.67 \end{bmatrix}^{T} \\
 & A^{part}_{q_2}=\begin{bmatrix}
            3 &  0 \\
            0  & 2
           \end{bmatrix}, 
 B^{part}_{q_2}=\begin{bmatrix} 0 \\ 1.5 \end{bmatrix},
 C^{part}_{q_2}=\begin{bmatrix} -0.67 &  0 \end{bmatrix}
\end{split}
\]
The \LSS\ $\Sigma_{part}$ is a \emph{partial realization} of the input-output map
$f_1$ from the initial state $x^{part}_{1}=(0,0)^{T}$, 
and 
of the input-output map $f_2$ from the initial state
$x^{part}_{2}=(-1,5,0)^{T}$.
That is, the  Markov parameters of $\Phi=\{f_1,f_2\}$ which are
 indexed by sequences from the set
$\{\epsilon,q_1,q_2\}$ coincide with those of the
input-output maps induced by the initial state $x^{part}_i$, $i=1,2$.
It is easy to see that $\Sigma_{part}$ is not a realization of the
input-output maps $f_i$ from the respective initial states $x^{part}_{i}$, $=1,2$.
That is, the input-output maps induced by the respective states $x^{part}_{i}$, $i=1,2$,
do not coincide with the maps $f_i$. One can either check it by direct calculation, or
by using uniqueness of a minimal realization. Using the latter approach,
it is enough to notice that $\Sigma_m$ is a minimal realization of $\Phi$ and 
it is of dimension $4$. Since $\Sigma_{part}$ of dimension $2$, hence smaller than 
the dimension of $\Sigma_m$, and all minimal realizations have to be of the
same dimension, it follows that $\Sigma_{part}$ cannot be a realization of $\Phi$.

In fact, by checking the Markov-parameters, one can see that $\Sigma_{part}$
recreates only the Markov-parameters indexed by sequences of length at most $1$, but 
there is a Markov-parameter of $\Phi$, indexed by a sequence of length $2$, 
which is not generated by $\Sigma_{part}$.
That is, it indeed happens  that a family of input-output maps generated by a 
\LSS\ satisfies the sufficient conditions for existence of a $N$-partial realization for some 
$N$, but the obtained partial realization is not a complete realization of the family of input-output 
maps.



\section{Linear switched systems}
\label{sect:switch_sys}
\label{sect:switch}
This section contains the definition of linear switched systems.
We will start with fixing notation and terminology
 which will be used throughout the paper.
 The notation used in this paper is mostly the standard one used in the field
 of control theory and formal language theory. In order to make the task of
 the reader easier, below we will list the most important notational
 conventions, grouped according to the disciplines.
 
 \subsubsection{ Notation from general mathematics and control theory}
\label{sect:prelim:not1}
 Denote by $T$ the set $[0, +\infty) \subseteq \mathbb{R}$ of all non-negative
 reals. The set $T$ will be the time-axis of the systems discussed in this paper.
 For any $m \ge 0$, denote by
 $PC(T,\mathbb{R}^{m})$ the class of
 piecewise-continuous maps from $T$ to $\mathbb{R}^{m}$. That is, 
 $f \in PC(T,\mathbb{R}^{m})$, if $f$ has finitely many points of discontinuity on
 each finite interval $[0,t]$, $t \in T$, and at each point of discontinuity the
 right- and left-hand side limits exist and they are finite.
Denote by $\mathbb{N}$ the set of natural number including $0$.
 By abuse of notation we will denote any constant function
 $f:T \rightarrow \mathbb{R}^{m}$ by its value. That is, if
 $f(t)=a \in \mathbb{R}^{m}$ for all $t \in T$, then $f$ will be
 denoted by $a$. 
 For any function $g$ the range of $g$ will be denoted by
 $\IM g$, i.e. if $g:A \rightarrow B$ for some sets $A$ and $B$, then
 $\IM g=\{ g(a) \in B \mid a \in A\}$.
  If $\mathcal{X}$ is a vector space and $Z$ is a subset of $\mathcal{X}$, then
  $\SPAN Z$ denotes the linear span of elements of $Z$ in $\mathcal{X}$.
   If $\mathcal{X},\mathcal{Y},\mathcal{Z}$ are vector spaces over $\mathbb{R}$, and 
   $F_{1}: \mathcal{X} \rightarrow \mathcal{Y}$,  $F_{2}: \mathcal{Y} \rightarrow \mathcal{Z}$ are linear maps,  then $F_{1}F_{2}$ denotes the composition 
   $F_{1}$ and $F_{2}$.

  Let $\phi: \mathbb{R}^{k} \rightarrow \mathbb{R}^{p \times m}$ be a smooth map.
  Consider a $k$ tuple of natural numbers
  $\alpha=(\alpha_{1},\alpha_{2},\ldots, \alpha_{k}) \in \mathbb{N}^{k}$, where
  $\alpha_1,\ldots,\alpha_k \in \mathbb{N}$.
   We will denote $D^{\alpha} \phi$ the partial derivative
   of $\phi(t_1,t_2,\ldots,t_k)$ evaluated at zero, such that
  the  order of derivation with respect to the variable $t_i$ is 
  $\alpha_{i}$ for $i=1,\ldots,k$. That is, 
   \[ 
    D^{\alpha} \phi= 
     \frac{d^{\alpha_{1}}}{dt_{1}^{\alpha_{1}}}\frac{d^{\alpha_{2}}}{dt_{2}^{\alpha_{2}}} \cdots \frac{d^{\alpha_{k}}}{dt_{k}^{\alpha_{k}}}
    \phi(t_{1},t_{2},\ldots, t_{k})|_{t_{1}=t_{2}=\cdots =t_{k}=0}.
   \] 

  For each $i=1,2,\ldots,n$, $e_j$ denotes the $j$th unit vector of $\mathbb{R}^{n}$, i.e.
  $e_j=(\delta_{1,j},\delta_{2,j},\ldots, \delta_{n,j})$, where 
  $\delta_{i,j}$ is the Kronecker symbol. 

\subsubsection{Infinite matrices}
\label{sect:prelim:not_matrix}
  In this paper we will use the notation of \cite{JacobAlg1} for matrices indexed by sets
  other than natural numbers. Let $I$ and $J$ be two arbitrary sets. A (real) matrix
  $M$ whose columns are indexed by the elements of $J$ and whose rows are indexed
  by the elements of $I$ is simply a map $M:I \times J \rightarrow \mathbb{R}$.
  The set of all such matrices is denoted by $\mathbb{R}^{I \times J}$.
  The entry of $M$ indexed by the row index $i \in I$ and column index $j \in J$ is 
  denoted by $M_{i,j}$ and it is defined as
  the value of $M$ at $(i,j)$, i.e. $M_{i,j}=M(i,j)$.
  The case of usual finite matrices can be recovered by viewing $n \times m$
  real matrices as matrices from $\mathbb{R}^{\{1,2,\ldots,n\} \times \{1,2,\ldots,m\}}$.
  In the sequel, when referring to the index set of a matrix, we will identify
  any natural number $n$ with the set $\{1,2,\ldots,n\}$. In other words, 
  $\mathbb{R}^{I \times n}$ denotes the set of matrices $\mathbb{R}^{I \times \{1,2,\ldots,n\}}$
  and $\mathbb{R}^{n \times J}$ denotes the set of matrices 
  $\mathbb{R}^{\{1,2,\ldots,n\} \times J}$.

  For a matrix $M \in \mathbb{R}^{I \times J}$,
  the columns of $M$ are simply maps of the form $I \rightarrow \mathbb{R}$ and
  the rows of $M$ are maps of the form $J \rightarrow \mathbb{R}$. 
  The set of maps of the form $I \rightarrow \mathbb{R}$ and
  $J \rightarrow \mathbb{R}$ will sometimes be denoted by $\mathbb{R}^{I}$ and
  $\mathbb{R}^{J}$ respectively.
   Furthermore, if $g \in \mathbb{R}^{I}$ (resp. $g \in \mathbb{R}^{J}$)
   then the value of $g$ at $i \in I$ (resp. $j \in J$) will be denoted by $g_i$ (resp. $g_j$).
  The column of $M$ indexed by $j \in J$ will
  be denoted by $M_{.,j}$ and is defined as $M_{.,j}(i)=M_{i,j}$, $i \in I$. Similarly,
  the row if $M$ indexed by $i \in I$ will be denoted by $M_{i,.}$ and is defined
  as $M_{i,.}(j)=M_{i,j}$ for all $j \in J$.
  If $M \in \mathbb{R}^{I \times J}$ and $S \in \mathbb{R}^{J \times K}$ and
  $J$ is finite, then the product of $M$ and $S$ is the matrix $MS \in \mathbb{R}^{I \times K}$
  such that $(MS)_{i,k}=\sum_{j \in J} M_{i,j}S_{j,k}$ for all $i \in I,k \in K$. In particular,
  if $M \in \mathbb{R}^{I \times r}$ and $S \in \mathbb{R}^{r \times K}$ for some
  natural number $r \in \mathbb{N}$, then their product $MS$ is well-defined and
  it belongs to $\mathbb{R}^{I \times K}$. 

  We will identify a map $f \in \mathbb{R}^{J}$ with the matrix
  $f \in \mathbb{R}^{J \times 1}$ defined as $f_{j,1}=f_j$ for all $j \in J$.
  Hence,  for a matrix $M \in \mathbb{R}^{I \times J}$,
   the product $Mf$ is defined as the following matrix in $\mathbb{R}^{I \times 1}$;
  $(Mf)_{i,1}=\sum_{j \in J} M_{i,j}f_j$.
  In addition, we will occasionally identify the rows of a matrix
  $M \in \mathbb{R}^{I \times J}$  with matrices $\mathbb{R}^{1 \times J}$. That is, 
  the row $M_{i,.}$ of $M$ indexed by $i \in  I$ will be viewed as the matrix
  $M_{i,.}:\{1\} \times J \ni (1,j) \mapsto M_{i,j}$.  
  With this identification, the product of the row $M_{i,.}$ with $f$ is a scalar
   $M_{i,.}f=\sum_{j \in J} M_{i,j}f_j \in \mathbb{R}$. Notice that here we tacitly
  assumed that $J$ is finite.
  
  Notice that the set of all maps $\mathbb{R}^{I}$ forms a vector space
  with respect to point-wise addition and multiplication by scalar. That is, if
  $f,g \in \mathbb{R}^{I}$ and $\alpha,\beta \in \mathbb{R}$, then
  the linear combination $\alpha f + \beta g: I \rightarrow \mathbb{R}$ is defined by 
  $(\alpha f+\beta g)(i)=\alpha f(i)+\beta g(i)$ for all $i \in I$.
  Consider a matrix $M \in \mathbb{R}^{I \times J}$ and recall that its columns are
  simply elements of $\mathbb{R}^{I}$. Hence,
  it makes sense to speak of the linear subspace spanned by the columns
  of a matrix $M \in \mathbb{R}^{I \times J}$. In the sequel, the \emph{rank of $M$},
  denoted by $\Rank M \in \mathbb{N} \cup \{\infty\}$ will mean the dimension of
  the linear space spanned by the columns of $M$. If this dimension is not finite,
  then the rank is taken to be $\infty$.
  We will denote by $\IM M$ the linear space spanned by the columns of $M$.

  If $M \in \mathbb{R}^{I \times J}$ and $J$ is finite, then $M$ can be viewed as a 
  linear map from $\mathbb{R}^{J}$ to $\mathbb{R}^{I}$, defined by 
  $(Mf)(i)=\sum_{j \in J} M_{i,j}f_j = M_{i,.}f$, $i \in I$,  
  for each $f \in \mathbb{R}^{J}$. 
  If $S \in \mathbb{R}^{J \times K}$ and $K$
  is finite, then the product $MS \in \mathbb{R}^{I \times K}$ corresponds
  to the linear map $\mathbb{R}^{K} \rightarrow \mathbb{R}^{I}$ obtained by
  composing the linear map corresponding to $S$ with the linear map corresponding
  to $M$.

 \subsubsection{ Notation from the theory of formal languages }
 \label{sect:prelim:lang}
 The notation described below is standard in formal languages and automata
 theory, see \cite{GecsPeak,AutoEilen}.
 Consider a finite set $X$ which will be called the \emph{alphabet}.
 Denote by $X^{*}$ the set of finite
 sequences of elements of $X$.  Finite sequences of elements of
 $X$ will be referred to as \emph{strings} or \emph{words} over
 the alphabet $X$.
 For a word $w=a_{1}a_{2} \cdots a_{k} \in X^{*}$, $a_1,a_2,\ldots,a_k \in X$, $k > 0$ the length
 of $w$ is denoted by $|w|$, i.e. $|w|=k$. 
 We will denote by $\epsilon$ the \emph{empty sequence (word)}. 
 The length of the empty sequence $\epsilon$ is zero: $|\epsilon|=0$. 
 We will denote by $X^{+}$ the set of
 of non-empty words over $X$. That is,
 $X^{+}=X^{*}\setminus \{\epsilon\}$.
 Consider two words $w \in X^{*}$ and $v \in X^{*}$ of the form
 $v=v_{1}v_2 \cdots v_{k}$, and  $w=w_{1}w_2 \cdots w_{m}$,
 $v_1,v_2,\ldots,v_k,w_1,w_2,\ldots,w_{m} \in X$.
 Define the concatenation $vw \in X^{*}$ of the words $v$ and $w$ as the 
 the word $vw=v_{1}v_2 \cdots v_{k}w_{1}w_2 \cdots w_{m}$. In particular, if
 $v=\epsilon$, i.e. if $k=0$, then $vw=w$. Similarly, if $w=\epsilon$, i.e. $m=0$, 
 then $vw=v$.
 If $w \in X^{+}$ is a word, then
 $w^{k}$ denotes the word $\underbrace{ww \cdots w}_{k-times}$.
 The word $w^{0}$ is just the empty word $\epsilon$.

\subsection{ Definition and basic properties of linear switched systems }
Below we present the formal definition of \SLSS. For a more detailed exposition, 
see \cite{Sun:Book,D:Lib,MP:RealForm,MP:Phd,MP:BigArticle}.
\begin{Definition}[Linear switched systems]
\label{lin_switch:def}
A \emph{linear switched system} (abbreviated by \LSS) is a control system $\Sigma$
 of the form
\begin{equation}
\label{lin_switch0}
\Sigma: \left \{ \begin{split}
   \dot x(t)=& A_{q(t)}x(t)+B_{q(t)}u(t)  \\
   y(t)=&C_{q(t)}x(t)
\end{split}
\right.
\end{equation}
\begin{itemize}
\item
 $x(t) \in \mathcal{X}$ is the continuous state at time $t \in T$,
\item
 $u(t) \in \mathcal{U}$ 
 denotes the continuous input at time $t \in T$, 
\item 
$q(t) \in Q$ denotes
 the discrete mode (state) at time $t$, 
\item
 $y(t) \in \mathbb{R}^{p}$ denotes the continuous output at time $t \in T$.  
\item
 The state-space is $\mathcal{X}=\mathbb{R}^{n}$, the
input-space is $\mathcal{U}=\mathbb{R}^{m}$, the output-space
is $\mathcal{Y}=\mathbb{R}^{p}$, and
$Q$ is the finite set of discrete modes (discrete states).
 Here $n,m,p$ are positive integers.
\item
For each discrete mode $q \in Q$,
the corresponding matrices are of the form
$A_{q} \in \mathbb{R}^{n \times n}$, 
$B_{q} \in \mathbb{R}^{n \times m}$ and $C_{q} \in \mathbb{R}^{p \times n}$. 
\end{itemize}  
    We will use 
    \( \SwitchSysLin \) as a
short-hand notation for \SLSS\  of the form (\ref{lin_switch0}).
\end{Definition}
\begin{Notation}[Notation for the spaces of inputs and outputs]
\label{switch_sys:not1}
In the sequel we denote by
 $\mathcal{U}$ the space $\mathbb{R}^{m}$ of continuous-valued input,
 by $\mathcal{Y}$ the space $\mathbb{R}^{p}$ of continuous-valued outputs, and
by $Q$ the set of discrete modes. 
\end{Notation}
Informally, the system (\ref{lin_switch0}) evolves as follows. For a 
piecewise-constant \emph{switching signal} $q(.): T \rightarrow Q$
and piecewise-continuous input $u: T \rightarrow \mathcal{U}$,
the \emph{state trajectory}
$x: T \rightarrow \mathcal{U}$ is a continuous piecewise-differentiable
function which satisfies the differential equation
(\ref{lin_switch0}). The output at time $t \in T$
is obtained by applying to $x(t)$ the readout map (matrix) 
$C_{q(t)}$.  
Below we define state- and output- trajectories more rigorously.
To this end, we define the notion of switching sequences. 
\begin{Definition}[Switching sequences]
A \emph{switching sequence} is a sequence of
the form $w=(q_1,t_1)(q_2,t_2)\cdots (q_k,t_k)$, where
$q_1,q_2,\ldots, q_k \in Q$ are discrete modes and
$t_1,t_2,\ldots, t_k$ denote the \emph{switching times} and $k \ge 0$.
The set of all switching  sequences is denoted by $(Q \times T)^{*}$. If
$k=0$ then we say that $w$ is the empty switching sequence and we denote it by $\epsilon$.
We denote the set of all \emph{non-empty switching sequences} by $(Q \times T)^{+}$.
\end{Definition}
The interpretation
of the sequence $w=(q_1,t_1)(q_2,t_2) \cdots (q_k,t_k)$ is the following. From time instance 
$0$ to time instance $t_1$
the active discrete mode is $q_1$, i.e. the value of the 
switching signal is $q_1$, from $t_1$ to $t_1+t_2$ the value
of the switching signal is $q_2$, from $t_1+t_2$ to
$t_1+t_2+t_3$ the value of the switching signal is $q_3$, and so on.
That is, the non-negative real $t_i$ indicates the time spent in the discrete mode
$q_i$, for all $i=1,2,\ldots,k$.
\emph{In this paper
the switching sequences are regarded as inputs
and we allow any switching sequence to occur.}
\begin{Definition}[State of \SLSS]
\label{def:state_traject}
Let $u \in PC(T,\mathcal{U})$ be a continuous-valued input and let
$w=(q_{1},t_{1})(q_{2},t_{2}) \cdots (q_{k},t_{k}) \in (Q \times T)^{+}$ be a
 non-empty switching sequence.
The \emph{state of $\Sigma$ reached from the initial state $x_0 \in \mathcal{X}$ with the
inputs $u$ and $w$} is denoted by $x_{\Sigma}(x_0,u,w)$ and 
it is defined recursively on $k$ as follows.
\begin{itemize}
\item
If $k=1$, then 
    $x_{\Sigma}(x_0,u,(q_1,t_1))$ is the solution at time $t_1$
    of the differential equation
    \[ \dot x(s)=A_{q_1}x(s)+B_{q_1}u(s) \]
    with the initial condition $x(0)=x_0$.
\item
    If $x_{\Sigma}(x_0,u,(q_1,t_1)(q_2,t_2) \cdots (q_{k-1},t_{k-1}))$ is already
    defined,
    then let $x_{\Sigma}(x_0,u,w)$ be the solution at time $t_k$ 
    of the differential equation
    \[ \dot x(s)=A_{q_k}x(s)+B_{q_k}u(s+\sum_{j=1}^{k-1} t_j) \]
    with the initial condition 
    $x(0)=x_{\Sigma}(x_0,u,(q_1,t_1)\cdots (q_{k-1},t_{k-1}))$.
\end{itemize}
\end{Definition}
    That is, the states evolves according to the differential
    equation determined by the discrete mode. If a switch
    occurs, then the state at the time of the switch is taken
    as the initial condition for generating a solution to
    the differential equation associated with the new discrete
    mode.
     As a next step, we will define the output trajectories of
     \SLSS. 
     \begin{Definition}[Output of \SLSS]
     Consider a continuous-valued input $u \in PC(T,\mathcal{U})$ and 
     a non-empty switching sequence 
     $w=(q_1,t_1)(q_2,t_2)\cdots (q_k,t_k) \in (Q \times T)^{+}$.
     The \emph{output 
     generated by the \LSS\ $\Sigma$ if it is started from
     initial state $x_0 \in \mathcal{X}$ and fed with the inputs $u$ and $w$},
     denoted by
     $y_{\Sigma}(x_0,u,w) \in \mathcal{Y}$, is defined by
     \begin{equation}
     \label{sect:switch:iodef1}
        y_{\Sigma}(x_0,u,w)=C_{q_{k}}x_{\Sigma}(x_0,u,w)
     \end{equation}	
    \end{Definition}
    That is, the current output is obtained from the current continuous
    state $x_{\Sigma}(x_0,u,w)$ by the application of the
    readout map $C_{q_{k}}$ associated with the current
    discrete mode $q_k$.
    We define the
    input-output map of a \LSS\  induced by a 
    particular initial state as follows.
    \begin{Definition}[Input-output maps of \SLSS]
    Consider a state $x_0 \in \mathcal{X}$ of $\Sigma$.
    Define the \emph{input-output map of $\Sigma$ induced by the state $x_0$} as the map
    \( y_{\Sigma}(x_0,.,.):PC(T,\mathcal{U}) \times (Q \times T)^{+}
     \rightarrow \mathcal{Y}
    \) such that for all continuous-valued input $u \in PC(T,\mathcal{U})$ and
     for all non-empty switching sequence $w \in (Q \times T)^{+}$
    \begin{equation}
    \label{sect:swicth:iodef2}
       y_{\Sigma}(x_0,.,.)(u,w)=y_{\Sigma}(x_0,u,w)
    \end{equation}   
    That is, the value of $y_{\Sigma}(x_0,.,.)$ at $(u,w)$ equals the output 
    generated by $\Sigma$ if started from the initial states $x_0$ and fed
    inputs $u$ and $w$.
    By abuse of notation we will denote
    $y_{\Sigma}(x_0,.,.)(u,w)$ by $y_{\Sigma}(x_0,u,w)$.
   \end{Definition}

\section{ Overview of realization theory for \LSS}
\label{sect:lin_sys:real_theo}
 The goal of this section is to present a brief overview of realization theory of \SLSS.
 In Subsection \ref{lin_sys:io:map} we recall the definition of
 input-output maps of \SLSS\ and the notion of a \LSS\ realization. 
In Subsection \ref{lin_sys:system_theo}.
we review a number of system-theoretic concepts such as
 observability, semi-reachability, minimality and \LSS\  morphisms. 
 In Subsection \ref{lin_sys:real_problem} the realization problem will be formulated.
 Subsection \ref{lin_sys:gen_ker} recalls the concept of generalized kernel
 representation of a family of
 input-output maps. 
 Subsection \ref{lin_sys:markov} presents  the definition of Markov-parameters
 for \SLSS. 
 Finally, Subsection \ref{lin_sys:real} presents the main results on realization
 theory of \SLSS.
 For a more details on the material of this section see \cite{MP:Phd,MP:BigArticle,MP:RealForm}.

\subsection{Input-output maps}
\label{lin_sys:io:map}
  In this section we define the class of maps
  which represent potential input-output maps of \SLSS. In addition, we introduce
  the notion of a \LSS\ realization. 
   \begin{Definition}[Input-output maps]
   \label{sect:switch_sys:def-1}
    In this paper, unless stated otherwise, an \emph{input-output map} will mean 
     a map of
    the form $f: PC(T, \mathcal{U}) \times (Q \times T)^{+}  \rightarrow \mathcal{Y}$. 
    The set of all such maps will be denoted by $\IO$. A \emph{family of
    input-output maps} 
    is just a (possibly infinite) subset of the set of all
    input-output maps $\IO$.
  \end{Definition}
   That is, an input-output map maps continuous-valued inputs and non-empty switching
   sequences to continuous-valued outputs.
    In order to
   formalize the notion of a realization by a \LSS\  of a family of
   input-output maps we will adopt the following formalism.
   \begin{Definition}[Realization of input-output maps]
   \label{switch_sys:real:def1}
   Consider a set $\Phi \subseteq  \IO$ of input-output maps.
      The family $\Phi$ is said to be \emph{realized} by 
        a \LSS\ 
   $\Sigma$  of the form (\ref{lin_switch0}) if there exists a map
   $\mu:\Phi \rightarrow \mathcal{X}$, which maps
   each input-output map $f$ from $\Phi$ to a state $
   \mu(f)$ of $\Sigma$,  such that 
    $f$ equals the input-output map 
   induced by $\mu(f)$, i.e.
   \begin{equation}
   \label{sect:problem_form:real0}
     \forall u \in PC(T,\mathcal{U}), w \in (Q \times T)^{+}: 
       y_{\Sigma}(\mu(f),u,w)=f(u,w) 
   \end{equation}  
   \end{Definition}
   One can think of the map $\mu$ as a way to determine the 
   initial state corresponding to each element of $\Phi$.
   \begin{Definition}[Realizations]
   \label{switch_sys:real:def2}
   Let $\Phi$ be a family of input-output maps.
   We will refer to
   the pairs $(\Sigma,\mu)$, where $\Sigma$ is a \LSS\  of the form
   (\ref{lin_switch0}) and $\mu:\Phi \rightarrow \mathcal{X}$
   is a map mapping elements of $\Phi$ to the states of $\Sigma$,
   as \LSS\ \emph{realizations} (\emph{realizations} for short). 
   The realization $(\Sigma,\mu)$
   is  said to be a \emph{realization of $\Phi$}, if 
   (\ref{sect:problem_form:real0}) holds for all 
   $f \in \Phi$. 
  \end{Definition}
   Note that not any
   realization $(\Sigma,\mu)$ with 
   $\mu: \Phi \rightarrow \mathcal{X}$ is a realization of $\Phi$. The statement that
   $(\Sigma,\mu)$ is a \emph{realization} only expresses the fact that we associate a 
   state of $\Sigma$ with each element of $\Phi$. However, we do not yet require that
   the input-output map induced by a designated state equals the corresponding element of $\Phi$.
   The latter is required only if we claim that $(\Sigma,\mu)$ \emph{is a realization of $\Phi$}.

 \subsection{System-theoretic concepts }
 \label{lin_sys:system_theo}
   The goal of this section is to define system theoretic concepts such as
   observability, span-reachability, system morphism, dimension and 
   minimality for \SLSS\  and for \LSS\ realizations. 
   Throughout the section $\Sigma$ denotes a \LSS\ of the form (\ref{lin_switch0}).

   The \emph{reachable set} of $\Sigma$ 
   from a set of initial states $\mathcal{X}_{0} \subseteq \mathcal{X}$ is 
   defined as
  \begin{equation*}
  \label{sect:problem_form:reach0}
    Reach(\Sigma, \mathcal{X}_{0}) = \{ x_{\Sigma}(x_{0},u,w) \in \mathcal{X} 
   \mid  u \in PC(T,\mathcal{U}), w \in (Q \times T)^{+}, x_{0} \in \mathcal{X}_{0} \} 
   \end{equation*}
   That is, $Reach(\Sigma,\mathcal{X}_{0})$ is the set of all those states
   which are obtained by starting the system from an initial state in $\mathcal{X}_{0}$,
   applying some
   continuous-valued input and some finite switching sequence, and considering the state 
   at the last switching time.
  \begin{Definition}[(Semi-)Reachability]  
    The \LSS\  $\Sigma$ is said to be \emph{reachable} from $\mathcal{X}_{0}$ if
    $Reach(\Sigma,\mathcal{X}_{0})=\mathcal{X}$ holds.
    The \LSS\ $\Sigma$ is \emph{\WR}  from $\mathcal{X}_{0}$ if 
    $\mathcal{X}$ is the smallest vector space containing 
    $Reach(\Sigma,\mathcal{X}_{0})$.
  \end{Definition}
   In other words, 
    $\Sigma$ is \WR\  from $\mathcal{X}_{0}$ if 
    the linear span of the elements of the reachable set 
    $Reach(\Sigma,\mathcal{X}_{0})$ yields the whole state-space
    $\mathcal{X}$.
   We proceed with defining the notion of observability for \SLSS.
   \begin{Definition}[Observability and Indistinguishability]
   Two states
   $x_{1}\neq x_{2} \in \mathcal{X}$ of
   the \LSS\  $\Sigma$
   are \emph{indistinguishable} if
   the input-output maps induced by $x_1$ and $x_2$ coincide,
   i.e. for all $u \in PC(T,\mathcal{U})$ and $w \in (Q \times T)^{+}$,
  $y_{\Sigma}(x_1,u,w)=y_{\Sigma}(x_2,u,w)$.
    The \LSS\  $\Sigma$ is called \emph{observable} if
   it has no pair of distinct indistinguishable states.
   \end{Definition}
   In other words, observability means that if we pick any two states of the system, then
   we are able to distinguish between them by feeding a 
  suitable continuous-valued input and
   a suitable switching sequence and then observing the resulting output.
   Below we will define the notion of dimension for \SLSS.
   \begin{Definition}[Dimension of \SLSS]
   \label{switch_sys:dim:def}
    Define the dimension of $\Sigma$, denoted by $\dim \Sigma$, as 
    the dimension $\dim \mathcal{X}=n$ of is state-space.
   \end{Definition}
Now we are ready to define the concept of minimality for \SLSS\ .
\begin{Definition}[Minimality of \SLSS\ ]
 Let $\Phi$ be a family of input-output maps and let $(\Sigma,\mu)$ be a 
 \LSS\  realization
 of $\Phi$.
 $(\Sigma,\mu)$ is \emph{a minimal realization of $\Phi$}, if
 for any \LSS\  realization
 $(\hat{\Sigma},\hat{\mu})$ of $\Phi$, $\dim \Sigma \le \dim \hat{\Sigma}$.
\end{Definition}
In simple words, a \LSS\  realization is a minimal realization of $\Phi$
if it has the smallest dimensional state-space among all the linear switched
systems which are realizations of $\Phi$.

   The notions of observability and semi-reachability
    can be extended to \LSS\  realizations as follows.
   \begin{Definition}[Observability and semi-reachability of realizations]
   \label{switch_sys:real:def3}
    Let $\Phi$ be a family
    of input-output maps and let $\mu:\Phi \rightarrow \mathcal{X}$ be a 
    map from $\Phi$ to the state-space of $\Sigma$.
    The realization
    $(\Sigma,\mu)$ is \emph{\WR\ }, 
    if $\Sigma$ is \WR\ from the range $\IM \mu$ of
    $\mu$. The realization 
    $(\Sigma,\mu)$ is \emph{observable}, if $\Sigma$ is observable.
   \end{Definition}
   As the next step,  we will define the notion of a \LSS\  morphism.  
  \begin{Definition}[Linear switched system morphism] 
  \label{sect:problem_form:lin:morphism}
   Consider a \LSS\  $\Sigma_1$ of the form (\ref{lin_switch0}) and
   a \LSS\  $\Sigma_2$ of the form 
   $\Sigma_{2}=\SwitchSysLin[a]$\footnote{Notice that the \SLSS\ $\Sigma_1$ and $\Sigma_2$ have the same set of discrete modes.}.
    A linear map
    $S:\mathcal{X} \rightarrow \mathcal{X}_{a}$
    is said to be a \emph{\LSS\  morphism}
    from $\Sigma_{1}$ to $\Sigma_{2}$,  and it is denoted by
    $S:\Sigma_{1} \rightarrow \Sigma_{2}$, if for all
    discrete modes $q \in Q$, 
   \begin{equation}
   \label{sect:problem_form:lin:morphism:eq0}
   \begin{array}{llll}
    A^{a}_{q}S=SA_{q}, & B_{q}^{a}=SB_{q}, &
    C_{q}^{a}S=C_{q} 
   \end{array}
  \end{equation}
  The map $S$ is called surjective ( injective ) if
  it is surjective ( injective ) as a linear map. The map 
  $S$ is said to be a \LSS\  isomorphism,
  if it is an isomorphism as a linear map.
  Consider two \LSS\  realizations $(\Sigma_{1},\mu_{1})$
  and $(\Sigma_2,\mu_2)$ such that the domain of definition of both $\mu_1$ and
  $\mu_2$ is a certain family $\Phi$ of input-output maps.
  A \LSS\  morphism 
  $S:\Sigma_{1} \rightarrow \Sigma_{2}$ 
  is  called a \emph{\LSS\  morphism } from
  realization $(\Sigma_{1},\mu_{1})$ to $(\Sigma_{2},\mu_{2})$, 
  if 
  $S \circ \mu_{1}=\mu_{2}$ holds, or, in other words, if for all $f \in \Phi$,
  $S(\mu_1(f))=\mu_2(f)$.
  The fact that $S$ is a \LSS\  morphism from $(\Sigma_1,\mu_1)$
  to $(\Sigma_2,\mu_2)$ will be denoted by
  $S:(\Sigma_{1},\mu_{1}) \rightarrow (\Sigma_{2},\mu_{2})$.
  The \SLSS\  realizations $(\Sigma_{1},\mu_{1})$ and 
  $(\Sigma_{2},\mu_{2})$ are
  said to be \emph{algebraically similar} or \emph{isomorphic}
  if there exists an \LSS\  isomorphism
  $S:(\Sigma_{1},\mu_{1}) \rightarrow (\Sigma_{2},\mu_{2})$. 
  \end{Definition}

\subsection{ Realization problem}
\label{lin_sys:real_problem}
The realization problem for \SLSS\  can be formulated as follows.
\begin{Problem}[Realization problem for \SLSS\ ]
\label{lin_switch:real_problem}
 Find necessary and sufficient conditions for existence of a \LSS\ 
 realization for a family of input-output maps $\Phi$.
 Find a characterization of minimal \LSS\  realizations of $\Phi$.
 Determine if minimal realizations of $\Phi$ are unique in any sense.
\end{Problem}

\subsection{ Generalized kernel representation}
\label{lin_sys:gen_ker}
 In this section we will recall the notion of
\emph{generalized kernel representation}.
It turns out that a family of input-output
maps can be realized by a \LSS\  only if
it admits a generalized kernel representation.
Informally, a family $\Phi$ of input-output maps has a 
generalized kernel representation if the following hold.
\begin{enumerate}
\item
  There exists an input-output map $y^{\Phi}$ such
  that for all $f \in \Phi$, 
  $f(u,w)=f(0,w)+y^{\Phi}(u,w)$ for all continuous-valued inputs $u$ and
  switching sequences $w$.
\item
   Each element $f$ of $\Phi$ is affine in continuous-valued
   inputs and analytic in switching times for all constant
   inputs.
\end{enumerate}
A good intuition for the notion of generalized kernel
representation can be derived by analogy with input-output maps
of linear systems. Recall from \cite{Cal:Des} that an input-output map
$y: PC(T,\mathcal{U}) \times T \rightarrow \mathcal{Y}$ can be realized
by a linear system $(A,B,C) \in \mathbb{R}^{n \times n} \times \mathbb{R}^{n \times m} \times \mathbb{R}^{p \times n}$ from the initial state 
$x_0 \in \mathbb{R}^{n}$, only if there exist analytic functions $K:T \rightarrow \mathbb{R}^{p}$
and $G:T \rightarrow \mathbb{R}^{p \times m}$ such that
\begin{equation} 
\label{gen_ker_int/0}
  y(u,t)=K(t)+\int_{0}^{t} G(t-s)u(s)ds 
\end{equation}
More precisely, in this case 
\begin{equation}
\label{gen_ker_int/1}
K(t)=Ce^{At}x_0 \mbox{ \ \ and \ \ } G(t)=Ce^{At}B
\end{equation}
Using the setting above, the maps $K_{w}^{f,\Phi}$ and $G_{w}^{\Phi}$ to be defined below are analogous to the map $K$ and $G$ respectively.
The formal definition is as follows.
\begin{Definition}[Generalized kernel-representation]
\label{sect:io:def1}
   Consider a family of input-output maps $\Phi$, as defined in 
   Definition \ref{sect:switch_sys:def-1}.
   The family $\Phi$ is said to have 
 \emph{generalized kernel representation},  if
   for all input-output maps $f \in \Phi$ and for all
   non-empty sequences of discrete modes
   $w=q_{1}q_{2} \cdots q_{k} \in Q^{+}$, 
   $q_{1},q_2\ldots, q_{k} \in Q$, $k > 0$, 
   there
  exist functions
  \[ K_{w}^{f,\Phi}:T^{k} \rightarrow \mathbb{R}^{p} \
 \mbox{ and }
  G_{w}^{\Phi}:T^{k} \rightarrow \mathbb{R}^{p \times m} \]
 such that the following holds.
 \begin{enumerate}
 \item
 \label{gen_ker/1}
     For each word $w \in Q^{+}$ and
     for each input-output map $f \in \Phi$,
     the functions
     $K_{w}^{f, \Phi}$ and $G^{\Phi}_{w}$ are analytic.
 \item
 \label{gen_ker/2}
   For each input-output map $f \in \Phi$ and for each (possibly empty) 
   sequences $w,v\in Q^{*}$, and for each discrete mode $q \in Q$, 
 it holds that for all
   $t_1,t_2,\ldots,t_{|w|}, t,\hat{t},t_{|w|+2},\ldots, t_{|w|+|v|+1} \in T$,
   \begin{equation*}
   \begin{split}
     K_{wqqv}^{f,\Phi}(t_{1},t_{2}, \ldots, t_{|w|},t,\hat{t},t_{|w|+2}, 
    \ldots
    t_{|w|+|v|+1})= \\
    =K^{f,\Phi}_{wqv}(t_{1},t_{2}, \ldots t_{|w|},t+\hat{t},t_{|w|+2} \ldots t_{|w|+|v|+1}) \\
     G_{wqqv}^{\Phi}(t_{1},t_{2}, \ldots, t_{|w|},t,\hat{t},t_{|w|+2}, 
    \ldots
    t_{|w|+|v|+1})=\\
    =G^{\Phi}_{wqv}(t_{1},t_{2}, \ldots t_{|w|},t+\hat{t},t_{|w|+2} \ldots t_{|w|+|v|+1})
 \end{split}
 \end{equation*}
  
 \item
 \label{gen_ker/3}
   For each pair of sequences $v,w \in Q^{*}$
   such that and $w$ is not the empty word, i.e.  $|w| > 0$,
    and for each
   discrete mode $q \in Q$,
   the following holds. For each input-output map $f \in \Phi$ and
   for each $t_1,t_2,\ldots, t_{|vw|} \in T$,
    \[  K^{f,\Phi}_{vqw}(t_{1},t_2,\ldots,t_{|v|}, 0, t_{|v|+1}, \ldots,t_{|wv|})=
      K^{f,\Phi}_{vw}(t_{1},t_{2},\ldots,t_{|vw|})
  \]    
  For each pair of words $v,w \in Q^{*}$
  such that both $v$ and $w$ are not empty, i.e.  $|v| > 0$, $|w| > 0$, 
  and for each
  discrete mode $q \in Q$, 
  the following holds.
  For each $t_1,t_2,\ldots, t_{|vw|} \in T$, 
  \[
      G_{vqw}^{\Phi}(t_{1},t_2,\ldots, t_{|v|}, 0, t_{|v|+1}, \ldots,t_{|wv|})=
      G_{vw}^{\Phi}(t_{1},t_2,\ldots ,t_{|vw|})
  \]    
 \item
 \label{gen_ker/4}
   For each input-output map $f \in \Phi$, for
   each non-empty switching sequence 
   $w=(q_{1},t_{1})(q_{2},t_{2}) \cdots (q_{k},t_{k}) (Q \times T)^{+}$,
   where $q_1,q_2,\ldots,q_k \in Q$  and $t_1,t_2,\ldots,t_k \in T$,
   each piecewise-continuous input $u \in PC(T,\mathcal{U})$, the following holds.
  \[
   \begin{split}
   & f(u,w)=
    K_{q_{1}q_{2} \cdots q_{k}}^{f,\Phi}(t_{1},t_{2},\ldots,t_{k})+ \\
  \sum_{i=1}^{k} \int_{0}^{t_{i}}
    &  G^{\Phi}_{q_{i}q_{i+1}
     \cdots q_{k}}(t_{i}-s,t_{i+1},\ldots,t_{k})
     u(s+ \sum_{j=1}^{i-1} t_{j})ds
   \end{split}
 \]
\end{enumerate}
\end{Definition}
The reader  may view the functions $K^{f,\Phi}_{w}$ as the part of the 
output which
depends on the initial condition and the functions $G^{\Phi}_{w}$ as
functions determining the dependence of the output on the continuous
inputs.

In fact, let $(\Sigma,\mu)$ be a \LSS\ realization of $\Phi$ and
assume that $\Sigma$ is of the form (\ref{lin_switch0}). It is easy to see
that for all $f \in \Phi$ and 
for any input $u \in PC(T,\mathcal{U})$ and
switching sequence $(q_1,t_1)(q_2,t_2) \cdots (q_k,t_k) \in (Q \times T)^{+}$,
 \begin{equation} 
   \label{sect:lin_switch:state:eq2}
   \begin{split}
      & y_{\Sigma}(\mu(f),u,w) = 
     C_{q_{k}}e^{A_{q_{k}}t_{k}}e^{A_{q_{k-1}}t_{k-1}} \cdots
      e^{A_{q_{1}}t_{1}}\mu(f)+  \\
      & +\int_{0}^{t_{k}} C_{q_{k}}e^{A_{q_{k}}(t_{k}-s)} 
      B_{q_{k}}
       u(s+\sum_{1}^{k-1} t_{i})ds
      +  \\
      & + C_{q_{k}}e^{A_{q_{k}}t_{k}} 
      \int_{0}^{t_{k-1}}
      e^{A_{q_{k-1}}(t_{k-1}-s)} 
      B_{q_{k-1}}u(s+\sum_{1}^{k-2} 
      t_{i})ds + \cdots \\
       & \cdots  
      + C_{q_{k}}e^{A_{q_{k}}t_{k}}e^{A_{q_{k-1}}t_{k-1}}
      \cdots e^{A_{q_{2}}t_{2}}
       \int_{0}^{t_{1}}
      e^{A_{q_{1}}(t_{1}-s)}B_{q_{1}}u(s)ds
      \end{split}
  \end{equation}
 From the equation above it is easy to see that $\Phi$ admits a 
 hybrid kernel representation of the form
 \begin{equation}
   \label{gen_ker_int/2}
   \begin{split}
   G^{\Phi}_{q_{1}q_{2}\cdots q_{k}}(t_{1},t_{2},\ldots,t_{k})
     = &
      C_{q_{k}}e^{A_{q_{k}}t_{k}}e^{A_{q_{k-1}}t_{k-1}} \cdots
      e^{A_{q_{1}}t_{1}}B_{q_{1}}  \\
    K^{f,\Phi}_{q_{1}q_{2} \cdots q_{k}}(t_{1},t_{2},\ldots,t_{k})= &
      C_{q_{k}}e^{A_{q_{k}}t_{k}}e^{A_{q_{k-1}}t_{k-1}} \cdots
      e^{A_{q_{1}}t_{1}}\mu(f).
   \end{split}
  \end{equation}
  for all non-empty sequences of discrete modes
  $q_1,q_2,\ldots,q_k \in Q$, $k \ge 1$.

\subsection{ Generalized Markov parameters for \LSS}
\label{lin_sys:markov}
Next we will define the notion of Markov parameters for input-output maps.
Markov-parameters play a central role in (partial) realization theory of \SLSS. 

Before proceeding further, recall from classical linear systems theory \cite{Cal:Des}
the notion of Markov parameter. 
Consider the linear input-output map 
of the form (\ref{gen_ker_int/0})
and define the Markov parameters of this map as derivatives of
the maps $K:T \rightarrow \mathbb{R}^{p}$ and $G:T \rightarrow \mathbb{R}^{p \times m}$;
$M_{k}=\frac{d^k}{dt^{k}} G(t)|_{t=0}$ and $N_{k}=\frac{d^{k}}{dt^{k}} K(t)|_{t=0}$
for all $k \ge 0$.
In turn, the derivatives of $G(t)$ and $K(t)$ can be expressed as the
derivatives of the input-output map as follows.
For a constant input $u \in \mathcal{U}$ define the map
\[ f_{u}:T \ni t \mapsto K(t)+(\int_{0}^{t} G(t-s)ds)u \]
That is, $f_{u}(t)$ is just the value of the input-output map (\ref{gen_ker_int/0}) at 
time $t$ if a constant input $u$ is fed in. 
Then $f_{0}(t)=K(t)$ and hence $N_{k}=\frac{d^k}{dt^{k}} f_{0}(t)|_{t=0}$.
Similarly, the $j$th column of $M_k$ can be written as
\( \frac{d^{k+1}}{dt^{k+1}} (f_{e_j}(t)-f_{0}(t))|_{t=0} \)
where $e_j$, $j=1,2,\ldots,m$ is the $j$th unit vector of $\mathcal{U}=\mathbb{R}^{m}$.
That is, \emph{the Markov-parameters of a linear input-output map are the
high-order time derivatives of the output trajectories induced by certain constant 
inputs, evaluated at zero}.

For \SLSS\ the Markov-parameters are defined in a similar way. Before proceeding
to the definition, we need the following notation
\begin{Notation}[Input-output maps as time functions]
\label{sect:io:io_der} 
 Consider an input-output map $f$ as in Definition \ref{sect:switch_sys:def-1}, 
 a non-empty sequence  of 
 discrete modes $w=q_1q_2\cdots q_k \in Q^{+}$,
 $q_1,q_2,\ldots,q_k \in Q$, $k \ge 1$ and
 a continuous-valued input $u \in PC(T,\mathcal{U})$.
 Define the map $f_{u,w}:T^{k} \rightarrow \mathcal{Y}$ as follows
 \begin{equation}
 \label{sect:io:io_der:eq1}
  f_{u,w}(t_1,t_2,\ldots,t_k)=f(u,(q_1,t_1)(q_2,t_2),\ldots,(q_k,t_k))
 \end{equation}
\end{Notation}
 That is, the values of $f_{u,w}$ are obtained from the values of 
 $f$ by fixing the a piecewise-continuous input $u$ and a sequence of discrete modes
 $w$ and varying the switching times only.
\begin{Definition}[Markov-parameters of $\Phi$]
\label{linswitch_intro:def:markov}
Let $\Phi$ be a family of input-output maps admitting a 
generalized kernel representation. 
 The \emph{Markov parameters} of $\Phi$ are the vectors
 $S_{q_0,q,j}(w) \in \mathbb{R}^{p}$ and $S_{f,q}(w) \in \mathbb{R}^{p}$, defined 
 for all discrete modes $q_0,q \in Q$ , all the
 sequences of discrete modes $w \in Q^{*}$ (including the empty sequence), 
all $j=1,2,\ldots,m$ and
 $f \in \Phi$ as follows. If $w=q_1q_2\cdots q_k$ for $k \ge 0$,
$q_1,q_2,\ldots,q_k \in Q$, then
\begin{equation*}
\label{main_results:lin:arb:pow1}
\begin{split}
 S_{f,q}(w)= & \frac{d}{dt_1}\frac{d}{dt_2}\cdots \frac{d}{dt_k} f_{0,q_1q_2\cdots q_kq}(t_1,t_2,\ldots,t_k,0)|_{t_1=t_2\cdots t_k=0} \\
 S_{q,q_0,j}(w) = &
           \frac{d}{dt_0}\frac{d}{dt_1}\frac{d}{dt_2}\cdots \frac{d}{dt_k} f_{e_j,q_0q_1q_2\cdots q_kq}(t_0,t_1,t_2,\ldots,t_k,0)|_{t_0=t_1=t_2\cdots t_k=0}  \\
        & - \frac{d}{dt_0}\frac{d}{dt_1}\frac{d}{dt_2} \cdots \frac{d}{dt_k}  f_{0,q_0q_1q_2\cdots q_kq}(t_0,t_1,t_2,\ldots,t_k,0)|_{t_0=t_1=t_2\cdots t_k=0}  \\
\end{split}
\end{equation*}
Here $e_j$ is the $j$th unit vector of $\mathcal{U}=\mathbb{R}^{m}$.
\end{Definition}
The vectors $S_{f,q}(w)$ correspond to the parameters $N_{k}$ for linear systems,
and $S_{q,q_0,j}(w)$ corresponds to the $j$th column of the $M_{k}$ component of
the Markov parameters.
As it was shown in \cite{MP:Phd,MP:BigArticle,MP:RealForm}, 
there is a close relationship between
the Markov parameters of $\Phi$ and products of system matrices.
\begin{Lemma}[\cite{MP:Phd,MP:RealForm,MP:BigArticle}]
\label{sect:io:prop1}
Let $\Sigma$ be of the form (\ref{lin_switch0}) and let $\mu:\Phi \rightarrow \mathcal{X}$.
The realization $(\Sigma,\mu)$ is a realization of $\Phi$,
if and only if the Markov-parameters of $\Phi$ are equal to the following
products of the system matrices. 
\begin{equation}
    \label{sect:real:prop1:eq1}
 \begin{split}
        & S_{q,q_0,j}(q_1q_2\cdots q_k)=
        C_{q}A_{q_{k}}A_{q_{k-1}} \cdots A_{q_{1}}B_{q_{0}}e_{j} \\
        & S_{f,q}(q_1q_2\cdots q_k)=
        C_{q}A_{q_{k}}A_{q_{k-1}} \cdots A_{q_{1}}\mu(f)
 \end{split}
 \end{equation}
 for each $q,q_0 \in Q$, for each 
 $q_1,q_2,\ldots,q_k \in Q$,  $k \ge 0$, $f \in \Phi$ and
 $j=1,2,\ldots,m$.
\end{Lemma}
Similarly to the linear case, the Markov-parameters of $\Phi$ can also be
represented as derivatives of the maps $K_{w}^{f,\Phi}$ and $G_{w}^{\Phi}$,
see \cite{MP:Phd,MP:BigArticle,MP:RealForm}.

\subsection{ Main results on 
   realization theory for \SLSS}
\label{lin_sys:real}
 The purpose of this section is to present formally the
 main results on realization theory of
 \SLSS\ .
 To this end, the notion of the
Hankel-matrix $H_{\Phi}$ of $\Phi$ is defined. Similarly to the linear case, the entries of
the Hankel-matrix will be formed by the Markov parameters.
\begin{Definition}[Hankel-matrix]
\label{main_result:lin:hank:arb:def}
 Assume that the set of discrete modes $Q$ has 
 $\QNUM$ elements and choose an enumeration
 of $Q$
 \begin{equation}
 \label{main_result:lin:hank:arb:def_enum}
   Q=\{\sigma_1,\sigma_2,\ldots,\sigma_{\QNUM}\}
 \end{equation}  
Recall the notation for infinite matrices presented in Section \ref{sect:prelim:not_matrix}.
Define the \emph{Hankel-matrix}
of $\Phi$ as the infinite real matrix, columns and rows of
which are indexed as follows.
The rows of $H_{\Phi}$ are
indexed by pairs $(v,i)$ where $v \in Q^{*}$ is a finite sequence
of discrete modes and $i \in \{1,2,\ldots,p\QNUM\}$. 
The columns of $H_{\Phi}$ are indexed by pairs
$(w,j)$, where $w \in Q^{*}$ is a sequence of discrete modes and 
$j$ is either an element of $\Phi$, or $j$ is a pair of 
the form $(q,z)$, where $q$ is a discrete mode and $z \in \{1,2,\ldots,m\}$.
That is, $H_{\Phi}$ is an infinite real matrix of the form
\begin{equation}
\label{hankel:lin:def0}
H_{\Phi} \in \mathbb{R}^{(Q^{*} \times I) \times (Q^{*} \times J_{\Phi})}
\end{equation}
where $I=\{1,2,\ldots, p\QNUM\}$ and $J_{\Phi}=\Phi \cup (Q \times \{1,2,\ldots,m\})$.
The entries of $H_{\Phi}$
are defined as follows. Fix sequences of discrete modes
$w,v \in Q^{*}$ and fix an element $j \in J_{\Phi}$.
For any $i$ of the form
$i=pK+r+1$ where $K=0,1,\ldots, \QNUM-1$ and $r=0,1,\ldots, p-1$,
the entry $(H_{\Phi})_{(v,i),(w,j)}$ is defined as follows
\begin{equation}
\label{hankel:lin:def0.3}
  (H_{\Phi})_{(v,i),(w,j)}=
  \left\{\begin{array}{rl}
       (S_{\sigma_{K+1},q,z}(wv))_{r+1} & \mbox{ if }
                           j=(q,z) \in Q \times \{1,\ldots,m\} \\
       (S_{f,\sigma_{K+1}}(wv))_{r+1}   & \mbox{ if }
                                     j=f \in \Phi
   \end{array}\right.				     
\end{equation}
Here $(S_{\sigma_{K+1},q,z}(wv))_{r+1}$ and
$(S_{f,\sigma_{K+1}}(wv))_{r+1}$ denote the $r+1$th element of the
vectors $S_{\sigma_{K+1},q,z}(wv) \in \mathbb{R}^{p}$ and
$S_{f,\sigma_{K+1}}(wv) \in \mathbb{R}^{p}$ respectively.
\end{Definition}
That is, $H_{\Phi}$ is constructed from certain high-order
derivatives of the input-output maps belonging to $\Phi$.
As it was noted in Section \ref{sect:prelim:not_matrix}, 
the columns of $H_{\Phi}$ belong to the vector 
space  of all maps $(Q^{*} \times \{1,2,\ldots,p\QNUM\}) \rightarrow \mathbb{R}$.
Hence, we can speak of the linear span of the columns of $H_{\Phi}$. In addition, 
according to the convention adopted in Section \ref{sect:prelim:not_matrix}, 
the \emph{rank of $H_{\Phi}$}, denoted by $\Rank H_{\Phi} \in \mathbb{N} \cup \{+\infty\}$, is
the dimension of the linear subspace spanned by the columns of
$H_{\Phi}$. 
Now we are ready to state the main theorem on the existence
of a \LSS\  realization for arbitrary switching.
  \begin{Theorem}[Realization of input-output maps, \cite{MP:Phd,MP:RealForm,MP:BigArticle}]
  \label{sect:real:theo2}
   Let $\Phi$ be a family of input-output maps.
      Then $\Phi$ has a realization by  a \LSS\  if and
      only if $\Phi$ has a generalized kernel representation
      and 
     the rank of the associated Hankel-matrix $H_{\Phi}$ of $\Phi$  is finite, i.e, 
     $\Rank H_{\Phi} < +\infty$.
  \end{Theorem}
 Next, we state state the main result of on
 minimality of \SLSS\ .
\begin{Theorem}[Minimality, \cite{MP:Phd,MP:RealForm,MP:BigArticle}]
\label{sect:real:theo3}
   If $(\Sigma,\mu)$ is a \LSS\  realization of $\Phi$, then
    the following are equivalent.
   \begin{itemize}
   \item[(i)]
     $(\Sigma,\mu)$ is a minimal \LSS\  realization of $\Phi$.
  \item[(ii)]
    The realization $(\Sigma,\mu)$ is \WR\  and it is observable.
  \item[(iii)]
     The state-space dimension of $\Sigma$ equals the rank
     of the Hankel-matrix of $\Phi$, i.e.
     $\dim \Sigma=\dim H_{\Phi}$.
 \end{itemize}   
   In addition, all minimal \LSS\  realizations of $\Phi$ are 
   isomorphic.
\end{Theorem}

\section{ Main results of the paper }
\label{sect:main_result}

The purpose of this section is to present the main results of the paper
formally.  The outline of this section is as follows.
In Subsection \ref{sect:main_result:part_real_problem} we state the partial-realization 
problem for \SLSS\ formally. In Subsection \ref{sect:main_result:hank}
we define the notion of Hankel sub-matrix, which will be needed for the statement
of the main results. In Subsection \ref{sect:main_result:alg0} we present the 
theorem characterizing the existence and minimality of partial realizations. In addition,
we present a simple algorithm for computing a minimal partial realization.
Finally, in Subsection \ref{sect:main_result:alg1} we present a Kalman-Ho-like
algorithm for computing a minimal partial realization.
\subsection{ Partial-realization problem} 
\label{sect:main_result:part_real_problem}
 Let $\Phi$ be a family of input-output maps admitting a generalized kernel representation.
    As we saw earlier, the Markov-parameters uniquely determine the input-output maps
    of a linear switched system. In addition, the realization problem for \LSS\ 
    can be reduced to the problem of finding a suitable representation of the Markov-parameters
    of the input-output maps.  However, in practice we can obtain only 
    Markov-parameters only up to some finite order. Hence, 
    the \emph{partial-realization problem} (with respect to the number of switches) arises.
    In order to present a formal problem formulation, the notion of (minimal)
    partial realization has to be introduced.
   \begin{Definition}[$N$-partial realization]
   \label{def:part_real}
    Let $\Phi$ be a family of input-output maps admitting a generalized kernel
    representation.
     Assume that $\Sigma$ is a \LSS\  of the form
     (\ref{lin_switch0}) and let $\mu: \Phi \rightarrow \mathcal{X}$.
     A realization $(\Sigma,\mu)$ is said to be
     an $N$-partial realization of $\Phi$, if for all $q,q_0 \in Q$,
     $f \in \Phi$, $j=1,2,\ldots,m$, and for all (possibly empty) sequences of
     discrete modes
     $q_1,q_2,\ldots,q_k \in Q$, of length at most $N$, i.e. $N \ge k \ge 0$, the
     following holds.
     \begin{equation}
     \label{sect:real:prop1:eq1.1}
     \begin{split}
        S_{q,q_0,j}(q_1q_2\cdots q_k)=&
        C_{q}A_{q_{k}}A_{q_{k-1}} \cdots A_{q_{1}}B_{q_{0}}e_{j} \\
        S_{f,q}(q_1q_2\cdots q_k)=&
        C_{q}A_{q_{k}}A_{q_{k-1}} \cdots A_{q_{1}}\mu(f)
     \end{split}
    \end{equation}
    The \LSS\  $\Sigma$ is said to be an $N$-partial realization of $\Phi$, if
    for there exists a map $\mu:\Phi \rightarrow \mathcal{X}$ such that the
    realization $(\Sigma,\mu)$ is a $N$-partial realization of $\Phi$.
    The realization $(\Sigma,\mu)$ is said to be a partial realization of $\Phi$,
    if it is a $N$-partial realization of $\Phi$ for some $N \in \mathbb{N}$.
   \end{Definition}
    From Lemma \ref{sect:io:prop1} it follows that $(\Sigma,\mu)$ is a realization
    of $\Phi$ if and only if $(\Sigma,\mu)$ is an $N$-partial realization of $\Phi$ for
    all $N \in \mathbb{N}$. Note that if $\Phi$ and $(\Sigma,\mu)$ is a $N$-partial
    realization of $\Phi$, then (\ref{sect:real:prop1:eq1.1}) holds for \emph{finitely
    many} Markov-parameters.
   \begin{Remark}[Terminology: Markov-parameters of order $N$]
    In the sequel, we will often refer to the Markov-parameters
    which are indexed by a sequence of discrete modes of
    length $N$ as \emph{Markov-parameters of order $N$}. That is, a Markov-parameter
    of order $N$ is either of the form $S_{q,q_0,j}(w)$ or $S_{f,q}(w)$ with $|w|=N$.
   \end{Remark}
    With this terminology, a \LSS is a $N$-partial realization of $\Phi$, if it
    recreates the Markov-parameters of $\Phi$ of order at most $N$. If $\Phi$ is finite,
    then the set of Markov-parameters of $\Phi$ of order at most $N$ is finite.
   \begin{Definition}[Minimal $N$-partial realization]
    A \LSS\ realization $(\Sigma,\mu)$ 
    is said to be a minimal $N$-partial realization of $\Phi$, if
    $(\Sigma,\mu)$ is a $N$-partial realization of $\Phi$ and for any other $N$-partial
    \LSS\ realization
    $(\Sigma^{'},\mu^{'})$  of $\Phi$, $\dim \Sigma \le \dim \Sigma^{'}$.
   \end{Definition}
    In plain words, a minimal $N$-partial realization of $\Phi$ is a $N$-partial
    realization of $\Phi$ with the smallest possible state-space dimension. 
   \begin{Problem}[Partial-realization problem]
   \label{prob:part:real}
    Let $\Phi$ be a family of input-output maps admitting a generalized kernel
    representation.
     The  partial-realization problem entails the following problems. 
     \begin{itemize} 
     \item
            Find conditions for existence of a $N$-partial realization
            of $\Phi$. Formulate an algorithm
            for computing a  $N$-partial realization from 
            generalized Markov-parameters of $\Phi$ of some bounded order. If $\Phi$
            is finite, this amounts to computing a $N$-partial realization from
            \emph{finitely many} Markov-parameters.
     \item
             Characterize minimal $N$-partial realizations of $\Phi$, find
             conditions for their existence and uniqueness.
     \item
            Find conditions under which a $N$-partial realization of $\Phi$
            is a complete realization of $\Phi$ in the sense of
            Definition \ref{switch_sys:real:def2}.
     \end{itemize}
   \end{Problem}
    We will devote the remaining part of the section to presenting the
    solution to the problem formulated above.
  We will show that it is possible to compute a minimal $N$-partial
  realization from a suitably chosen sub-matrix of the Hankel matrix of $\Phi$.
   In fact, this sub-matrix is
  formed by Markov-parameters of bounded order, and hence it is finite, if $\Phi$ is
  finite.
  If $N$ is large enough, the resulting $N$-partial realization will be a 
  complete realization of $\Phi$.
  In addition, if $\Phi$ is finite, then this sub-matrix will be finite as well.

 \subsection{Finite Hankel-matrices}
 \label{sect:main_result:hank}
  In order to state the main results rigorously, we need to define formally.
  certain finite sub-matrices of the Hankel-matrix. 
  To this end, recall from Definition \ref{main_result:lin:hank:arb:def} 
 the definition the Hankel-matrix of $\Phi$ and of the sets
 $I=\{1,2,\ldots,p\QNUM\}$, $\QNUM=|Q|$ and 
$J_{\Phi}=\Phi \cup (Q \times \{1,2,\ldots,m\})$.
 Fix  natural numbers $L,M \in \mathbb{N}$.
 We will use the following notation for sequences of discrete modes of 
 length at most $L$.
 \begin{Notation}
  We will denote by $Q^{ \le L}$ the set of all sequences of discrete modes of length at
  most $L$, i.e. $Q^{ \le L}=\{ w \in Q^{*} \mid |w| \le L\}$.
 \end{Notation}
  Notice that $Q^{\le L}$ is a finite set of cardinality $(|Q|^{L+1}-1)/(|Q|-1)$.
 We will define the sub-matrix of the Hankel-matrix $H_{\Phi}$ indexed by 
 sequences of bounded length as follows.
 \begin{Definition}[$H_{\Phi,L,M}$ sub-matrices of the Hankel-matrix]
 \label{sect:main_result:sub-hank:def1}
  Let $H_{\Phi,L,M}$ be the 
  sub-matrix of $H_{\Phi}$ formed by the intersections of all the rows indexed by 
  an index of the form $(v,i)$ with $i \in I$ and $v \in Q^{ \le L}$ and
  all the columns indexed by an index of the form $(w,j)$ with $j \in J_{\Phi}$ and 
  $w \in Q^{\le M}$. That is, 
  $H_{\Phi,L,M} \in \mathbb{R}^{(Q^{\le L} \times I) \times (Q^{\le M} \times J_{\Phi})}$ 
  and for each $(v,i) \in I \times Q^{\le L}$ and $(w,j) \in (Q^{\le M} \times J_{\Phi})$,
 \[ (H_{\Phi,L,M})_{(v,i),(w,j)}=(H_{\Phi})_{(v,i),(w,j)} \]
 \end{Definition}
   Notice that the entries of$H_{\Phi,L,M}$ are
   made up of Markov-parameters of order at most $K+L$.
   Hence, if $\Phi$ is finite, then the set $J_{\Phi}$ is finite and the
   matrix $H_{\Phi,L,M}$ is a finite matrix.
   Since the entries of the Hankel-matrix were defined in terms of the generalized
   Markov-parameters  of $\Phi$, the same is true for the entries of $H_{\Phi,N,M}$.
   Writing out the definition of the entries of $H_{\Phi,L,M}$ and
   using the enumeration of $Q$ fixed in (\ref{main_result:lin:hank:arb:def_enum}) of 
   Definition \ref{main_result:lin:hank:arb:def}, we get
   that for any $i$ of the form
$i=Kp+r+1$ where $K=0,1,\ldots, \QNUM-1$ and $r=0,1,\ldots, p-1$,
the entry $(H_{\Phi,L,M})_{(v,i),(w,j)}$ is defined as follows
\begin{equation}
\label{hankel:lin:part_def0.3}
  (H_{\Phi,L,M})_{(v,i),(w,j)}=
  \left\{\begin{array}{rl}
       (S_{\sigma_{K+1},q,z}(wv))_{r+1} & \mbox{ if }
                           j=(q,z) \in Q \times \{1,\ldots,m\} \\
       (S_{f,\sigma_{K+1}}(wv))_{r+1}   & \mbox{ if }
                                     j=f \in \Phi
   \end{array}\right.
\end{equation}
Here $(S_{\sigma_{K+1},q,z}(wv))_{r+1}$ and
$(S_{f,\sigma_{K+1}}(wv))_{r+1}$ denote the $r+1$th element of the
vectors $S_{\sigma_{K+1},q,z}(wv) \in \mathbb{R}^{p}$ and
$S_{f,\sigma_{K+1}}(wv) \in \mathbb{R}^{p}$ respectively.
    That is, if $\Phi$ is finite, then $H_{\Phi,L,M}$ represents a finite collection of
    Markov-parameters of $\Phi$. In turn
    the Markov-parameters of $\Phi$ are defined via the high-order derivatives of
    the elements of $\Phi$ with respect to the switching times. Hence, we get that the
    entries of $H_{\Phi,L,M}$ are just high-order derivatives with respect to the
    switching times of the elements of $\Phi$, and the order of these derivatives
    is bounded by $L+M$.

  \subsection{ Solution of the partial-realization problem}
  \label{sect:main_result:alg0}
    As we remarked in the introduction, there are essentially two ways to
    construct a partial \LSS\  realization. In 
    Algorithm \ref{alg0}  we present the first one which is conceptually the simplest one. 
   \begin{algorithm}
   \caption{\label{alg0}}
    \begin{algorithmic}[1]
    \STATE
        Assume that $n=\Rank H_{\Phi,N,N+1}$ 
        and let $\MORPH$ be a linear isomorphism
        from the column space of $H_{\Phi,N,N+1}$ to $\mathbb{R}^{n}$. For each
        column index $(w,j) \in Q^{\le N+1} \times J_{\Phi}$, 
        \emph{denote by
        $\COL_{w,j}$ the column of $H_{\Phi,N,N+1}$ indexed by $(w,j)$}.

    \STATE
       Define the \LSS\  realization $(\Sigma_N,\mu_N)$ as follows.
       
       $\Sigma_N$ is a \LSS\  of the form (\ref{lin_switch0}) such that
       \begin{itemize}
       \item
          The state-space of $\Sigma_N$ is $\mathbb{R}^{n}$.
       \item
          For each mode $q \in Q$, the matrix
          $B_q \in \mathbb{R}^{n \times m}$ satisfies
          \begin{equation}
          \label{alg0:eq-2}
              B_{q}=\begin{bmatrix} 
                \MORPH(\COL_{\epsilon,(q,1)}), &
                \MORPH(\COL_{\epsilon,(q,2)}), &  \cdots, &
                \MORPH(\COL_{\epsilon, (q,m)})
              \end{bmatrix}
          \end{equation}
       \item
         For each mode $q \in Q$, the matrix $C_{q} \in \mathbb{R}^{p \times n}$ satisfies
         \begin{equation}
         \label{alg0:eq-1}
         \begin{split}
             & C_q\MORPH(\COL_{w,j})=  \\
            & \begin{bmatrix} (H_{\Phi,N,N+1})_{(\epsilon,1),(w,j)}, &
                            (H_{\Phi,N,N+1})_{(\epsilon,2),(w,j)}, &
                            \cdots, &
                            (H_{\Phi,N,N+1})_{(\epsilon,p),(w,j)}
            \end{bmatrix}^{T}
         \end{split}
         \end{equation}
	 for each $w \in Q^{\le N+1}$.
       \item 
         For each mode $q \in Q$, the matrix $A_{q} \in \mathbb{R}^{n \times n}$ is
         the solution of the system of linear equations
         \begin{equation}
         \label{alg0:eq1}
            \MORPH^{-1}A_{q}\MORPH(\COL_{w,j})=\COL_{wq,j} 
            \mbox{ for each } w \in Q^{\le N}, j \in J_{\Phi} 
         \end{equation}
          If the system of equations (\ref{alg0:eq1}) does not have a solution, then
          abort.
       \end{itemize} 
    \STATE
        For each input-output map $f \in \Phi$, 
        \begin{equation}
        \label{alg0:eq-3}
          \mu_{N}(f)=\MORPH(\COL_{\epsilon,f})
        \end{equation}
    \STATE
       Return $(\Sigma_N,\mu_N)$, with $\Sigma_N=(\mathbb{R}^{n},\mathcal{U},\mathcal{Y},Q, \{(A_{q},B_q,C_q)\mid q \in Q\})$.
    \end{algorithmic}
   \end{algorithm}
   Informally, the realization $(\Sigma_N,\mu_N)$ returned by Algorithm \ref{alg0}
   is defined on an isomorphic copy of the column space of $H_{\Phi,N,N+1}$.
   The $j$th column of the matrix $B_{q}$ is the isomorphic copy of the
   column of $H_{\Phi,N,N+1}$ indexed by $(q,j)$. The matrices $C_{q}$ are such
   that if they are interpreted as linear maps, then each $C_{q}$ maps the isomorphic
   copy of a column of $H_{\Phi,N,N+1}$ to its first $p$ rows, i.e. rows indexed by
   $(\epsilon,1),(\epsilon,2), \ldots, (\epsilon,p)$, where $\epsilon$ denotes the
   empty word.
   The matrices $A_{q}$ realize a shift on the columns of $H_{\Phi,N,N+1}$;
   by interpreting $A_{q}$ as a map on the columns of $H_{\Phi,N,N+1}$, $A_{q}$ maps
   the column indexed by an index $(w,j) \in Q^{\le {N}} \times J_{\Phi}$ to the column
   indexed by $(wq,j)$. Finally, the value $\mu_{N}(f)$ of the map $\mu_N$ for the
   input-output map $f \in \Phi$ equals the column of $H_{\Phi,N,N+1}$ indexed by 
   $(\epsilon,f)$.
     Notice that Algorithm \ref{alg0} may fail to return
   a realization, as  (\ref{alg0:eq1}) need not always have a solution.
   \begin{Remark}
    If $\Phi$ consists of finitely many input-output maps, then Algorithm \ref{alg0}
    is a indeed an effective procedure and it can be implemented as a numerical algorithm.
   \end{Remark}

   We state the partial-realization theorem as a theorem providing sufficient
   conditions for Algorithm \ref{alg0} to yield a $2N+1$-partial, \WR\ and observable
   realization of $\Phi$.
  \begin{Theorem}[Existence of partial realization of \SLSS]
  \label{part_real_lin:theo1}
   Let $\Phi$ be a family of input-output maps admitting a generalized kernel
   representation. 
  \begin{itemize}
  \item{\textbf{Existence  and computability of a partial realization} \\}
  Assume that for some $\mathbb{N} \ni N>0$,
  \begin{equation}
  \label{part_real_lin:theo1:eq1}
  \Rank H_{\Phi,N,N}=\Rank H_{\Phi,N+1,N}=\Rank H_{\Phi,N,N+1}
  \end{equation}
   holds.  
  Then Algorithm \ref{alg0} returns a realization 
  $(\Sigma_{N},\mu_{N})$, and $(\Sigma_N,\mu_N)$ is a 
  \WR\  and observable $2N+1$-partial realization of $\Phi$.

  \item{\textbf{Existence of a complete realization} \\}
  If $\Rank H_{\Phi,N,N}=\Rank H_{\Phi}$, then 
  (\ref{part_real_lin:theo1:eq1}) holds and 
   the $2N+1$-partial realization $(\Sigma_N,\mu_N)$ 
  is a minimal realization of $\Phi$.
  The condition $\Rank H_{\Phi,N,N}=\Rank H_{\Phi}$ holds for a given $N$, if
  there exists a \LSS\  realization $(\Sigma,\mu)$ of $\Phi$ such that
  $\dim \Sigma \le N+1$.
  \end{itemize}
\end{Theorem}
 The proof of Theorem \ref{part_real_lin:theo1} is presented in Section \ref{sect:lin_alg}.
 Theorem \ref{part_real_lin:theo1} implies that  
 if $\Phi$ has a realization by a \LSS\ , 
 then the parameters of a minimal \LSS\  realization of
  $\Phi$ are computable
 from, and hence are determined by,
 finitely many time derivatives with respect to the switching times. 
\begin{Remark}[Partial versus complete realization]
\label{part_real_lin:theo1:remark}
 The numerical example in Section \ref{sect:num} 
provides an example of a family of input-output maps
 $\Phi$ realizable by a \LSS, such that $\Phi$ has the following property. 
 For some $N$ the family $\Phi$ satisfies 
 (\ref{part_real_lin:theo1:eq1}) but not the condition $\Rank H_{\Phi,N,N}=\Rank H_{\Phi}$.
 That is, there is an $N \in \mathbb{N}$ 
 such that Theorem \ref{part_real_lin:theo1} yields a
 partial realization of $\Phi$ which is not a complete realization. 
 That is, the two statements of Theorem \ref{part_real_lin:theo1} indeed describe two
 separate cases.
\end{Remark}

 Theorem \ref{part_real_lin:theo1} allows us to formulate the following
 characterization of minimal partial \LSS\ realizations.
 \begin{Theorem}[Minimal partial realization]
 \label{part_real_lin:theo1.1}
  With the notation of Theorem \ref{part_real_lin:theo1}, if $\Phi$ satisfies
 (\ref{part_real_lin:theo1:eq1}), then the following holds.
 \begin{enumerate}
 \item
  \label{part_real_lin:theo1.1:res1}
   A minimal $2N+1$ partial realization of $\Phi$ exists, in fact, 
   the realization $(\Sigma_N,\mu_N)$ returned by Algorithm \ref{alg0}
   is a minimal $2N+1$ partial realization of $\Phi$.
 \item
  \label{part_real_lin:theo1.1:res2}
   Any minimal $2N+1$ partial realization of $\Phi$ is \WR\  and observable and
   it is of dimension $\Rank H_{\Phi,N,N}$. 
 \item
 \label{part_real_lin:theo1.1:res3}
  All minimal $2N+1$ partial realizations of $\Phi$ are isomorphic.
 \end{enumerate}
 \end{Theorem}
 The proof of Theorem \ref{part_real_lin:theo1.1} is presented in Section \ref{sect:lin_alg}.
 \begin{Remark}
  The reader might wonder why we talk about $2N+1$-partial realizations in
  Theorem \ref{part_real_lin:theo1} and \ref{part_real_lin:theo1.1}. The reason behind it
  is that the finite Hankel-matrix $H_{\Phi,N,N+1}$ is
  formed by values of the Markov-parameters of order at most $2N+1$.
   We would like the \LSS\ realization obtained
   from $H_{\Phi,N,N+1}$ to recreate at least those Markov-parameters which are the entries of
   the matrix $H_{\Phi,N,N+1}$. But this means that the
   \LSS\ realization obtained from $H_{\Phi,N,N+1}$ must be
   a $2N+1$-partial realization of $\Phi$.
 \end{Remark}


\subsection{ A Kalman-Ho-like partial-realization algorithm for \LSS}
 \label{sect:main_result:alg1}
 While Algorithm \ref{alg0} is theoretically attractive 
because of its simplicity, it is not necessarily the most suitable one for numerical
implementation. 
In Algorithm \ref{PartReal} 
we present an alternative algorithm for computing a 
\WR\ and observable $2N+1$-partial realization of $\Phi$. The algorithm
is based on factorization of the Hankel-matrix $H_{\Phi,N+1,N}$ and
it is similar to the Kalman-Ho algorithm. 

\begin{algorithm}
\caption{\textbf{ComputePartialRealization($H_{\Phi,N+1,N}$)}}
\label{PartReal}
\begin{algorithmic}[1]
 \STATE
 Compute a decomposition of $H_{\Phi,N+1,N}$
    \[ H_{\Phi,N+1,N}=OR \]
    where $O \in \mathbb{R}^{(Q^{\le N+1} \times \{1,2,\ldots,p\QNUM\}) \times r}$ and
    $R \in \mathbb{R}^{r \times (Q^{\le N} \times J_{\Phi})}$ are matrices such that
    $r=\Rank R=\Rank O=\Rank H_{\Phi,N+1,N}$.
 
 \STATE
  Recall from (\ref{main_result:lin:hank:arb:def_enum}) the enumeration of $Q$.
  For each $q \in Q$, $q=\sigma_i$, for some $i=1,2,\ldots,\QNUM$
  define the matrix $\widetilde{C}_q \in \mathbb{R}^{p \times r}$ by
 \[  \widetilde{C}_q=\begin{bmatrix} O_{(\epsilon,(i-1)p+1),.}^{T}, & O_{(\epsilon,(i-1)p+2),.}^{T}, & 
    \ldots, & O_{(\epsilon,ip),.}^{T} \end{bmatrix}^{T} 
  \]
   where $O_{k,.}$ denotes the row of $O$ indexed by $k$. 

\STATE
  For each $q \in Q$, define the matrix $\widetilde{B}_{q} \in \mathbb{R}^{r \times m}$
  by
  \[ \widetilde{B}_{q}=\begin{bmatrix}
                        R_{.,(\epsilon,(q,1))}, & R_{.,(\epsilon,(q,2))}, & \cdots &
                        R_{.,(\epsilon,(q,m))}
                        \end{bmatrix}
  \]
  where $R_{.,(\epsilon,(q,j))}$ stands for the column of $R$ indexed by 
  $(\epsilon,(q,j))$ for $q \in Q$ and $j=1,2,\ldots,m$.
  
\STATE
  Define the map $\widetilde{\mu}_{N}: \Phi \rightarrow \mathbb{R}^{r}$ as
  \[ \forall f \in \Phi: \widetilde{\mu}_{N}(f)=R_{., (\epsilon,f)} \]
  where $R_{.,(\epsilon,f)}$ stands for the column of $R$ indexed by 
  $(\epsilon,f)$ for $f \in \Phi$.
\STATE  
For each $q \in Q$ let the matrix $\widetilde{A}_{q} \in \mathbb{R}^{r \times r}$  
be the solution of equation
\begin{equation}
\label{ComputePartialReal}
\bar{\Gamma}\widetilde{A}_{q}=\bar{\Gamma}_{q} 
\end{equation}
where
$\bar{\Gamma}, \bar{\Gamma}_{q} \in \mathbb{R}^{(Q^{\le N} \times \{1,2,\ldots,p\QNUM\}) \times r}$ are matrices of the
form
\begin{equation*} 
\begin{split}
 \bar{\Gamma}_{(v,i),j}=O_{(v,i),j} \mbox{ and } 
  (\bar{\Gamma}_{q})_{(v,i),j}=O_{(qv,i),j} 
\end{split}
\end{equation*}
for all  \( (v,i) \in Q^{\le N} \times \{1,2,\ldots,p\QNUM\}, j=1,2,\ldots,r \).
That is, $\bar{\Gamma}$ is obtained from $O$ by deleting all the rows indexed by
pairs $(v,i)$ where $v$ is a sequence of discrete modes of length $N+1$. The matrix
$\bar{\Gamma}_{q}$ is the shifted version of $O$, i.e. its row indexed by $(v,i)$
is the row of $O$ indexed by $(qv,i)$.
\STATE
  If there no unique solution to (\ref{ComputePartialReal}) then return
  $NoRealization$. Otherwise return $(\widetilde{\Sigma}_{N},\widetilde{\mu}_{N})$
  where $\widetilde{\Sigma}_{N}$ is a \LSS\  of the form
  \[ \widetilde{\Sigma}_{N}=(\mathbb{R}^{r}, \mathcal{U},\mathcal{Y}, \{( \widetilde{A}_q,\widetilde{B}_q,\widetilde{C}_q) \mid q \in Q\}) \]
\end{algorithmic}
\end{algorithm}
The result of Algorithm \ref{PartReal} is described in the theorem below.
\begin{Theorem}[Partial realization algorithm]
\label{part_real_lin:theo2}
With the notation above the following holds.
\begin{enumerate}
\item 
 \label{part_real_lin:theo2:res1}
  Assume that form some $N>0$, (\ref{part_real_lin:theo1:eq1}) of
  Theorem \ref{part_real_lin:theo1} holds.
  Then Algorithm \ref{PartReal} always returns a 
  \LSS\  realization $(\widetilde{\Sigma}_{N},\widetilde{\mu}_N)$ 
  and $(\widetilde{\Sigma}_{N},\widetilde{\mu}_N)$ is a
  minimal $2N+1$-realization of $\Phi$. In fact, 
  $(\widetilde{\Sigma}_{N},\widetilde{\mu}_{N})$ is isomorphic to the
  \LSS\  realization $(\Sigma_N,\mu_N)$ of
  Theorem \ref{part_real_lin:theo1}.

\item
 \label{part_real_lin:theo2:res2}
   If for some $N >0$,
   \( \Rank H_{\Phi,N,N}=\Rank H_{\Phi} \),
    then (\ref{part_real_lin:theo1:eq1}) holds, and Algorithm \ref{PartReal} returns 
   a minimal realization $(\widetilde{\Sigma}_{N},\widetilde{\mu}_{N})$ of $\Phi$.
   If there exists a 
   \LSS\  realization $(\Sigma,\mu)$ of 
   $\Phi$, such that $\dim \Sigma \le N+1$, 
   then $\Rank H_{\Phi,N,N}=\Rank H_{\Phi}$
   holds and 
   the realization $(\widetilde{\Sigma}_{N},\widetilde{\mu}_N)$ 
   returned by Algorithm \ref{PartReal} is a minimal
   realization of $\Phi$.
\end{enumerate}
\end{Theorem}
 The proof of Theorem \ref{part_real_lin:theo2} is presented in Section \ref{sect:lin_alg}.

 Theorem \ref{part_real_lin:theo2}
 has immediate   consequences for systems identification and model reduction of linear switched
 systems, by guaranteeing correctness of the realization algorithm described in
 Algorithm \ref{PartReal}.
In Section \ref{sect:concl} we discuss the potential applications of Theorem
\ref{part_real_lin:theo1} and Theorem \ref{part_real_lin:theo2} in more detail.

\begin{Remark}[Factorization of $H_{\Phi,N+1,N}$]
 There are many algorithms to compute a factorization of $H_{\Phi,N+1,N}$ and
 Algorithm \ref{PartReal} can work with any of these algorithms. In particular,
 the factorization of $H_{\Phi,N+1,N}$ can be done through singular value
 decomposition. More precisely, assume that $\Phi$ is finite, and let 
 $H_{\Phi,N+1,N}=U\Sigma V$ be the singular
 value decomposition of $H_{\Phi,N+1,N}$. 
 Notice that if $\Phi$ is finite, then  the index sets of $H_{\Phi,N+1,N}$ are finite and 
 the singular value
 decomposition can be computed by enumerating the index sets of $H_{\Phi,N+1,N}$
 and viewing $H_{\Phi,N+1,N}$ as a usual finite matrix.
 The choice $O=U\Sigma^{1/2}$ and
 $R=\Sigma^{1/2}V$ yields a decomposition $H_{\Phi,N+1,N}=OR$ satisfying the conditions
 of Algorithm \ref{PartReal}. 
 Algorithm \ref{PartReal} has been implemented, and the factorization used in the
 implementation is precisely the singular value decomposition of $H_{\Phi,N+1,N}$.
\end{Remark}

\section{ Proof of the partial realization theorems for LSS }
\label{sect:lin_alg}
The goal of this section is to present the proof of Theorem \ref{part_real_lin:theo1}, \ref{part_real_lin:theo1.1}, and \ref{part_real_lin:theo2}.
Both proofs 
rely heavily on the relationship
between rational representations and \SLSS\  described in
\cite{MP:RealForm,MP:BigArticle}.
The outline of the section is the following. In Subsection \ref{sect:lin:real}
we will present a short overview of the relationship between rational representations and
\SLSS\ .
Subsection \ref{sect:lin:part_real} presents the proof of Theorem \ref{part_real_lin:theo1}, 
\ref{part_real_lin:theo1.1}, and \ref{part_real_lin:theo2}.
The reader is advised to review Appendix \ref{sect:pow} and Appendix \ref{sect:pow:part}
before reading this section.
Throughout the section we will tacitly use the notation and terminology of Appendix \ref{sect:pow} and Appendix \ref{sect:pow:part}.

\subsection{ Rational representations and \SLSS\ }
\label{sect:lin:real}
 Below we will present a brief overview of the relationship between formal power series
 and their representations and \SLSS\ . For a more detailed presentation see
 \cite{MP:RealForm,MP:BigArticle}.
 Informally, the relationship is as follows.
 \begin{itemize}
 \item
   Let $\Phi$ be a family of input-output maps and assume that $\Phi$ has a 
   generalized kernel representation. Then we can construct a family of formal
   power series $\Psi_{\Phi}$ from the generalized Markov-parameters of $\Phi$.
   The details of the construction will be presented below. The family $\Psi_{\Phi}$
   will be referred to as \emph{the family of formal power series associated with $\Phi$}.
\item
   If $(\Sigma,\mu)$ is a \LSS\  realization of $\Phi$,
   then we can construct a rational representation $R_{\Sigma,\mu}$ of $\Psi_{\Phi}$ from
   the parameters of $(\Sigma,\mu)$. The converse is also true; if $R$ is a representation
   of $\Psi_{\Phi}$, then we can construct a realization $(\Sigma_{R},\mu_{R})$ of 
   $\Phi$ from the parameters of $R$. 
   The details of the construction will be presented below.
   We will call $R_{\Sigma,\mu}$ the
   \emph{representation associated with $(\Sigma,\mu)$} and we will call
    $(\Sigma_{R},\mu_{R})$ the \emph{realization associated with $R$}.
\end{itemize} 
 We will start with defining the family of formal power series $\Psi_{\Phi}$ associated
 with $\Phi$.
 To this end, recall from Definition \ref{linswitch_intro:def:markov}
 the definition of Markov-parameters
 $S_{q,q_0,j}(w)$, $S_{f,q}(w)$ of $\Phi$, for each pair of discrete modes $q,q_0 \in Q$, 
 input-output map $f \in \Phi$, index $j=1,2,\ldots,m$, and a word $w \in Q^{*}$.
 It is easy to see that the maps
 $S_{q,q_0,j}: Q^{*} \ni w \rightarrow S_{q,q_0,j}(w) \in \mathbb{R}^{p}$
 and $S_{f,q}:Q^{*} \ni w \rightarrow S_{f,q}(w) \in \mathbb{R}^{p}$
 define formal power series $S_{q,q_0,j}$ and $S_{f,q}$ in $\mathbb{R}^{p}\ll Q^{*} \gg$.
  Recall the enumeration of the set of discrete modes
  $Q$ defined in (\ref{main_result:lin:hank:arb:def_enum}).
  That is, $Q$ assumed
  to have $\QNUM$ elements given by the distinct
  elements $\sigma_1,\sigma_2,\ldots, \sigma_{\QNUM}$.
  For each discrete mode $q \in Q$, index $j=1,2,\ldots, m$, and input-output map $f \in \Phi$ define
  the formal power series
  $S_{q,j}, S_{f} \in \mathbb{R}^{p\QNUM}\ll Q^{*} \gg$ 
  as follows; for each word $w \in Q^{*}$ let
  \begin{equation}
  \label{sect:real:form-1}
  \begin{array}{ll}
    S_{q,j}(w)&=\begin{bmatrix}
       (S_{\sigma_{1},q,j}(w))^{T}, & (S_{\sigma_{2},q,j}(w))^{T}, & \cdots, & 
       (S_{\sigma_{\QNUM},q,j}(w))^{T}
       \end{bmatrix}^{T} \in \mathbb{R}^{p\QNUM},  \\
    S_{f}(w)&=\begin{bmatrix}
                 (S_{f,\sigma_{1}}(w))^{T}, & 
		 (S_{f,\sigma_{2}}(w))^{T}, & \cdots, &
		 (S_{f,\sigma_{\QNUM}}(w))^{T}
	       \end{bmatrix}^{T} \in \mathbb{R}^{p\QNUM}	 
    \end{array}
  \end{equation}	       
 That is, the values of the
 formal power series $S_{q,j}$ are obtained by stacking up
 the values of $S_{\sigma_{i},q,j}$ for $i=1,2,\ldots, \QNUM$.
 Similarly, the values of $S_{f}$ are obtained by stacking up
 the values of $S_{f,\sigma_{i}}$ for $i=1,2,\ldots, \QNUM$.
  Define the set $J_{\Phi}=\Phi \cup \{ (q,z) \mid q \in Q, z=1,2,\ldots,m \}$.
  Define the \emph{indexed set of formal power series associated with $\Phi$}
  as the following family of formal power series 
    \begin{equation}
    \label{sect:real:form-2}
    \Psi_{\Phi}=\{ S_{j} \in \mathbb{R}^{p\QNUM} \ll Q^{*} \gg \mid j \in J_{\Phi} \}
     \end{equation}
  \begin{Remark}[Equivalence of definitions of the Hankel-matrix]
  \label{sect:real:hank}
  It is easy to see that
  the Hankel-matrix $H_{\Psi_{\Phi}}$ of the family of formal power series $\Psi_{\Phi}$ 
  is identical to the
  Hankel-matrix $H_{\Phi}$ of $\Phi$ as defined in
  Definition \ref{main_result:lin:hank:arb:def}
  \end{Remark}

  Next we define the rational representation
  $R_{\Sigma,\mu}$ associated with a \LSS\ realization $(\Sigma,\mu)$.
  Let $\Sigma$ be a \LSS\  of the form
  (\ref{lin_switch0}) and assume that 
  $\mu:\Phi \rightarrow \mathcal{X}$ is a map assigning to
  each element $f \in \Phi$ an initial state of $\Sigma$.
  \begin{Construction}[Representation associated with a realization]
  \label{sect:real:const1}
  Define the 
  \emph{representation associated with $(\Sigma,\mu)$} as the rational $p\QNUM-J_{\Phi}$ 
  representation
  \[
   R_{\Sigma,\mu}=
      (\mathcal{X}, \{ A_{q} \}_{q \in Q}, \widetilde{B},\widetilde{C})
  \]
  The various components of $R_{\Sigma,\mu}$ are defined as 
  follows.  
  \begin{itemize}
  \item{\textbf{State-space $\mathcal{X}$. }}
  The state-space of the rational representation 
  $R_{\Sigma,\mu}$ is the same as the state-space
  of $\Sigma$, i.e. $\mathbb{R}^{n}=\mathcal{X}$.
  \item{\textbf{Alphabet. }}
  The representation $R_{\Sigma,\mu}$ is defined over the
  alphabet which equals the set of discrete modes $Q$. 
  \item{\textbf{State-transition maps (matrix) $\{A_{q} \in \mathbb{R}^{n}\}_{q \in Q}$.  }}
  For each discrete mode $q \in Q$, the corresponding state-transition matrix 
  $A_{q}$ of $R_{\Sigma,\mu}$ is identical to the
  matrix $A_{q}$ of $\Sigma$. 
  \item{\textbf{Readout matrix (map) $\widetilde{C} \in \mathbb{R}^{p\QNUM \times n}$. }}
  The readout matrix $\widetilde{C}$ 
   is obtained by vertically
   "stacking up" the matrices $C_{\sigma_1},\ldots, C_{\sigma_{\QNUM}}$
   in this order from top to bottom. 
   That is, the $p\QNUM \times n$ matrix
   $\widetilde{C}$ is of the form
   \[ \widetilde{C}=\begin{bmatrix} 
     C_{\sigma_{1}}^{T} & C_{\sigma_{2}}^{T}, & \cdots & C_{\sigma_{\QNUM}}^{T} 
     \end{bmatrix}^{T}.
   \]
   \item{\textbf{Initial states}. }
   The indexed set of the initial states of $R_{\Sigma,\mu}$
   is of the form
   \[ \widetilde{B}=\{ \widetilde{B}_{j} \in \mathcal{X} \mid j \in J_{\Phi}\}, \] 
   i.e. it is indexed by the elements of the index set $J_{\Phi}$.
   Its elements are defined by
   $\widetilde{B}_{f}=\mu(f)$ if $f$ is an element of
   the family $\Phi$, and 
   $\widetilde{B}_{(q,l)}=(B_{q})_{.,l}$, i.e.
   $\widetilde{B}_{(q,l)}$ is the $l$th column of $B_{q}$, for all
   $q \in Q$ and $l=1,2,\ldots m$. 
   \end{itemize}
  \end{Construction}
    Conversely, below we will construct a \LSS\ realization 
    $(\Sigma_R,\mu_R)$ from any rational representation $R$
    satisfying some mild conditions. 
    We would like to note
    that these conditions are automatically satisfied by
    a suitably chosen isomorphic copy of any rational representation of 
    $\Psi_{\Phi}$. 
   The definition goes as follows.
   \begin{Construction}[Realization associated with a representation \\]
   \label{sect:real:const2}
    \ \  \\
   Consider a $p\QNUM-J_{\Phi}$ representation $R$ of the following form
    \begin{equation} 
  \label{sect:real:repr2sys}
    R=
      (\mathcal{X}, \{ A_{q} \}_{q \in Q}, \widetilde{B},\widetilde{C})
   \end{equation}
   and assume that the following holds.
   \begin{enumerate}
   \item
    The state-space is of the form $\mathcal{X}=\mathbb{R}^{n}$ for some $n >0$.
    \footnote{
    If  $\mathcal{X}=\mathbb{R}^{n}$ does not
    hold, then replace $R$ with the isomorphic copy
    $\RMORPH R$ defined in (\ref{repr:euclid:eq1}), Remark
    \ref{repr:state_space:rem} whose state-space is 
    $\mathbb{R}^{n}$. Since $\RMORPH R$ and $R$ are isomorphic, 
    if $R$ is a representation of $\Psi_{\Phi}$, then
    $\RMORPH R$ will be a representation of $\Psi_{\Phi}$ as well.
    }
   \item The readout map $\widetilde{C}$ takes its
   values in $\mathbb{R}^{p\QNUM}$.
   \item
    The set of initial states $\widetilde{B}$ of $R$ is indexed by the 
   index set $J_{\Phi}=\Phi \cup \{ (q,j) \mid q \in Q, j=1,2,\ldots,m\}$.
    \end{enumerate}
    Define the \emph{\LSS\  realization
    $(\Sigma_{R},\mu_{R})$ associated with R}
    as follows. Let $\Sigma_{R}$ be of the form 
    (\ref{lin_switch0}), where
    \begin{itemize}
    \item{\textbf{State-space. }}
    The
    state-space $\mathbb{R}^{n}$ of $\Sigma_{R}$ is the same as that of $R$.
    \item{\textbf{System matrices $\{A_q \in \mathbb{R}^{n \times n}\}_{q \in Q}$ of $\Sigma_{R}$. }}
     For each discrete mode $q \in Q$,
      the matrix $A_{q}$ of $\Sigma_{R}$ is identical to
     the state-transition matrix $A_{q}$ of $R$.
    \item{\textbf{System matrices $\{C_q \in \mathbb{R}^{p \times n}\}_{q \in Q}$ of $\Sigma_{R}$. }}
    For any discrete state $q \in Q$ of the form 
    $q=\sigma_i$ for some $i=1,2,\ldots,\QNUM$ and for any  
    $l=1,2,\ldots,p$, the $l$th row of $C_{q}$ equals the $p(i-1)+l$-th row of $\widetilde{C}$, i.e.
   \[ 
      \widetilde{C}=\begin{bmatrix} C_{\sigma_{1}}^{T}, & C_{\sigma_{2}}^{T},  &  \cdots & C_{\sigma_{N}}^{T} \end{bmatrix}^{T}. \]
   \item{\textbf{System matrices $\{B_{q} \in \mathbb{R}^{n \times m}\}_{q \in Q}$ of $\Sigma_{R}$. }}
    For each discrete mode $q \in Q$ the $n \times m$
    matrix $B_{q}$ is obtained as follows; the $l$th column
    of $B_{q}$ equals the initial state $\widetilde{B}_{q,l}$
    for all $l=1,2,\ldots, m$. That is
    \[ B_{q}=\begin{bmatrix} \widetilde{B}_{(q,1)}, & 
         \widetilde{B}_{(q,2)}, & \cdots, & \widetilde{B}_{(q,m)}
       \end{bmatrix}
    \]
   \end{itemize}
    The map $\mu_{R}: \Phi \rightarrow \mathcal{X}$ assigns to
    each element $f$ of $\Phi$ the initial state of $R$ indexed
    by $f$, i.e. 
    \[ \mu_{R}(f)=\widetilde{B}_{f} \mbox{ for all } f \in \Phi. \]
  \end{Construction}  
 The following theorem states the relationship between representations
and realizations formally.
 \begin{Remark}
   If we apply Construction \ref{sect:real:const1} to $(\Sigma_R,\mu_R)$ then
    the resulting representation $R_{\Sigma_{R},\mu_{R}}$ coincides with $R$, i.e.
    $R_{\Sigma_{R},\mu_{R}}=R$. Conversely, if we apply Construction \ref{sect:real:const2},
    to the representation $R_{\Sigma,\mu}$ associated with a \LSS\ realization $(\Sigma,\mu)$,
    then we get $(\Sigma,\mu)$ back, i.e.
   $\Sigma_{R_{\Sigma,\mu}}=\Sigma$, $\mu_{R_{\Sigma,\mu}}=\mu$.
 \end{Remark}

\begin{Theorem}[\cite{MP:Phd,MP:BigArticle,MP:RealForm}]
\label{sect:real:theo1}
   Let $\Phi$ be a family of input-output maps and assume
   that $\Phi$ admits a generalized kernel representation. 
   With the notation and assumptions above the following holds.
     \begin{itemize}
     \item
      The realization $(\Sigma,\mu)$ is a \LSS\ 
    realization of $\Phi$ if and only if $R_{\Sigma,\mu}$
    is a rational representation of $\Psi_{\Phi}$.
    \item
    The representation $R$
    is a representation of $\Psi_{\Phi}$ if and only if 
    $(\Sigma_{R},\mu_{R})$ is a realization  of $\Phi$.
  
  \item
    The realization $(\Sigma,\mu)$ is \WR\ if and only if
    $R_{\Sigma,\mu}$ is reachable, and $(\Sigma,\mu)$ is observable if and only if
    $R_{\Sigma,\mu}$ is observable.
    Similarly, $(\Sigma_{R},\mu_{R})$ is \WR\ if and only if $R$ is reachable, and
    $(\Sigma_R,\mu_R)$ is observable if and only if $R$ is observable.

   \item  
     If $(\Sigma,\mu)$ is a minimal realization of $\Phi$, then
     $R_{\Sigma,\mu}$ is a minimal representation of
     $\Psi_{\Phi}$. Conversely, if $R$ is a minimal
     representation of $\Psi_{\Phi}$, then
     $(\Sigma_{R},\mu_{R})$ is a minimal realization of $\Phi$.
   \item
  The map
  $\MORPH:(\Sigma,\mu) \rightarrow (\Sigma^{'},\mu^{'})$ is a
  \LSS\  morphism if and only if
  $\MORPH:R_{\Sigma,\mu} \rightarrow R_{\Sigma^{'},\mu^{'}}$ is a
  representation morphism.\footnote{Notice that the state-space of $(\Sigma,\mu)$ coincides
  with that of $R_{\Sigma,\mu}$ and the state-space of $(\Sigma^{'},\mu^{'})$ coincides 
  with that of $R_{\Sigma^{'},\mu^{'}}$. Hence $S$ can indeed be viewed as a representation
  morphism, provided that it commutes with the matrices of the representations involved}
  In addition, $S$ is a representation isomorphism if
  and only if $S$ is a \LSS\  isomorphism.
   Conversely, if the representations $R$ and $R^{'}$ satisfy the assumptions of
   Construction \ref{sect:real:const2}, then 
   $\RMORPH: R \rightarrow R^{'}$ is a representation morphism if and only if
   $\RMORPH: (\Sigma_R,\mu_R) \rightarrow (\Sigma_{R^{'}},\mu_{R^{'}})$ is a
   \LSS\ morphism. Moreover, in this case $\RMORPH$ is a representation
    isomorphism if and only if it is a \LSS\ isomorphism.
    \end{itemize}
 \end{Theorem}

\subsection{ Proofs of the results on partial realization of \SLSS\  }
\label{sect:lin:part_real}
 The goal of this section is to present the proof of Theorem 
 \ref{part_real_lin:theo1}, Theorem \ref{part_real_lin:theo1.1} and
 Theorem \ref{part_real_lin:theo2}. Before proceeding to the proof, we will
 need the following result on the relationship between the partial-realization
 problem for \SLSS\  and the partial-realization problem for rational representations.
 \begin{Theorem}
 \label{part_real_lin:theo1/2}
  Assume that $\Phi$ is a family of input-output maps admitting a generalized kernel
  representation.  A \LSS\  realization $(\Sigma,\mu)$ is an 
  $2N+1$-partial realization of
  $\Phi$ if and only of the associated representation $R_{\Sigma,\mu}$ is a $2N+1$-partial
  representation of the family of formal power series $\Psi_{\Phi}$.
  Conversely, a representation $R$ is a $2N+1$-partial representation of the family
  $\Psi_{\Phi}$ of formal power series associated with $\Phi$,
  if and only if the \LSS\ realization
  $(\Sigma_{R},\mu_{R})$ associated with $R$ is a $2N+1$-partial realization of $\Phi$.
 \end{Theorem}
 \begin{pf}
   The second statement of the theorem follows from the first one by noticing that
   $R=R_{\Sigma_{R},\mu_{R}}$. Hence, it is enough to prove the first statement of
   the theorem.
   It follows easily from Definition \ref{def:part_real} and the definition of
   the formal power series $S_{f}, S_{q_0,j}$, $f \in \Phi$, $q_0 \in Q$, $j=1,2,\ldots,m$,
   that $(\Sigma,\mu)$ is an $2N+1$-partial realization of $\Phi$, if and only if
   for all $w \in Q^{\le 2N+1}$, $q_0 \in Q$, $f \in \Phi$, $j=1,2,\ldots,m$
   \begin{equation}
   \label{sect:real:pow:eq3}
   \begin{split}
    S_{q_{0},j}(w)&=
    \begin{bmatrix} C_{\sigma_{1}}^{T} & C_{\sigma_{2}}^{T} & \cdots & C_{\sigma_{\QNUM}}^{T}
    \end{bmatrix}^{T}
    A_{w}B_{q_{0}}e_{j}  \\
    S_{f} (w)&=
    \begin{bmatrix} C_{\sigma_{1}}^{T} & C_{\sigma_{2}}^{T} & \cdots & C_{\sigma_{\QNUM}}^{T}
    \end{bmatrix}^{T}
    A_{w}\mu(f) \\
 \end{split}
\end{equation}
 From the definition of $R_{\Sigma,\mu}$ it follows directly that 
 (\ref{sect:real:pow:eq3}) is equivalent to $R_{\Sigma,\mu}$ being an $2N+1$-representation
 of $\Psi_{\Phi}$.
 \end{pf}
 Now we are ready to present the proof of Theorem \ref{part_real_lin:theo1}.
\begin{pf}[Proof of Theorem \ref{part_real_lin:theo1}]
 First, recall from Definition \ref{part_real_pow:submatrix} the definition of the
 matrix $H_{\Phi_{\Psi},K,L}$ for some $K,L > 0$. It is easy to see that the matrix
 $H_{\Phi,K,L}$ as defined in Definition \ref{sect:main_result:sub-hank:def1} and the matrix
 $H_{\Phi_{\Psi},K,L}$ coincide, in particular $\Rank H_{\Phi,K,L}=\Rank H_{\Psi_{\Phi},K,L}$.
 Recall from Remark \ref{sect:real:hank} that 
 the matrices $H_{\Phi}$ and $H_{\Phi_{\Psi}}$ coincide and hence their ranks are equal as well.
 We will prove the two statements of the theorem separately.

 \textbf{Proof of existence  and computability of a partial realization} \\
 We will apply Theorem \ref{part_real_pow:theo} to $\Psi_{\Phi}$.
 Using the remark above,  the condition 
 (\ref{part_real_lin:theo1:eq1})
 can be rewritten as  (\ref{part_theo:eq1}) of Theorem \ref{part_real_pow:theo}.
Consider the $p\QNUM-J_{\Phi}$ representation 
$R_{N}=(\IM H_{\Phi,N,N+1}, \{A_{q}\}_{q \in Q}, B,C)$ defined in
Theorem \ref{part_real_pow:theo}.
Since (\ref{part_theo:eq1}) holds, we get that $R_{N}$ is well-defined and it is a
$2N+1$-partial representation of $\Psi_{\Phi}$, and $R_{N}$ is reachable and observable.
Consider Algorithm \ref{alg0} and the isomorphism $\MORPH$ defined there. 
It is clear that 
$\MORPH$ maps the column space of $H_{\Psi_{\Phi},N,N+1}=H_{\Phi,N,N+1}$ to
$\mathbb{R}^{n}$ with $n=\Rank H_{\Phi,N,N+1}$.
Consider the isomorphic copy $\MORPH R_N$ of $R_N$ via the morphism $\MORPH$, i.e.
$\MORPH R_N = (\mathbb{R}^{n}, \{ \MORPH A_{q} \MORPH^{-1}\}_{q \in Q}, \{ \MORPH(B_j) \mid j \in J_{\Phi}\}, C\MORPH^{-1})$.  Obviously, $\MORPH R_N$ is also a $2N+1$-partial
representation of $\Psi_{\Phi}$, and it is reachable and observable.

Consider now the \LSS\  realization 
$(\Sigma_{\MORPH R_{N}},\mu_{\MORPH R_{N}})$ associated with $\MORPH R_{N}$, as
defined by Construction \ref{sect:real:const2}. 

Notice that if (\ref{part_real_lin:theo1:eq1})
holds, then (\ref{alg0:eq1}) of Algorithm \ref{alg0} has a unique solution. Indeed,
(\ref{alg0:eq1}) can be rewritten as 
\( A_qR=R_q \). Here $R,R_q \in \mathbb{R}^{n \times (Q^{\le N} \times J_{\Phi})}$ and
for each $(w,j) \in Q^{\le N} \times J_{\Phi}$, the column of $R$ indexed by $(w,j)$ is
the image $\MORPH(\COL_{w,j})$ of the column of $H_{\Phi,N,N}$ indexed by $(w,j)$;
the column of $R_{q}$ indexed by $(w,j)$ is the image $\MORPH(\COL_{wq,j})$ of the
column of $H_{\Phi,N,N+1}$ indexed by $(wq,j)$. The column space of $R$ is 
the column space of $H_{\Phi,N,N}$
by $\MORPH$, hence if (\ref{part_real_lin:theo1:eq1}) holds, then $\Rank R=n$. That is,
the solution of $A_{q}R=R_q$, if it exists, is unique.

Hence, if (\ref{part_real_lin:theo1:eq1}) holds, then the \LSS\ realization $(\Sigma_N,\mu_N)$
is uniquely defined by (\ref{alg0:eq1}),(\ref{alg0:eq-1}),(\ref{alg0:eq-2}) 
and (\ref{alg0:eq-3}), if it exists.
It is easy to see that the \LSS\ realization $(\Sigma_{\MORPH R_{N}},\mu_{\MORPH R_{N}})$ 
satisfies (\ref{alg0:eq1}),(\ref{alg0:eq-1}),(\ref{alg0:eq-2}) and (\ref{alg0:eq-3}) and hence
it coincides with $(\Sigma_N,\mu_N)$. In other words, the \LSS\ realization  $(\Sigma_N,\mu_N)$
returned by Algorithm \ref{alg0} equals the \LSS\ realization 
$(\Sigma_{\MORPH R_N},\mu_{\MORPH R_N})$ associated with the isomorphic copy 
$\MORPH R_N$ of $R_N$, i.e.
\( (\Sigma_{\MORPH R_N}, \mu_{\MORPH R_N}) = (\Sigma_N,\mu_N). \)
By Theorem \ref{part_real_lin:theo1/2}
it follows then that $(\Sigma_N,\mu_N)$ is a $2N+1$-partial realization of $\Phi$.
In addition, from Theorem \ref{sect:real:theo1}
 it follows that $(\Sigma_N,\mu_N)$ is \WR\ and observable.

\textbf{Proof of existence of a complete realization} \\
In addition, from Theorem \ref{part_real_pow:theo} it follows that if
$\Rank H_{\Phi,N,N}=\Rank H_{\Phi}$, i.e. $\Rank H_{\Psi_{\Phi},N,N}=\Rank H_{\Psi_{\Phi}}$,
then $R_{N}$ is a minimal representation of $\Psi_{\Phi}$, and hence so is $\MORPH R_N$.
But then by Theorem
\ref{sect:real:theo1}, $(\Sigma_{\MORPH R_N}, \mu_{\MORPH R_N})=(\Sigma_N,\mu_N)$ is a minimal realization of $\Psi_{\Phi}$.
Finally, assume that there
exists a realization $(\Sigma,\mu)$ of $\Phi$, such that $\dim \Sigma \le N$.
Consider the representation $R=R_{\Sigma,\mu}$ associated with $(\Sigma,\mu)$ as
defined in Construction \ref{sect:real:const1}. 
Notice that $\dim R=\dim \Sigma \le N$. Then by Theorem \ref{sect:real:theo1}, 
$R$ is a 
representation of $\Psi_{\Phi}$. Applying Theorem \ref{part_real_pow:theo} we get that
$\Rank H_{\Psi_{\Phi}}=\Rank H_{\Psi_{\Phi},N,N}$ and hence
$\Rank H_{\Phi}=\Rank H_{\Phi,N,N}$ holds, and that 
$R_N$, and hence $\MORPH R_N$,  is a  minimal representation of $\Psi_{\Phi}$.
Again, from Theorem \ref{sect:real:theo1} 
it follows then that $(\Sigma_N,\mu_N)=(\Sigma_{\MORPH R_N},\mu_{\MORPH R_N})$ is
 a minimal realization of $\Phi$.

\end{pf}
Next, we will present the proof of Theorem \ref{part_real_lin:theo1.1}

\begin{pf}[Proof of Theorem \ref{part_real_lin:theo1.1}]
 We will prove the statements of Theorem \ref{part_real_lin:theo1.1} one by one.
 Notice that if $\Phi$ satisfies (\ref{part_real_lin:theo1:eq1}), then 
 the associated family of formal power series $\Psi_{\Phi}$ satisfies
 (\ref{part_theo:eq1}).
  This allows us to use Theorem \ref{part_real_pow:theo2}
 and Theorem \ref{part_real_lin:theo1/2} to prove the statements of the theorem.
 
 \textbf{Proof of Part \ref{part_real_lin:theo1.1:res1}} \\
 From the proof of Theorem \ref{part_real_lin:theo1} it follows that 
 the \LSS\ realization $(\Sigma_N,\mu_N)$ is the realization associated with an
 isomorphic copy of the
 representation $R_N$. By the remark above the family $\Psi_{\Phi}$ associated with $\Phi$
 satisfies the condition of Theorem \ref{part_real_pow:theo2} and hence by Theorem
 \ref{part_real_pow:theo2} the representation $R_N$ is a minimal 
  $2N+1$-partial representation of
 $\Psi_{\Phi}$. Then Theorem \ref{part_real_lin:theo1/2} implies that
 $(\Sigma_N,\mu_N)$ is a $N$-partial realization of $\Phi$ and its dimension is $\dim R_{N}$.
 It is easy to see that $(\Sigma_N,\mu_N)$ is a minimal $2N+1$-partial realization of $\Phi$.
 For assume that $(\Sigma,\mu)$ is a $2N+1$-partial realization of $\Phi$. Then by
 Theorem \ref{part_real_lin:theo1/2} the associated representation $R_{\Sigma,\mu}$ is a
 $2N+1$-partial representation of $\Psi_{\Phi}$, and hence, 
 $\dim \Sigma_N=\dim R_N \le \dim R_{\Sigma,\mu}=\dim \Sigma$.
 
 \textbf{Proof of Part \ref{part_real_lin:theo1.1:res2}} \\
   If $(\Sigma,\mu)$ is a minimal $2N+1$-partial realization of $\Phi$, then from Theorem 
   \ref{part_real_lin:theo1/2} it follows that $R_{\Sigma,\mu}$ is a minimal $2N+1$-partial
   representation of $\Psi_{\Phi}$. Indeed, consider a  $2N+1$-partial representation $R$
   of $\Psi_{\Phi}$. Using Remark \ref{repr:state_space:rem}, we can replace $R$ by a suitable
   isomorphic copy, state-space of which is $\mathbb{R}^{n}$ with $n=\dim R$.
   This isomorphic copy will also be a $2N+1$-partial representation of $\Psi_{\Phi}$.
   Hence, without loss of generality we can assume that the state-space of $R$
   is $\mathbb{R}^{n}$ for $n=\dim R$ and hence Construction \ref{sect:real:const2} can be 
   applied.
   Then it follows from Theorem \ref{part_real_lin:theo1/2} that 
   $(\Sigma_{R},\mu_{R})$ is a $2N+1$-partial realization of $\Phi$. Hence,
   $\dim R_{\Sigma,\mu} = \dim \Sigma \le \dim \Sigma_R=\dim R$ by minimality of 
   $(\Sigma,\mu)$.  

   But then by Theorem \ref{part_real_pow:theo2} the representation
   $R_{\Sigma,\mu}$ has to be reachable, observable, and of dimension
   $\Rank H_{\Psi_{\Phi},N,N}=\Rank H_{\Phi,N,N}$. By Theorem \ref{sect:real:theo1},
   the latter means that $(\Sigma,\mu)$ is \WR, observable, and of dimension
   $\Rank H_{\Phi,N,N}$.
  
 \textbf{Proof of Part \ref{part_real_lin:theo1.1:res3}} \\
   Let $(\Sigma,\mu)$ and $(\hat{\Sigma},\hat{\mu})$ be two minimal
   $2N+1$-partial realizations of $\Phi$. Using an argument analogous to the one 
   presented in the proof of Part \ref{part_real_lin:theo1.1:res2}, one can show
   that the associated representations $R_{\Sigma,\mu}$ and $R_{\hat{\Sigma},\hat{\mu}}$
   are minimal $2N+1$-partial representations of the family $\Psi_{\Phi}$ of
   formal power series associated with $\Phi$.  From Theorem \ref{part_real_pow:theo2}
   it follows that $R_{\Sigma,\mu}$ and $R_{\hat{\Sigma},\hat{\mu}}$ are isomorphic.
   From Theorem \ref{sect:real:theo1} it then follows that $(\Sigma,\mu)$ and
   $(\hat{\Sigma},\hat{\mu})$ are isomorphic.

\end{pf}
We will continue with the proof of Theorem \ref{part_real_lin:theo2}
\begin{pf}[Proof of Theorem \ref{part_real_lin:theo2}]
 The proof of the theorem relies on the following observation. 
 Assume that Algorithm \ref{PartRepr} is 
 applied to $\Psi_{\Phi}$. 
  If we use the same factorization algorithm for factorizing the matrix 
  $H_{\Phi,N+1,N}=H_{\Psi_{\Phi},N+1,N}$ in both Algorithm \ref{PartReal} and 
  Algorithm \ref{PartRepr}, then the following holds.
  The equation (\ref{ComputePartialReal}) has a unique solution if and only 
  if (\ref{ComputePartialRepr}) has a unique
  solution. In addition,
  Algorithm \ref{PartReal} returns a \LSS\  realization
 if and only if Algorithm \ref{PartRepr} applied to 
 $H_{\Phi,N+1,N}=H_{\Psi_{\Phi},N+1,N}$ returns a rational representation 
 $\widetilde{R}_{N}$, and the realization $(\widetilde{\Sigma}_N,\widetilde{\mu}_N)$
 returned by Algorithm \ref{PartReal} is the \LSS\  realization
 associated with the representation $\widetilde{R}_{N}$, as defined
 in Construction \ref{sect:real:const2}. That is, 
\( 
(\widetilde{\Sigma}_{N},\widetilde{\mu}_{N})=
      (\Sigma_{\widetilde{R}_N},\mu_{\widetilde{R}_{N}})\).

\textbf{ Proof of Part \ref{part_real_lin:theo2:res1}} \\
Recall that for all $L,M >0$, the matrix $H_{\Phi,L,M}$ is identical to
$H_{\Psi_{\Phi},L,M}$, and that the Hankel-matrix $H_{\Phi}$ is identical to 
$H_{\Psi_{\Phi}}$. Hence, by applying Theorem \ref{sect:svd:theo1} to 
$H_{\Psi_{\Phi},N+1,N}$ we get that if (\ref{part_real_lin:theo1:eq1}) holds, then 
Algorithm \ref{PartRepr} returns a $2N+1$-partial representation 
$\widetilde{R}_{N}$ of $\Psi_{\Phi}$. Then from the discussion above
it follows that Algorithm \ref{PartReal} returns  
a realization $(\widetilde{\Sigma}_{N},\widetilde{\mu}_{N})$ and
\( (\widetilde{\Sigma}_{N},\widetilde{\mu}_{N})=
  (\Sigma_{\widetilde{R}_N},\mu_{\widetilde{R}_{N}}) \).
Then from
Theorem \ref{part_real_lin:theo1/2} it follows that 
 $(\widetilde{\Sigma}_{N},\widetilde{\mu}_{N})$ is a $2N+1$-partial realization of
$\Phi$. Moreover, Theorem \ref{sect:svd:theo1} implies that in this case
$\widetilde{R}_{N}$ is isomorphic to the representation $R_{N}$ from
Theorem \ref{part_real_pow:theo}. Recall from the proof of Theorem
\ref{part_real_lin:theo1} that the realization $(\Sigma_N,\mu_N)$ from 
Theorem \ref{part_real_lin:theo1} is the realization associated with
the isomorphic copy $\MORPH R_N$ of the representation $R_{N}$,
 as defined in Construction \ref{sect:real:const2}.
Since $R_N$ and $\widetilde{R}_{N}$ are isomorphic, the representations 
$\MORPH R_N$ and $\widetilde{R}_{N}$ are isomorphic as well. Then by Theorem
\ref{sect:real:theo1} the realizations $(\widetilde{\Sigma}_N,\widetilde{\mu}_N)$
and $(\Sigma_N,\mu_N)$ are isomorphic as well.

\textbf{ Proof of Part \ref{part_real_lin:theo2:res2}} \\
If $\Rank H_{\Phi}=\Rank H_{\Phi,N,N}$,  then by Theorem 
\ref{part_real_lin:theo1}, (\ref{part_real_lin:theo1:eq1}) holds and
$(\Sigma_N,\mu_N)$ is a minimal realization of $\Phi$. Similarly, if
there exists a \LSS\  realization $(\Sigma,\mu)$ of $\Phi$ such
that $\dim \Sigma \le N$, then by Theorem \ref{part_real_lin:theo1},
$\Rank H_{\Phi}=\Rank H_{\Phi,N,N}$, and (\ref{part_real_lin:theo1:eq1}) holds, 
and $(\Sigma_N,\mu_N)$ is a minimal
realization of $\Phi$. By the preceding argument, in both cases Algorithm \ref{PartReal}
returns a realization $(\widetilde{\Sigma}_N,\widetilde{\mu}_{N})$ and
$(\widetilde{\Sigma}_N,\widetilde{\mu}_{N})$ is isomorphic to $(\Sigma_N,\mu_N)$.
Hence, in both cases the realization $(\widetilde{\Sigma}_N,\widetilde{\mu}_{N})$
returned by Algorithm \ref{PartReal} is a minimal realization of $\Phi$.
\end{pf}

\section{ Discussions and Conclusions}
\label{sect:concl}
In this paper we presented basic results on partial-realization theory 
\SLSS. We also discussed the algorithmic aspects of the theory and presented two
algorithms for computing a partial \LSS\ realization. We have also shown that under
suitable conditions, the obtained partial realizations are in fact (complete) minimal
realizations of the input-output maps at hand. 

The paper represents only the first steps towards partial-realization theory of
hybrid systems. A number of problems remains open. In particular, we would like to
extend the results of the current paper to other classes of hybrid systems. Another
important research direction is the improvement of the computational complexity of the
presented algorithms. 

As we mentioned it earlier, the results of the paper can 
potentially be useful for the solution of a
number of problems, namely, for system identification,
model reduction and possibly fault detection and computer vision.
Below we will discuss these potential application domains in more detail.

 \textbf{Systems identification algorithms} \\
  As the history of development of system identification demonstrates, partial-realization 
   theory can be used to design algorithms for system identification.
  In particular, the famous subspace identification algorithms
   \cite{BilSubMoor,BookOverschee} rely heavily
  on partial-realization theory.

  In fact, 
  based on the results of the paper one can derive 
  the following crude identification algorithm. 
  Consider a  finite family $\Phi=\{f_1,f_2,\ldots,f_{d}\}$ of input-output maps.
  We would like to find a \LSS\  realization of $\Phi$.
  Fix a natural number $N$. The number $N$ may be chosen 
  either based on our belief on the dimension of a 
  potential \LSS\  realization of $\Phi$, or our ability/willingness
  to make  measurements of the system responses for switching sequences up to 
  length $2N+1$.
  \begin{algorithmic}[1]
  \STATE
        For each $f \in \Phi$, estimate the first-order time derivatives of
        $D^{(1,1,1,\ldots,1,0)} f_{0,wq}$ and 
          \( D^{(1,1,1,\ldots,1,0)} f_{e_j,q_0wq}-D^{(1,1,1,\ldots,1,0)} f_{0,q_0wq}
        \) for all sequences of discrete modes $w$ of length at most $2N+1$ and all
        discrete modes $q_0,q \in Q$.
        Use these estimates to construct the Markov-parameters of order up to $2N+1$.
        Subsequently, use the thus obtained Markov-parameters to
        construct the Hankel-matrices $H_{\Phi,N,N}$, $H_{\Phi,N+1,N}$ and
        $H_{\Phi,N,N+1}$.
  \STATE
      If $\Rank H_{\Phi,N,N}=\Rank H_{\Phi,N+1,N}=\Rank H_{\Phi,N,N+1}$,
      then compute the $2N+1$ partial realization $(\Sigma_N,\mu_{N})$ 
      either by Algorithm \ref{alg0} or Algorithm \ref{PartReal}.
  \end{algorithmic}
  It follows immediately from Theorem \ref{part_real_lin:theo1} and Theorem \ref{part_real_lin:theo2} that if
  $\Phi$ can indeed be realized by a \LSS\ , then for large
  enough $N$ the procedure above will return a minimal \LSS\ 
  realization of $\Phi$. In fact, if $N+1$ is the dimension of a potential \LSS\ realization
  of $\Phi$, then the algorithm above always returns a minimal \LSS\ realization of $\Phi$.
 
  The algorithm described above has several drawbacks. For example, it is unclear how
  to obtain in practice the derivatives needed to compute the generalized Markov-parameters.
  Another problem that the algorithm does not take into account measurement noise.
  The authors are hopeful that these shortcoming can be overcome, as similar
  problems were successfully dealt with in the case of linear systems. 
  However, we would like to remark that turning partial-realization theory to
  practically usable identification algorithms might take considerable amount of time;
  in case of linear systems it took several decades.

 \textbf{Spaces of \SLSS\  } \\
   Partial realization theory can potentially be useful for studying the geometry and
   topology of the space of \SLSS. In turn, the latter is useful for parametric
   system identification, as it helps to formalize the well-known paradigm
   of system identification: 
   "find a model of a certain class which matches the measured data best".
   In addition, the geometric insights and the resulting notions of distance
   between systems can be useful for fault detection and computer vision \cite{Martin1,Moore1,DorettoCWS03,VishwanathanSV07,PetreczkyV07}.
   For bilinear systems, the first steps were made in \cite{SontagMod}.
   For linear systems, the topology and geometry of Hankel-matrices were
   investigated in \cite{Lindquist1,Brockett1,Helmke0,Helmke1,Helmke2,Hinrichsen1,Hinrichsen2}. Note that spaces of linear systems were
   also studied using a different approach, namely, by looking at equivalence
   classes of minimal linear systems under algebraic similarity. This approach
   resulted in deep results and insights in the geometry of spaces of systems and
   important algorithms for parametric identification, see
   \cite{Hazewinkel1,Hazewinkel2,Hanzon1,Hanzon2,HanzonBook,RalfPeeters,Byrnes1,Byrnes3}.

\textbf{Model reduction } \\
  The presented results on partial-realization theory can also be used for
  model reduction. Recall from \cite{GrimmePhd,Grimme1,Antoulas1} 
the moment matching approach
  to model reduction of linear systems. The core of this approach is to
  approximate a high-order to system with a lower order one, such that the
  lower order system is a a partial realization of a finite number of
   Markov-parameter of the original system.
  For \SLSS\  the role of Markov-parameters (or moments) is played by the
  generalized Markov-parameters. Hence, one could try to extend the moment matching
  approach to \SLSS\  and replace a high-order system with a lower order one which
  is a partial-realization of a finite sequence of generalized Markov-parameters.
  Of course, a great number of issues needs to be solved. In particular, obtaining
  good numerical algorithms  will be a challenge.

\textbf{Extension to more general classes of hybrid systems} \\
  We believe that partial-realization theory for \SLSS\ 
  can be useful for obtaining realization theory, system identification and model
  reduction for piecewise-affine hybrid systems with guards. One reason for this is
  that a piecewise-affine hybrid system with guards can be viewed as a feedback 
  interconnection of a \LSS\  with an event generator.
   The former
 is just the collection of the continuous subsystems of the piecewise-affine hybrid system, the
 latter is constructed from the guards of the piecewise-affine hybrid system.
  Understanding the realization theory of one of the components of this feedback loop should
  help to understand the realization theory of piecewise-affine hybrid systems. 
  The second reason is that
  \SLSS\  form a subclass of piecewise-affine hybrid systems.
   by identifying each \SLSS\  with a
 piecewise-affine hybrid system whose guards depend only on the choice of inputs.
  Hence, any result on (partial-) realization theory of piecewise-affine hybrid systems
  should be consistent with the corresponding results for \SLSS.

\textbf{Acknowledgment}
The authors thank
Pieter Collins and Luc Habets for useful discussions and suggestions.
	       

\begin{thebibliography}{10}

\bibitem{Antoulas1}
Athanasios~C. Antoulas and Dan~C. Sorensen.
\newblock Approximation of large-scale dynamical systems: an overview.
\newblock {\em Int. J. Appl. Math. Comput. Sci.}, 11(5):1093--1121, 2001.
\newblock Numerical analysis and systems theory (Perpignan, 2000).

\bibitem{Reut:Book}
J.~Berstel and C.~Reutenauer.
\newblock {\em Rational series and Their Languages}.
\newblock Springer-Verlag, 1984.

\bibitem{Brockett1}
Roger~W. Brockett.
\newblock The geometry of the partial realization problem.
\newblock {\em Decision and Control including the 17th Symposium on Adaptive
  Processes, 1978 IEEE Conference on}, 17:1048--1052, Jan. 1978.

\bibitem{Byrnes1}
Christopher~I. Byrnes and Norman~E. Hurt.
\newblock On the moduli of linear dynamical systems.
\newblock 4:83--122, 1979.

\bibitem{Byrnes3}
Christopher~I. Byrnes and Anders Lindquist.
\newblock The stability and instability of partial realizations.
\newblock {\em Systems Control Lett.}, 2(2):99--105, 1982/83.

\bibitem{Cal:Des}
Frank~M. Callier and Charles~A. Desoer.
\newblock {\em Linear System Theory}.
\newblock Springer-Verlag, 1991.

\bibitem{MacChenBilSub}
H.~Chen and J.M. Maciejowski.
\newblock New subspace identification method for bilinea systems.
\newblock Technical Report CUED/FINFENG/TR.337, Department of Engineering,
  University of Cambridge, 2000.

\bibitem{Moore1}
K.~De~Cock and B.~De~Moor.
\newblock Subspace angles and distances between {ARMA} models.
\newblock In {\em Proc. of the Int. Symposium of Math. Theory of Networks and
  Systems}, 2000.

\bibitem{DorettoCWS03}
Gianfranco Doretto, Alessandro Chiuso, Ying~Nian Wu, and Stefano Soatto.
\newblock Dynamic textures.
\newblock {\em International Journal of Computer Vision}, 51(2):91--109, 2003.

\bibitem{AutoEilen}
Samuel Eilenberg.
\newblock {\em Automata, Languages and Machines}.
\newblock Academic Press, New York, London, 1974.

\bibitem{BilSubMoor}
W.~Favoreel, B.~De~Moor, and P.~Van~Overschee.
\newblock Subspace identification of bilinear systems subject to white inputs.
\newblock Technical Report TR 1996-531, ESAT, Katholieke Univesiteit Leuven,
  1996.

\bibitem{MFliessHank}
M.~Fliess.
\newblock Matrices de hankel.
\newblock {\em J. Math. Pures Appl.}, (23):197 -- 224, 1973.

\bibitem{MFliessFormPow}
M.~Fliess.
\newblock Functionnelles causales non lin\'eaires et ind\'etermin\'ees non
  commutatives.
\newblock {\em Bull. Soc. Math. France}, (109):2 -- 40, 1981.

\bibitem{GecsPeak}
F.~G\'ecseg and I~Pe\'ak.
\newblock {\em Algebraic theory of automata}.
\newblock Akad\'emiai Kiad\'o, Budapest, 1972.

\bibitem{BilPart}
Dieter Gollmann.
\newblock Partial realization by discrete-time internally bilinear systems: An
  algorithm.
\newblock In {\em Mathematical theory of networks and systems, Proc. int.
  Symp., Beer Sheva/Isr. 1983}, 1983.

\bibitem{Lindquist1}
William~B. Gragg and Anders Lindquist.
\newblock On the partial realization problem.
\newblock {\em Linear Algebra Appl.}, 50:277--319, 1983.

\bibitem{GrimmePhd}
E.~J. Grimme.
\newblock {\em Krylov Projection Methods}.
\newblock PhD thesis, Univ. Illinois, Urbana Champaign, 1997.

\bibitem{Grimme1}
E.~J. Grimme, D.~C. Sorensen, and P.~Van~Dooren.
\newblock Model reduction of state space systems via an implicitly restarted
  {L}anczos method.
\newblock {\em Numer. Algorithms}, 12(1-2):1--31, 1996.

\bibitem{Hanzon2}
B.~Hanzon.
\newblock Riemannian geometry on families of linear systems, the deterministic
  case.
\newblock Technical Report 88-62, Delft University of Technology, Faculty of
  Mathematics and Informatics, 1988.

\bibitem{HanzonBook}
B.~Hanzon.
\newblock {\em Identifiability, recursive identification and spaces of linear
  dynamical systems, Part I and Part II}, volume~63 of {\em CWI Tract}.
\newblock CWI, Amsterdam, 1989.

\bibitem{Hanzon1}
B.~Hanzon.
\newblock On the differentiable manifold of fixed order stable linear systems.
\newblock {\em Systems and Control Letters}, 13:345 -- 352, 1989.

\bibitem{Hazewinkel2}
M.~Hazewinkel.
\newblock Moduli and canonical forms for linear dynamical systems {III}: The
  algebraic-geometric case.
\newblock In C.~Martin and R.~Hermann, editors, {\em The 1976 AMES Research
  Center (NASA) Conference on Geometric Control Theory}, pages 229--276.
  Brookline, Mass.: Math Sci Press 1977.

\bibitem{Hazewinkel1}
M.~Hazewinkel.
\newblock Moduli and canonical forms for linear dynamical systems {II}: The
  topological case.
\newblock {\em Mathematical Systems Theory}, 10:363--385, 1977.

\bibitem{Helmke2}
U.~Helmke and D.~Hinrichsen.
\newblock Canonical forms and orbit spaces of linear systems.
\newblock {\em IMA J Math Control Info}, 3:167 -- 184, 1986.

\bibitem{Helmke1}
U.~Helmke, D.~Hinrichsen, and W.~Manthey.
\newblock A cell decomposition of the space of real {H}ankel matrices of rank
  {$\le n$} and some applications.
\newblock {\em Linear Algebra Appl.}, 122/123/124:331--355, 1989.

\bibitem{Hinrichsen1}
D.~Hinrichsen, W.~Manthey, and D.~Pr{\"a}tzel-Wolters.
\newblock The {B}ruhat decomposition of finite {H}ankel matrices.
\newblock {\em Systems Control Lett.}, 7(3):173--182, 1986.

\bibitem{Hinrichsen2}
Diederich Hinrichsen and Wilfried Manthey.
\newblock On a cell decomposition for {H}ankel matrices and rational functions.
\newblock {\em J. Reine Angew. Math.}, 451:15--50, 1994.

\bibitem{MR0245360}
B.~L. Ho and R.~E. Kalman.
\newblock Effective construction of linear state-variable models from
  input/output data.
\newblock In {\em Proc. Third Annual Allerton Conf. on Circuit and System
  Theory}, pages 449--459. Univ. Illinois, Urbana, Ill., 1965.

\bibitem{Isi:Tac}
A.~Isidori.
\newblock Direct construction of minimal bilinear realizations from nonlinear
  input-output maps.
\newblock {\em IEEE Transactions on Automatic Control}, pages 626--631, 1973.

\bibitem{Isi:Nonlin}
Alberto Isidori.
\newblock {\em Nonlinear Control Systems}.
\newblock Springer Verlag, 1989.

\bibitem{Isi:Bilin}
Alberto Isidori, Paolo D'Alessandro, and Antonio Ruberti.
\newblock Realization and structure theory of bilinear dynamical systems.
\newblock {\em SIAM J. Control}, 12(3), 1974.

\bibitem{JacobAlg1}
Nathan Jacobson.
\newblock {\em Lectures in Abstract Algebra}, volume II: linear algebra.
\newblock D. van Nostrand Company, Inc. New York, 1953.

\bibitem{MR0255260}
R.~E. Kalman, P.~L. Falb, and M.~A. Arbib.
\newblock {\em Topics in mathematical system theory}.
\newblock McGraw-Hill Book Co., New York, 1969.

\bibitem{Salomaa:Book}
W.~Kuich and A.~Salomaa.
\newblock {\em Semirings, Automata, Languages}.
\newblock Springer-Verlag, 1986.

\bibitem{D:Lib}
Daniel Liberzon.
\newblock {\em Switching in Systems and Control}.
\newblock Birkh\"auser, Boston, 2003.

\bibitem{Helmke0}
Wilfried Manthey and Uwe Helmke.
\newblock Bruhat canonical form for linear systems.
\newblock {\em Linear Algebra Appl.}, 425(2-3):261--282, 2007.

\bibitem{Martin1}
R.~Martin.
\newblock A metric for {ARMA} processes.
\newblock {\em IEEE Trans. on SIgnal Processing}, 48(4):1164--1170, 2000.

\bibitem{RalfPeeters}
Ralf Peeters.
\newblock {\em System Identification Based on Riemannian Geometry: Theory and
  Algorithms}.
\newblock PhD thesis, Free University, Amsterdam, 1994.

\bibitem{MP:HybBilinReal}
Mihaly Petreczky.
\newblock Realization theory for bilinear hybrid systems.
\newblock In {\em In Proceedings 11th IEEE International Conference on Methods
  and Models in Automation and Robotics}, 2005.

\bibitem{MP:RealBilin}
Mihaly Petreczky.
\newblock Realization theory for bilinear switched systems.
\newblock In {\em Proc. of 44th IEEE Conference on Decision and Control}, 2005.

\bibitem{MP:HybLinBilinTechReport}
Mihaly Petreczky.
\newblock Realization theory for linear and bilinear hybrid systems.
\newblock Technical Report MAS-R0502, CWI, 2005.

\bibitem{MP:BigArticle}
Mihaly Petreczky.
\newblock Realization theory of linear and bilinear switched systems: A formal
  power series approach.
\newblock Technical Report MAS-R0403, CWI, 2005.
\newblock To appear in ESAIM Control, Optimization and Calculus
  of Variations, DOI 10.1051/cocv/2010014 and DOI 10.1051/cocv/2010015.

\bibitem{MP:HybPow}
Mih\'aly Petreczky.
\newblock Hybrid formal power series and their application to realization
  theory of hybrid systems.
\newblock In {\em Proceedings 17th International Symposium on Mathematical
  Theory of Networks and Systems}, 2006.

\bibitem{MP:Phd}
Mihaly Petreczky.
\newblock {\em Realization Theory of Hybrid Systems}.
\newblock PhD thesis, Vrije Universiteit, Amsterdam, 2006.
\newblock Available at \verb'http://www.cwi.nl/~mpetrec'.

\bibitem{MP:RealForm}
Mihaly Petreczky.
\newblock Realization theory for linear switched systems: Formal power series
  approach.
\newblock {\em Systems and Control Letters}, 56(9 -- 10):588--595, 2007.

\bibitem{PetreczkyV07}
Mih{\'a}ly Petreczky and Ren{\'e} Vidal.
\newblock Metrics and topology for nonlinear and hybrid systems.
\newblock In {\em Hybrid Systems: Computation and Control}, pages 459--472,
  2007.

\bibitem{MPRV:JumpMarkov}
Mihaly Petreczky and Rene Vidal.
\newblock Realization theory of stochastic jump-markov linear systems.
\newblock In {\em Proceedings 46th IEEE Conference on Decision and Control},
  2007.

\bibitem{SontagMod}
D.~Eduardo Sontag.
\newblock A remark on bilinear systems and moduli spaces of instantons.
\newblock {\em Systems and Control Letters}, 9(5):361--367, 1987.

\bibitem{Son:Resp}
Eduardo~D. Sontag.
\newblock {\em Polynomial Response Maps}, volume~13 of {\em Lecture Notes in
  Control and Information Sciences}.
\newblock Springer Verlag, 1979.

\bibitem{Son:Real}
Eduardo~D. Sontag.
\newblock Realization theory of discrete-time nonlinear systems: Part {I} --
  the bounded case.
\newblock {\em IEEE Transaction on Circuits and Systems}, CAS-26(4), April
  1979.

\bibitem{Sun:Book}
Zhendong Sun and Shuzhi~S. Ge.
\newblock {\em Switched linear systems : control and design}.
\newblock Springer, London, 2005.

\bibitem{TetherPartLin}
A.~Tether.
\newblock Construction of minimal linear state-variable models from finite
  input-output data.
\newblock {\em IEEE Transactions Automatic Control}, 15(4):427-- 436, 1970.

\bibitem{BookOverschee}
P.~van Overschee and B.~De Moor.
\newblock {\em Subspace Identification for Linear Systems}.
\newblock Kluwer Academic Publishers, 1996.

\bibitem{VishwanathanSV07}
S.~V.~N. Vishwanathan, Alexander~J. Smola, and Ren{\'e} Vidal.
\newblock Binet-cauchy kernels on dynamical systems and its application to the
  analysis of dynamic scenes.
\newblock {\em International Journal of Computer Vision}, 73(1):95--119, 2007.

\end{thebibliography}

\appendix
 \section{ Formal Power Series } 
 \label{sect:pow} 
The section recalls basic results on formal power series. 
The material of this section is based on the classical
theory of formal power series, see
\cite{Reut:Book,Salomaa:Book}.  However, a number of
concepts and results are extensions of the standard ones to families of
formal power series.
The outline of the section is the following. In Subsection \ref{sect:pow:def}
we present the basic concepts for formal power series.
In Subsection \ref{sect:pow:rat_repr} we present the notion of a rational
representation of a family of formal power series. 
In Subsection \ref{sect:pow:main} we state the main result on existence
and minimality of a rational representation of a family of formal power series.

\subsection{ Formal Power Series: Definition and Basic Concepts }
\label{sect:pow:def}

 Let $X$ be a finite set, which we will refer to as the alphabet.
 Recall from Section \ref{sect:prelim:lang} the notion of a
 finite word over an alphabet and the related notation.
  A \emph{formal power series} $S$ with coefficients in 
  $\mathbb{R}^{p}$ is a map
  \[ S: X^{*} \rightarrow \mathbb{R}^{p} \]
  There are many ways to give an intuition for the definition
  of a formal power series. For the purposes of this paper,
  the most suitable one is to think of a formal power series
  as an output of a machine defined as follows. The machine
  reads symbols belonging to $X$ from its input tape and writes elements
  of $\mathbb{R}^{p}$ onto its output tape. 
   We denote by
   $\mathbb{R}^{p} \ll   X^{*}\gg$ the set of all formal power series
   with coefficients in $\mathbb{R}^{p}$.
  The set of all formal power series over $X$ with
  coefficients in $\mathbb{R}^{p}$ forms a vector space
  with respect to point-wise addition and multiplication.
  That is, if $\alpha, \beta \in \mathbb{R}$ are two
  scalars and $S,T \in \mathbb{R}^{p}\ll X^{*} \gg$ are
  two formal power series, then the linear combination
  $\alpha S+ \beta T$ is defined as the formal power
  series assigning each word $w \in X^{*}$ the value
  $\alpha S(w)+\beta T(w)$. 

   In the sequel we will mostly be interested in \emph{families
   of formal power series}. 
   \begin{Definition}[Familly of formal power series]
   Let $J$ be an arbitrary 
   (possibly infinite) set. 
   A \emph{family of formal power series in 
    $\mathbb{R}^{p}\ll X^{*}\gg$ indexed by $J$} is simply a 
    collection 
   $\Psi=\{ S_{j} \in \mathbb{R}^{p}\ll X^{*} \gg \mid j \in J\}$
   of formal power series from $\mathbb{R}^{p}\ll X^{*} \gg$ indexed
   by elements of $J$.
   \end{Definition}
    Notice that the definition above does not require 
    $S_{j}$, $j \in J$ to be all
    distinct formal power series, i.e. $S_{l}=S_{j}$ for some indices
    $j,l \in J$, $j \ne l$ is allowed.
   One can think of a family of formal power series as a
   family of input-output maps of the machine 
   (discrete-time system) 
   described above, realized from a set of initial states
   indexed by elements of $J$. We would like to
   point out that the notion of (family of) formal power series
   is a purely formal one, the interpretation of them as
   input-output maps is just one of the many 
   possible interpretations.

\subsection{ Rational Representations and 
              Rational Formal Power Series}
\label{sect:pow:rat_repr}	      
   Above we have defined the notion of a formal power series
   and a family of formal power series and we have related
   these notions to input-output maps of some systems.
   Below we will recall the notion of a \emph{rational
   representation}, which can be thought of as a special
   subclass of these systems.  We will also define morphisms
   between rational representations along with a 
   notion of observability and reachability.

   \begin{Definition}[Representation]
   \label{sect:pow:rat_rep:def0}
   Let $J$ be an arbitrary set and let $\mathbb{N} \ni p >0$.  A 
   \emph{rational representation of
   type $p$-$J$ over the alphabet $X$}
   is a tuple
   \begin{equation}
   \label{repr:def0}
    R=(\mathcal{X},\{A_{\sigma}\}_{\sigma \in X},B,C)
   \end{equation}
   where
   \begin{itemize}
   \item
    The space $\mathcal{X}$ is a finite dimensional 
   vector space over $\mathbb{R}$, called \emph{state-space} of $R$,
   \item
    For each letter
   $\sigma \in X$, 
   $A_{\sigma}:\mathcal{X} \rightarrow \mathcal{X}$ is a 
   linear map, referred to as the \emph{state-transition map}
   \item
   The map $C:\mathcal{X} \rightarrow \mathbb{R}^{p}$
   is a linear map, referred to as the \emph{readout map},
  \item
   The family $B=\{ B_{j} \in \mathcal{X} \mid j \in J\}$
   is a collection of (not necessarily distinct) elements of $\mathcal{X}$ indexed
   by $J$. 
  \end{itemize}
  If $p$ and $J$ are clear from the context, then we will
   refer to $R$ simply as a \emph{rational representation}.
   
   The dimension $\dim \mathcal{X}$ of the state-space is called the 
   \emph{dimension} of the representation $R$ and
   it is denoted by $\dim R$. 
  \end{Definition}
  \begin{Remark}
   Notice that if a basis of $\mathcal{X}$ is fixed and
   $n=\dim \mathcal{X}$, then the state-transition maps $A_{\sigma}$, $\sigma \in X$ and
   the readout map $C$ can be identified
   with their matrix representations in this basis, and
   for each $j \in J$, $B_{j}$ can be identified with the
   $\mathbb{R}^{n}$ column vector of its coordinates. 
   If $\mathcal{X}=\mathbb{R}^{n}$, then we will identify
   the linear maps $A_{\sigma}$, $\sigma \in X$ and
   $C$ with their matrix representations in the standard
   orthogonal basis of $\mathbb{R}^{n}$. In this case 
   we will call them the \emph{state-transition matrices} and the
   \emph{readout matrix} respectively.
   \end{Remark}
  \begin{Definition}[Rational family of formal power series]
  \label{sect:pow:rat_rep:def-1}
  Let 
  $\Psi=\{ S_{j} \in \mathbb{R}^{p}\ll X^{*} \gg  \mid j \in J\}$
  be a family of formal power series indexed by $J$. The
  representation $R$ from (\ref{repr:def0}) is said to be
  a \emph{representation of $\Psi$}, if for each index $j \in J$,
  for all sequences $\sigma_1,\sigma_2,\ldots, \sigma_k \in X$, $k \ge 0$,
   \begin{equation}
   \label{repr:def1}
      S_{j}(\sigma_{1}\sigma_{2} \cdots \sigma_{k})=
      CA_{\sigma_{k}}A_{\sigma_{k-1}} \cdots A_{\sigma_{1}}B_{j}.
   \end{equation}
  We will say that a family of formal power series
  $\Psi=\{ S_{j} \in \mathbb{R}^{p}\ll X^{*} \gg \mid j \in J\}$
  is \emph{rational}, if there exists a representation
  $R$ such that $R$ is a representation of $\Psi$.
  \end{Definition}
 \begin{Notation}
 \label{repr:not1}
   The following notation will greatly simplify the expressions
   used for rational representations.
   Let $A_{\sigma}: \mathcal{X} \rightarrow \mathcal{X},\sigma \in X$ 
   be linear maps and let $w=\sigma_1\sigma_2 \cdots \sigma_{k} \in X^{*}$,
   $\sigma_1, \sigma_2,\cdots \sigma_k \in X$, $k \ge 0$,  be a word over 
   $X$. Then $A_{w}$ denotes the 
   composition of the linear maps $A_{\sigma_1}, A_{\sigma_{2}}, \ldots A_{\sigma_k}$ in that order, that is
   \begin{equation} 
   \label{repr:eq2}
   A_{w}=A_{\sigma_{k}}A_{\sigma_{k-1}} \cdots A_{\sigma_{1}}
   \end{equation}
   If $w=\epsilon$ is the empty words, then $A_{\epsilon}$ is taken to be the identity map.
\end{Notation}
   With the notation above, (\ref{repr:def1}) can be 
   rewritten as 
   $S_{j}(w)=CA_{w}B_{j}$ for all $w \in X^{*}$, $j \in J$.
   Notice that if $\mathcal{X}=\mathbb{R}^{n}$ for some
   $\mathbb{N} \ni n > 0$, and hence the state-transition maps $A_{\sigma}$
   can be viewed as matrices, then (\ref{repr:eq2})
   defines a notation for products of the state-transition
   matrices taken along the word $w$.
\begin{Definition}[Minimality]
 A representation $R_{min}$ of $\Psi$ is called \emph{minimal}, 
 if 
 for each representation $R$ of $\Psi$, 
 \( \dim R_{min} \le \dim R \), i.e. $R_{min}$ is a
 rational representation of $\Psi$ with the smallest possible
 state-space dimension.
\end{Definition} 
 We will continue with presenting the
 notions of \emph{observability} and \emph{reachability} for rational representations.
 Define the subspaces $W_{R}$ and $O_{R}$ of
 $\mathcal{X}$ by
  \begin{eqnarray}
  \label{sect:pow:reachobs:eq1.1}
  W_{R} & = & \SPAN\{ A_{w}B_{j} \in \mathcal{X} \mid w \in X^{*}, j \in J\}    \\
  \label{sect:pow:reachobs:eq2.1}
     O_{R}& = & \bigcap_{w \in X^{*}} \ker CA_{w}
  \end{eqnarray}
  That is, the subspace $W_{R}$ is the linear span of
  the elements of the state-space of the form $A_{w}B_{j}$, where
  $w$ runs through all words over $X$ and $j$ runs through all
  the indices in $J$. The space $O_{R}$ is the intersection
  of the null-spaces (kernels) of all the linear maps
  $CA_{w}:\mathcal{X} \rightarrow \mathbb{R}^{p}$, where
  $w$ runs through the set of all words $X^{*}$.
 \begin{Definition}[Reachability]
 \label{reach:rep:def}
  We will say that the representation $R$ is \emph{reachable}
  if $\dim W_{R}=\dim R$.
  The subspace $W_{R}$ will be referred to as the \emph{reachability subspace of $R$}.
 \end{Definition}
\begin{Definition}[Observability]
 \label{obs:rep:def}
 We will say that $R$ is
  \emph{observable} if $O_{R}=\{0\}$, i.e. $O_{R}$ consists of the zero element only.
  The subspace $O_{R}$ will be referred to as the \emph{observability subspace of $R$}.
\end{Definition}
  Next, we define the notion of morphism between rational
  representations. 
\begin{Definition}[Representation morphism]
  Let $R=(\mathcal{X},\{ A_{\sigma} \}_{\sigma \in X},B,C)$,
  \( \widetilde{R}=(\widetilde{\mathcal{X}},\{ \widetilde{A}_{\sigma} \}_{\sigma \in X}, \widetilde{B}, \widetilde{C}) \) be two 
  $p$-$J$ rational representations. 
  A linear map
  $\RMORPH: \mathcal{X} \rightarrow \widetilde{\mathcal{X}}$ is
  called a \emph{representation morphism} from $R$ to 
  $\widetilde{R}$ and is denoted by $\RMORPH:R \rightarrow \widetilde{R}$
  if $\RMORPH$ commutes with $A_{\sigma}$, $B_{j}$ and
  $C$ for all $j \in J$, $\sigma \in X$, that is, if
  the following equalities hold
  \begin{equation}
     \begin{array}{rcl}
     \label{repr:morph:eq2}
     \RMORPH A_{\sigma}=\widetilde{A}_{\sigma}\RMORPH, \forall \sigma \in X, 
     & \RMORPH B_{j}=\widetilde{B}_{j}, \forall j \in J,  &
     C=\widetilde{C}\RMORPH 
   \end{array}
  \end{equation}
    The representation morphism $\RMORPH$ is called \emph{surjective, injective, isomorphism}  
    if $\RMORPH$ is a surjective, injective or isomorphism respectively if viewed as a linear map.
\end{Definition}
  \begin{Remark}
  \label{repr:state_space:rem}
If $R=(\mathcal{X},\{A_{\sigma}\}_{\sigma \in X}, B,C)$ is a
representation of $\Psi$, then for any vector space isomorphism
$\RMORPH:\mathcal{X} \rightarrow \mathbb{\mathbb{R}}^{n}$, $n=\dim R$, the
tuple
\begin{equation} 
\label{repr:euclid:eq1}
  \RMORPH R=(\mathbb{\mathbb{R}}^{n}, \{ \RMORPH A_{\sigma}\RMORPH^{-1} \}_{\sigma \in X}, \RMORPH B, C\RMORPH^{-1}) 
\end{equation}  
where $\RMORPH B=\{ \RMORPH B_{j} \in \mathbb{R}^{n} \mid j \in J\}$
  is also a representation of $\Psi$.
Moreover, for all $\sigma \in X$, $\RMORPH A_{\sigma}\RMORPH^{-1}$ 
can be naturally viewed as a matrix by taking its matrix
representation with respect to the natural basis of
$\mathbb{R}^{n}$. Similarly, by taking
matrix and vector representations
of $C\RMORPH^{-1}$ and $\RMORPH B_{j}$, $j \in J$ in the natural
basis of $\mathbb{R}^{n}$, we can view $C\RMORPH^{-1}$ and $\RMORPH B_{j}$, $j \in J$ as a $p \times n$ matrix and $n \times 1$ vector respectively.
Moreover, $\RMORPH:R \rightarrow \RMORPH R$ is a representation
isomorphism. That is, we can always replace a representation of
$\Psi$ with an isomorphic representation, state-space of which
is $\mathbb{R}^{n}$ for some $n$, and the parameters of
which are matrices and real vectors, as opposed to linear maps
and elements of abstract vector spaces. Moreover,
isomorphisms clearly preserve such properties as observability,
reachability, and minimality.
\end{Remark}

\subsection{ Existence and Minimality of Rational
             Representations: Main Results}
\label{sect:pow:main}
 The purpose of the section is to state the main
 results on existence and minimality of representations
 of families of rational formal power series. 

   We will start by stating the main result on
   the existence of a rational representation. However, in 
   order to state the main theorem, we need to define the
   concept of the \emph{Hankel matrix} of a family of
   formal power series.
    Let $\Psi=\{ S_{j} \in \mathbb{R}^{p}\ll X^{*} \gg \mid j \in J\}$ be a family of formal power series.
    Recall the notation of Section \ref{sect:prelim:not_matrix}.
    \begin{Construction}[Hankel-matrix] 
    \label{sect:pow:hank:def0}
    Define the \emph{Hankel-matrix} of $\Psi$ as the 
    following infinite matrix $H_{\Psi}$.
    The rows of $H_{\Psi}$ are indexed
    by pairs $(v,i)$ where $v \in X^{*}$ is an arbitrary word and
    $i=1,\ldots, p$. The columns of $H_{\Psi}$ are indexed
    by pairs $(w,j)$ where $w \in X^{*}$ is a word over
    $X$ and $j \in J$ runs through the elements of $J$.
    The set of all real infinite matrices with the indexing 
    the rows and columns described above will be denoted by
    $\mathbb{R}^{(X^{*} \times I) \times (X^{*} \times J)}$
    where $I=\{1,\ldots, p\}$.
    Hence,  we can write 
    $H_{\Psi} \in \mathbb{R}^{(X^{*} \times I) \times (X^{*}\times J)}$.
    The entry of $H_{\Psi}$ lying on the intersection of
    the row indexed by $(v,i)$ and the column indexed
    by $(w,j)$ is defined as
    \begin{equation}
    \label{repr:eq3}
      (H_{\Psi})_{(v,i)(w,j)}= (S_{j}(wv))_{i}
    \end{equation}
    where $(S_{j}(wv))_{i}$ denotes the $ith$ entry of
    the column vector $S_{j}(wv) \in \mathbb{R}^{p}$.
 \end{Construction}   
    According to the convention adopted in Section \ref{sect:prelim:not_matrix},
    we define the \emph{rank of $H_{\Psi}$, denoted by
    $\Rank H_{\Psi}$},  as the dimension of the vector space
    spanned by the columns of $H_{\Psi}$.
    With the notation above the following theorem holds. 
 \begin{Theorem}[Existence of a representation,\cite{MP:Phd,MP:BigArticle}]
 \label{sect:form:theo1}
  The family $\Psi$ is rational, i.e. admits a rational
  representation, if and only if $\Rank H_{\Psi} <+\infty$,
  i.e. the rank of the  Hankel-matrix $H_{\Psi}$ is finite.
\end{Theorem}
 The proof of the above theorem can be found in \cite{MP:Phd,MP:BigArticle}.

 As the next step we will present below the main result on
 minimality of rational representations.
 \begin{Theorem}[Minimal representation, \cite{MP:Phd,MP:BigArticle}]
 \label{sect:form:theo3}
 Let $\Psi=\{ S_{j} \in \mathbb{R}^{p} \ll X^{*} \gg \mid j \in J\}$
 be a family of formal power series.
 The following are equivalent.
 \begin{itemize}
 \item[(i)]
    $R_{min}$
    is a minimal representation of $\Psi$.
 \item[(ii)]
  $R_{min}$ is reachable and observable.
 \item[(iii)]
    If $R$ is a reachable  representation of $\Psi$,
    then there exists a surjective representation morphism
    $\RMORPH:R \rightarrow R_{min}$.
  \item[(iv)]
    $\Rank H_{\Psi}=\dim R_{min}$.
 \end{itemize}   
 In addition, all minimal representations of $\Psi$ are isomorphic.
\end{Theorem}
 The proof of the above theorem can be found in \cite{MP:Phd,MP:BigArticle}.
\begin{Remark}[Related work]
 The counterpart of the above two theorems, i.e. Theorem \ref{sect:form:theo1} and
 \ref{sect:form:theo3},  for a single formal
 power series is a classical result, see
 \cite{MFliessHank,Reut:Book,Salomaa:Book,Son:Resp,Son:Real}.
\end{Remark}

\section{Partial-realization theory of formal power series}
\label{sect:pow:part}
Below we will formulate and solve the counterpart of the partial realization problem for 
families of formal power series. We use the obtained result to derive a 
solution to the partial realization problem for linear switched systems formulated
in Problem \ref{prob:part:real}. We start with defining the notion of
a $N$-partial representation.
\begin{Definition}[$N$-partial representation]
 Let $\Psi=\{ S_{j} \in \mathbb{R}^{p}\ll X^{*} \gg \mid j \in J\}$
be a family of formal power series indexed by $J$.
A $J$-$p$ representation
$R=(\mathcal{X}, \{ A_{x} \}_{x \in X}, C,B)$, with $B=\{ B_{j} \in \mathcal{X} \mid j \in J\}$
is said to be an
\emph{$N$-partial representation of $\Psi$} if for each index $j \in J$
and each word $w \in X^{*}$ of length at most $N$, i.e. $|w| \le N$,
\[ S_{j}(w)=CA_{w}B_{j} \]
\end{Definition}
That is, if $R$ is a $N$-partial representation of $\Psi$, then $R$ recreates the
values of the elements of $\Psi$ for all the words of length at most $N$.
Now we are ready to formulate the partial-realization problem for formal power series.
\begin{Problem}[Partial-realization problem for formal power series]
 Let $\Psi$ be a family of formal power series indexed by $J$. 
 \begin{itemize}
 \item
     Find conditions for existence of a $N$-partial representation of $\Psi$. Formulate
     an algorithm for computing a $N$-partial representation of $\Psi$.
 \item
     Characterize minimal $N$-partial representations of $\Psi$, their existence and
     uniqueness.
 \item
     Find conditions under which the $N$-partial representation above becomes a 
     rational representation of $\Psi$ in the sense of Definition \ref{sect:pow:rat_rep:def0}.
 \end{itemize}
\end{Problem}
 We will devote the rest of the section to solving the problem described above.
 The outline of the section is the following. 
 Subsection \ref{part:form:main_results} 
 presents the main results on partial-realization 
 theory of families of formal power series, along with a Kalman-Ho-like algorithm
 for computing a minimal partial representation.
 Subsection \ref{app:proof} presents the proof of the results presented in Subsection
  \ref{part:form:main_results}.
 Throughout the section, $\Psi$ will denote the family of formal power series
 $\Psi=\{S_j \in \mathbb{R}^{p}\ll X^{*} \gg \mid j \in J\}$.

\subsection{Main results on partial-realization theory of  formal power series}
\label{part:form:main_results}
 The goal of the section is to present the main results on partial-realization
 theory of families of formal power series. 
 Note that partial-realization theory of a single formal power series is
 more or less  equivalent to partial-realization theory of bilinear and
 state-affine systems. The latter was already investigated in \cite{Isi:Tac,BilPart,Isi:Bilin,Son:Resp,Son:Real,BilSubMoor,MacChenBilSub}. The results to be presented below represent an
 extension of the results of
 \cite{Isi:Tac,BilPart,Isi:Bilin,Son:Resp,Son:Real,BilSubMoor,MacChenBilSub}.
 In addition,
 here we state partial-realization theory directly for formal power series.
 This is in contrast to
 \cite{Isi:Tac,BilPart,Isi:Bilin,Son:Resp,Son:Real,BilSubMoor,MacChenBilSub}, where
 the partial-realization problem and solution were stated in terms of state-affine or
 bilinear systems.

 The outline of the section is the
 following.
 Subsection \ref{part:form:theo}
 presents the statement of the main results on partial realization theory for 
 formal power series.
 Subsection \ref{part:form:svd} 
 presents an algorithm for computing a partial representation,
 which is similar to the well-known Kalman-Ho algorithm for the linear systems.

\subsubsection{ Partial realization theory }
\label{part:form:theo}

We will start with defining the upper-left block matrix of $H_{\Psi}$
which is indexed by words over $X$ of finite length.
More precisely, fix natural numbers $M,K > 0$ and define the following
sets
\begin{equation}
 \begin{split}
 I_{M}&=\{(v,i) \mid v \in X^{*}, |v| \le M,i=1,\ldots,p \} \\
 J_{K}&=\{ (w,j) \mid j \in J, w \in X^{*}, |w| \le K \}
\end{split}
\end{equation}
 Intuitively, $I_M$ is a subset of the row indices of the Hankel-matrix $H_{\Psi}$,
 made up of indices of the form $(v,i)$, where $v$ runs through all the words of length
 at most $M$. Similarly, the set $J_K$ is a subset of column indices of $H_{\Psi}$,
 made up of column indices $(w,j)$ where $w$ runs through all the words of length at
 most $K$.
 The set $J_K$ is the set of columns indices of the matrix to be defined, and
 $I_{M}$ is the set of row indices of this matrix.

\begin{Definition}[Sub-matrices of the Hankel-matrix $H_{\Psi}$]
\label{part_real_pow:submatrix}
Define the matrix $H_{\Psi,M,K} \in \mathbb{R}^{I_{M} \times J_{K}}$ by
\[ (H_{\Psi,M,K})_{(v,i),(w,j)}=(H_{\Psi})_{(v,i),(w,j)}=(S_{j}(wv))_{i}
\]
for all $(v,i) \in I_{M}$ and $(w,j) \in J_K$.
\end{Definition}
That is, $H_{\Psi,M,K}$ is the left upper
$I_{M} \times J_{K}$ block matrix of $H_{\Psi}$. Notice that if
$J$ is finite, then $|J_{K}| < +\infty$, that is, 
$H_{\Psi,M,K}$ is a \emph{finite matrix}.

It turns out that under certain circumstances partial representations
not only exist but they also yield a minimal representation of the
whole family of formal power series. Moreover, 
such partial representations can be constructed from finite data.

\begin{Theorem}[Existence of partial representation \\]
\label{part_real_pow:theo}
With the notation above the following holds.
\begin{enumerate}
\item
\label{part_theo:part1}
      If for some $N > 0$,
      \begin{equation}
      \label{part_theo:eq1}
        \Rank H_{\Psi,N,N}=\Rank H_{\Psi,N,N+1}=\Rank H_{\Psi,N+1,N} 
      \end{equation}
      then there exists a $2N+1$-partial representation
      $R_{N}$ of $\Psi$, of the form
      \begin{equation}
      \label{part_theo:eq3}
       R_{N}=(\IM H_{\Psi,N,N+1}, \{ A_{\sigma} \}_{\sigma \in X}, B,C)
      \end{equation}
      where the parameters of $R_{N}$ are defined as follows.
      For each word $w \in X^{*}$, $|w| \le N+1$, and each index $j \in J$ denote by
      $(H_{\Psi,N,N+1})_{.,(w,j)}$ the column of $H_{\Psi,N,N+1}$ indexed by
      $(w,j)$. With this notation,
      \begin{itemize}
      \item
         $\IM H_{\Psi,N,N+1}$ denotes the linear space spanned by the columns of 
         $H_{\Psi,N,N+1}$.
      \item
          For each $\sigma \in X$, the linear map $A_{\sigma}:\IM H_{\Psi,N,N+1} \rightarrow \IM H_{\Psi,N,N+1}$ has the
          property that for each word $w \in X^{*}$, $|w| \le N$ and each index $j \in J$
          \begin{equation}
          \label{part_theo:eq4}
           A_{\sigma}((H_{\Psi,N,N+1})_{.,(w,j)})=
               (H_{\Psi,N,N+1})_{.,(w \sigma, j)}
          \end{equation}
          i.e. $A_{\sigma}$ maps the column indexed by $(w,j) \in J_N$ to the column
          indexed by $(w\sigma,j)$.
      \item
          For each $w \in X^{*}$, $|w| \le N+1$ and $j \in J$, the linear
          map $C:\IM H_{\Psi,N,N+1} \rightarrow \mathbb{R}^{p}$ satisfies,
          \begin{equation}
          \label{part_theo:eq5}
          \begin{split}
           & C((H_{\Psi,N,N+1})_{.,(w,j)}) = \\
           & \begin{bmatrix}
            (H_{\Psi,N,N+1})_{(\epsilon,1),(w,j)}, & 
            (H_{\Psi,N,N+1})_{(\epsilon,2),(w,j)}, & \cdots, &
            (H_{\Psi,N,N+1})_{(\epsilon,p),(w,j)}
          \end{bmatrix}^{T}
          \end{split}
          \end{equation}
          That is, $C$ maps each column to the vector in $\mathbb{R}^{p}$ formed
          by the first $p$ entries associated with $v=\epsilon$
          of the column indexed by $(w,j)$.
      \item
         The set $B=\{ B_{j} \in \IM H_{\Psi,N,N+1} \mid j \in J \}$ is defined by
         \begin{equation}
         \label{part_theo:eq6} 
           B_{j}=(H_{\Psi,N,N+1})_{.,(\epsilon,j)} \mbox{ for all } j \in J 
         \end{equation}
         i.e. $B_{j}$ is simply the column of $H_{\Psi,N,N+1}$ indexed by
         $(\epsilon,j)$.
      \end{itemize}
      In addition, the representation $R_{N}$ is reachable and observable.
\item
\label{part_theo:part2}
   If for some $N >0$
   \begin{equation}   
   \label{part_theo:eq2}
     \Rank H_{\Psi,N,N}=\Rank H_{\Psi}
   \end{equation}
   then (\ref{part_theo:eq1}) holds 
   and the representation $R_{N}$ from (\ref{part_theo:eq3})
   is a minimal representation of $\Psi$.
\item
\label{part_theo:part3}
   If  $R$ is a representation of $\Psi$, $\dim R \le N+1$, then 
   (\ref{part_theo:eq2}) holds for $N$ and
   the representation $R_{N}$ from (\ref{part_theo:eq3}) exists and
   it is a minimal representation of $\Psi$.
\end{enumerate}   
\end{Theorem}
 The proof of Theorem \ref{part_real_pow:theo} is presented in Appendix \ref{app:proof}.
 Using the results of Theorem \ref{part_real_pow:theo} we can state the following
 characterization of minimal partial representations of $\Psi$.

\begin{Theorem}[Minimal partial representation]
\label{part_real_pow:theo2}
 With the notation of Theorem \ref{part_real_pow:theo}, if $\Psi$ satisfies
 (\ref{part_theo:eq1}), then the following holds.
 \begin{enumerate}
 \item
  \label{part_real_pow:theo2:res1}
   A minimal $2N+1$ partial representation of $\Psi$ exists, in fact, the representation $R_N$
   of Theorem \ref{part_real_pow:theo} is a minimal $2N+1$ partial representation of $\Psi$.
 \item
  \label{part_real_pow:theo2:res2}
   Any minimal $2N+1$ partial representation of $\Psi$ is reachable and observable and
   it is of dimension $\Rank H_{\Psi,N,N}$. 
 \item
 \label{part_real_pow:theo2:res3}
  All minimal $2N+1$ partial representations of $\Psi$ are isomorphic.
 \end{enumerate}
\end{Theorem}
 The proof of Theorem \ref{part_real_pow:theo2} is presented in Appendix \ref{app:proof}.
 \begin{Remark}
  The reader might wonder why we speak of  $2N+1$-partial representations in 
  Theorem \ref{part_real_pow:theo} and \ref{part_real_pow:theo2}. The reason behind it
  is that the finite Hankel-matrix $H_{\Psi,N,N+1}$ is
  formed by values of the formal power series from $\Psi$ for words of length at most
  $2N+1$. 
   We would like the representation obtained
   from $H_{\Psi,N,N+1}$ to recreate at least the entries of
   the matrix $H_{\Psi,N,N+1}$. But this means precisely that the
   representation obtained from $H_{\Psi,N,N+1}$ should be
   a $2N+1$-partial representation of $\Psi$. 
 \end{Remark}

\subsubsection{ Partial realization algorithm }
\label{part:form:svd}

 In this section we present an algorithm, described in Algorithm \ref{PartRepr}, which
 computes a partial representation by factorizing the Hankel-matrix.
 The technique of Hankel-matrix factorization has been used in 
 realization 
 theory and systems identification for several decades. 
 It forms the theoretical basis of algorithms for subspace
 identification, see for example
 \cite{BilSubMoor,MacChenBilSub}. Throughout the section, $\Psi$ stands for the family of
 formal power series $\Psi=\{S_j \in \mathbb{R}^{p}\ll X^{*} \gg \mid j \in J\}$.

\begin{algorithm}
\caption{\textbf{ComputePartialRepresentation($H_{\Psi,N+1,N}$)}}
\label{PartRepr}
\begin{algorithmic}[1]
 \STATE
\label{alg:svd:step1}
 Compute a decomposition of $H_{\Psi,N+1,N}$
    \[ H_{\Psi,N+1,N}=OR \]
    $O \in \mathbb{R}^{I_{N+1} \times r}$, 
    $R \in \mathbb{R}^{r \times J_{N}}$,
    $\Rank R=\Rank O=\Rank H_{\Psi,N+1,N}=r$
 
 \STATE
  Define the matrix $\widetilde{C} \in \mathbb{R}^{p \times r}$ by
 \begin{equation} 
  \label{ComputePartialRepr-1}
    \widetilde{C}=\begin{bmatrix} O_{(\epsilon,1),.}^{T}, & O_{(\epsilon,2),.}^{T}, & 
    \cdots & O_{(\epsilon,p),.}^{T} \end{bmatrix}^{T} 
  \end{equation}
   where $O_{k,.}$ denotes the row of $O$ indexed by $k \in I_{N+1}$.

\STATE
  Define the family of vectors
 \( \widetilde{B}=\{ \widetilde{B}_{j} \in \mathbb{R}^{r} \mid j \in J \} \), 
  where for each $j \in J$, 
  \begin{equation} 
  \label{ComputePartialRepr-2}
    \widetilde{B}_{j}=R_{.,(\epsilon,j)}
  \end{equation}
  where $R_{.,(\epsilon,j)}$ stands for the column of $R$ indexed by $(\epsilon,j)$.
\STATE  
For each $\sigma \in X$ let $\widetilde{A}_{\sigma} \in \mathbb{R}^{r \times r}$ 
be the solution of
\begin{equation}
\label{ComputePartialRepr}
\bar{\Gamma}\widetilde{A}_{\sigma}=\bar{\Gamma}_{\sigma} 
\end{equation}
where
$\bar{\Gamma}, \bar{\Gamma}_{\sigma} \in \mathbb{R}^{I_{N} \times r}$ are matrices of the
form
\begin{equation*} 
\begin{split}
 \bar{\Gamma}_{(u,i),j}=O_{(u,i),j} \mbox{ and } 
  (\bar{\Gamma}_{\sigma})_{(u,i),j}=O_{(\sigma u,i),j} 
\end{split}
\end{equation*}
for all  \( (u,i) \in I_{N}, j=1,2,\ldots,r \).
\STATE
  If there no unique solution to (\ref{ComputePartialRepr}) then return
  $NoRepresentation$. Otherwise return
   \[ \widetilde{R}_{N}=(\mathbb{R}^{r}, \{ \widetilde{A}_{x} \}_{x \in X}, \widetilde{B},\widetilde{C}) \]
\end{algorithmic}
\end{algorithm}
Algorithm \ref{PartRepr} represents an algorithm based on matrix factorization 
of $H_{\Psi,N+1,N}$. In addition, if $N$
is large enough, Algorithm \ref{PartRepr} in fact yields a minimal representation of
$\Psi$.
Algorithm \ref{PartRepr} above may return two different types of data. It returns
a rational representation if 
(\ref{ComputePartialRepr}) has a unique solution, and the symbol $NoRepresentation$ otherwise.

\begin{Remark}[Implementation of the matrix factorization]
In step \ref{alg:svd:step1} of Algorithm \ref{PartRepr}
above one can use any algorithm for computing a factorization.
For example, one could use SVD decomposition, in which case
$H_{\Psi,N+1,N}=U\Sigma V^{T}$, and
$O=U(\Sigma^{1/2})$, $R=(\Sigma^{1/2})V^{T}$ is a valid choice
for decomposition. 
\end{Remark}
The following theorem characterizes the outcome Algorithm \ref{PartRepr}.
\begin{Theorem}
\label{sect:svd:theo1}
 Let $\Psi=\{S_{j} \in \mathbb{R}^{p}\ll X^{*} \gg \mid j \in J\}$ be a family of
 formal power series.
With the notation above the following holds.
\begin{enumerate}
\item 
\label{sect:svd:theo1:part1}
  Assume that form some $N>0$
  \begin{equation} 
  \label{sect:svd:theo1:eq1}
  \Rank H_{\Psi,N,N+1}=\Rank H_{\Psi,N+1,N}=\Rank H_{\Psi,N,N}
  \end{equation}
  Then Algorithm \ref{PartRepr} always returns a 
  formal power series representation $\widetilde{R}_{N}$ and $\widetilde{R}_{N}$ is an 
  $2N+1$-partial representation of $\Psi$. In fact,
  the representation $R_{N}$ from Theorem \ref{part_real_pow:theo} and
  $\widetilde{R}_{N}$ are isomorphic. Hence, $\widetilde{R}_N$ is a 
  minimal $2N+1$-partial representation of $\Psi$, it is reachable and observable.

\item
\label{sect:svd:theo1:part2}
   If for some $N >0$,
   \begin{equation}
   \label{sect:svd:theo1:eq2}
   \Rank H_{\Psi,N,N}=\Rank H_{\Psi}
   \end{equation}
    then (\ref{sect:svd:theo1:eq1}) holds, and Algorithm \ref{PartRepr} returns 
   a minimal representation 
   $\widetilde{R}_{N}$ 
   of $\Psi$.
\item
\label{sect:svd:theo1:part3}
   Assume $\Rank H_{\Psi} \le N+1$, or , equivalently, there exists a 
   representation $R$ of $\Psi$, such that $\dim R \le N+1$.
   Then (\ref{sect:svd:theo1:eq2}) holds and 
   the representation $\widetilde{R}_{N}$ returned by Algorithm \ref{PartRepr} is a minimal
   representation of $\Psi$.
\end{enumerate}
\end{Theorem}
\begin{Remark}[Solution to (\ref{ComputePartialRepr})]
  Although solution to (\ref{ComputePartialRepr}) need not always exist, one can always 
  take a matrix $\widetilde{A}_{\sigma}$, $\sigma \in X$ as a solution to the minimization problem
   \( \min_{\widetilde{A}_{\sigma}} ||\bar{\Gamma}\widetilde{A}_{\sigma}-\bar{\Gamma}_{\sigma}||_{2} \).
   Then $\widetilde{A}_{\sigma}$ can
  be obtained using standard numerical techniques for 
  solving approximation problems.
  With this modification, the algorithm can be applied even if (\ref{sect:svd:theo1:eq1}) 
  fails. However,
  the representation returned by the modified 
  algorithm need not be a $2N+1$-partial representation
  in this case.
\end{Remark}

 \subsection{ Proof of the partial-realization results for formal power series}
\label{app:proof}

The goal of the section is to present the proof of Theorem \ref{part_real_pow:theo}, \ref{part_real_pow:theo2} and \ref{sect:svd:theo1}. 
In Subsection \ref{part_pow_pf:prelim} we will introduce some notation
and state some preliminary results, which will be needed for the proof of the theorems.
In Subsection \ref{part_pow_pf:pf} we will present the proof of 
Theorem \ref{part_real_pow:theo}--\ref{sect:svd:theo1}. Finally, in
Subsection \ref{part_pow_pf:aux} we will present the proof of the technical 
results which are used for the proof of the
Theorems \ref{part_real_pow:theo}--\ref{sect:svd:theo1}.


\subsubsection{ Auxiliary definitions and results}
\label{part_pow_pf:prelim}
 To begin with, for the purposes of partial realization theory, we need to recall some
   basic steps of the proof of Theorem \ref{sect:form:theo1}.
To this end, we have to introduce additional notation and terminology.
    Let $w \in X^{*}$  be a word over $X^{*}$ and let
    $S \in \mathbb{R}^{p}\ll X^{*} \gg$ be a formal power
    series.
    Define the formal power series
    $w \circ S \in \mathbb{R}^{p}\ll X^{*} \gg$,
    called the \emph{left shift of $S$ by $w$}, as follows; we require that for all
    $v \in X^{*}$ the value of $w \circ S$ at $v$ is as follows
     \begin{equation}
     \label{repr:eq6}
         (w \circ S)(v)=S(wv) 
     \end{equation}
  i.e. the value of $w \circ S$ at $v$ equals the value of $S$ at $wv$.
  Notice that for any word $w \in X^{*}$ of the
  form $w=\sigma_{1}\sigma_2 \cdots \sigma_{k}$, 
  $\sigma_{1},\sigma_2 \ldots, \sigma_{k} \in X$ and for any formal
  power series $T \in \mathbb{R}^{p}\ll X^{*} \gg$, the
  following equality holds
  \begin{equation} 
  \label{repr:eq7}
    w \circ T =\sigma_{k} \circ (\sigma_{k-1} \circ ( \cdots ( \sigma_{1} \circ T) \cdots )))
  \end{equation}
  Moreover, notice that the shift operation is linear, that
  is, for any $T,S \in \mathbb{R}^{p}\ll X^{*}\gg$,
  and for any scalars $\alpha,\beta \in \mathbb{R}$, and for
  any word $w \in X^{*}$,
  $w \circ (\alpha S + \beta T)=\alpha (w \circ S)+\beta (w \circ T)$.
    
  \begin{Definition}[Smallest shift invariant space]
  \label{sect:form:def1}
     Let $\Psi=\{ S_{j} \in \mathbb{R}^{p}\ll X^{*} \gg \mid j \in J\}$ be a family of formal power series.
     Define the \emph{the smallest shift invariant
     linear space  containing all the elements of 
     $\Psi$}, denoted by $W_{\Psi}$, 
     as the following subspace of $\mathbb{R}^{p}\ll X^{*} \gg$, 
     \begin{equation}
     \label{repr:eq8}
       W_{\Psi}= \SPAN \{ w \circ S_{j} \in \mathbb{R}^{p}\ll  X^{*}\gg \mid
        j \in J, 
               w \in X^{*} \}
     \end{equation}
 That is, $W_{\Psi}$ is composed of the linear combinations of
 all formal power series
 of the form $w \circ S_{j}$ for some $j \in J$ and $w \in X^{*}$.
 \end{Definition}
 \begin{Remark}
 \label{sect:form:def1:rem1}
 It is easy to see that there is one-to-one correspondence
 between the formal power series $w \circ S_{j}$ and
 the column of $H_{\Psi}$ indexed by $(w,j)$ for any word
 $w \in X^{*}$ and index $j \in J$. In particular, it
 follows that
    $W_{\Psi}$ is isomorphic to 
    the span of columns of $H_{\Psi}$ and hence
    \[ \dim W_{\Psi} = \Rank H_{\Psi} \]
\end{Remark}
 The statement of Theorem \ref{sect:form:theo1}  follows from the following
 two auxiliary statements, proofs of which can be found in \cite{MP:Phd,MP:BigArticle}
 \begin{Lemma}[\cite{MP:Phd,MP:BigArticle}]
 \label{sect:form:theo1:lemma1}
  Assume that $\dim W_{\Psi} < +\infty$ holds.
  Then a representation $R_{\Psi}$
  of $\Psi$ is given by
  \begin{equation} 
  \label{repr:eq9}
   R_{\Psi}=(W_{\Psi},\{A_{\sigma}\}_{\sigma \in X},B,C)
  \end{equation}
   where
   for each letter $\sigma \in X$, the map  
    $A_{\sigma}:W_{\Psi} \rightarrow W_{\Psi}$ is defined
    as the shift by $\sigma$, i.e. for each $T \in W_{\Psi}$, 
    $A_{\sigma}(T)=\sigma \circ T $;
    the collection 
    $B=\{ B_{j} \in W_{\Psi} \mid j \in J \}$ is formed
    by the elements of $\Psi$, that is
    $B_{j}=S_{j}$ for each $j \in J$;
    the linear map
    $C:W_{\Psi} \rightarrow \mathbb{R}^{p}$ is defined
    as the evaluation at the empty word, i.e.
    $C(T)=T(\epsilon)$ for all $T \in W_{\Psi}$.
\end{Lemma}    
\begin{Lemma}[\cite{MP:Phd,MP:BigArticle}]
\label{sect:form:theo1:lemma2}
  If $\Psi$ is rational, then $\dim W_{\Psi} < +\infty$. 
  More precisely, for each representation $R$ of $\Psi$,
  $\dim W_{\Psi} \le \dim R$.
\end{Lemma}
 Using the lemmas above the proof of Theorem
 \ref{sect:form:theo1} becomes trivial. For the sake of completeness we present
 it below.
 From Remark \ref{sect:form:def1:rem1} it follows that 
 $\dim W_{\Psi} = \Rank H_{\Psi}$. If $\Rank H_{\Psi} < +\infty$,
 then Lemma \ref{sect:form:theo1:lemma1} implies that
 $R_{\Psi}$ is a well-defined representation of $\Psi$,
 hence $\Psi$ is rational. Conversely, if $\Psi$ is rational
 then Lemma \ref{sect:form:theo1:lemma2} implies that
 $\dim W_{\Psi}=\Rank H_{\Psi} < +\infty$.
\begin{Remark}
 The representation $R_{\Psi}$ defined in (\ref{repr:eq9})
 is called \emph{free}.
\end{Remark}

Below we will define various spaces which play a role analogous to $W_{\Psi}$.
However, before proceeding further, some additional notation needs to be set up.
\begin{Notation}[Words of length at most $N$]
 We will denote the set of all words over $X$ of length at most $N$ by $X^{\le N}$, i.e.
 \( X^{\le N}=\{ w \in X^{*} \mid |w| \le N\} \).
\end{Notation}
We will start with defining the space $W_{\Psi,.,N}$
of formal power series which corresponds
to the subset of columns of $H_{\Psi}$ indexed by indices of the form
$(w,j)$ with $j \in J$ and $w \in X^{\le N}$.
  Define the vector space $W_{\Psi,.,N}$ as the sub-space of $W_{\Psi}$
  spanned by all formal power series
  of the form $w \circ S_{j}$ with $j \in J$ and $w \in X^{\le N}$, i.e.
  \begin{equation}
  \label{part_pow:pf1}
  W_{\Psi,.,N}=\SPAN\{ (w \circ S_{j}) \mid j \in J, w \in X^{*}, |w| \le N\}
  \end{equation}
\emph{Notice that if $J$ is finite, then $W_{\Psi,.,N}$ is generated by 
 finitely many elements}.
The motivation for considering the space $W_{\Psi,.,N}$ is revealed by the
following lemma and its corollary.
\begin{Lemma}
\label{part_pow:lemma1}
 Assume that $\Rank H_{\Psi}  \le N$. 
 For any formal power series $T \in W_{\Psi}$  and any word $w \in X^{*}$,
 there
 exists scalars $\alpha_{w,v} \in \mathbb{R}$, and words $v \in X^{\le N-1}$
 such that 
\[ w \circ T=\sum_{v \in X^{\le N-1}} \alpha_{w,v} ( v \circ T). \]
 That is, the shift of $T$ by $w$ is a linear combination of the shifts of $T$ by
 words of the length at most $N-1$.
\end{Lemma}
 The proof of Lemma \ref{part_pow:lemma1} will be presented in
 Subsection \ref{part_pow_pf:aux}.
\begin{Corollary}
\label{part_pow:pf1:col1}
 If $\Rank H_{\Psi} \le N+1$, then 
 \( W_{\Psi}=W_{\Psi,.,N} \).
\end{Corollary}
 The proof of Corollary \ref{part_pow:pf1:col1} will be presented in 
 Subsection \ref{part_pow_pf:aux}.
In other words, if $N$ is big enough, then $W_{\Psi,.,.N}$ generates the whole
space $W_{\Psi}$.  
Alternatively, stated in the language of Hankel-matrices,
the columns of $H_{\Psi}$ indexed by words of length at most $N$ span the
whole image of $H_{\Psi}$.

Although $W_{\Psi,.,N}$ is generated by finitely many elements, if $J$ is finite,
its elements are formal power series which contain infinite amount of data.
In order to replace $W_{\Psi,.,N}$ by a space generated by finite vectors, we
proceed as follows.
We will define a linear space $W_{\Psi,M,N}$ 
which plays a similar role for $H_{\Psi,M,N}$ as
$W_{\Psi}$ does for $H_{\Psi}$. 
Similarly to $W_{\Psi}$, the elements of $W_{\Psi,M,N}$ are maps mapping words over $X$
to vectors in $\mathbb{R}^{p}$. However, in contrast to $W_{\Psi}$, the elements of
$W_{\Psi,M,N}$ will be defined only on the words of length at most $M$. 
Similarly to $W_{\Psi}$, the space $W_{\Psi,M,N}$ will allow us to prove our main
results for the finite Hankel sub-matrices $H_{\Psi,M,N}$ in a more intuitive way.

To this end, we will introduce notation for restrictions of formal power series
$\mathbb{R}^{p}\ll X^{*} \gg$ to words of length at most $M$.
\begin{Notation}[Restriction of formal power series to $X^{\le M}$]
Denote by $\mathbb{R}^{p}\ll X^{\le M}\gg$ the set of 
functions $T: X^{\le M} \rightarrow \mathbb{R}^{p}$.
\end{Notation}
It is clear
that $\mathbb{R}^{p}\ll X^{\le M} \gg$ forms a vector space
with point-wise addition and point-wise multiplication by
scalar. 

As the next step, define the map
\[ \eta_{M}: \mathbb{R}^{p}\ll X^{*} \gg  \rightarrow \mathbb{R}^{p}\ll X^{\le M}\gg \]
 projecting any formal power series to its restriction to $X^{\le M}$, i.e.
for each formal power series
   $T \in \mathbb{R}^{p}\ll X^{*} \gg$, 
\begin{equation}
\label{project:def}
 \eta_{M}(T)(w)=T(w) \mbox{ for all } w \in X^{\le M} 
\end{equation}
It is easy to see that $\eta_{M}$ is a surjective linear map.
Define the vector space
$W_{\Psi,N,M}$ by
\begin{equation}
\label{part_pow:def2}
 W_{\Psi,M,N}=\SPAN \{ \eta_{M}(w \circ S_{j}) \mid w \in X^{\le N}, j \in J\}
\end{equation}
It is easy to see that $W_{\Psi,M,N}$ is the image of $W_{\Psi,.,N}$ by 
$\eta_{M}$.
The relationship between $W_{\Psi,M,N}$ and $H_{\Psi,M,N}$ can be
best described as follows; the generator set of $W_{\Psi,M,N}$ is formed
by the columns of $H_{\Psi,M,N}$. More precisely,
recall from Subsection \ref{sect:prelim:not_matrix} that 
$\mathbb{R}^{I_{M}}$ denotes the vector space
of all maps from $I_{M}$ to $\mathbb{R}$.  Define the map 
$\psi_{M}: \mathbb{R}^{p}\ll X^{\le M} \gg \rightarrow \mathbb{R}^{I_{M}}$
by 
\begin{equation} 
\label{part_pow:pf3.1}
(\psi_{M}(T))_{(v,i)}=(T(v))_{i} \mbox{ for all  } (v,i) \in I_{M}
\end{equation}
where $(T(v))_i$ stands for the $i$th entry of the vector $T(v) \in \mathbb{R}^{p}$.
It is easy to see that $\psi_{M}$ is a linear isomorphism. Moreover,
$\IM H_{\Psi,M,N}=\psi_{M}(W_{\Psi,M,N})$. 
In particular,
\begin{Corollary}
\label{part_pow:pf3:col3}
$\dim W_{\Psi,M,N}=\Rank H_{\Psi,M,N}$.
\end{Corollary}
Finally, we will show that if $\Psi$ is rational and $M$ is big enough, then the
restriction of $\psi_M$ to $W_{\Psi}$ is injective. It means then that $W_{\Psi,M,N}$ and
$W_{\Psi,.,N}$ are isomorphic for all $N$, if $M$ is large enough. 
\begin{Lemma}
\label{part_pow:lemma3}
  Assume that $\Rank H_{\Psi}=\dim W_{\Psi} \le M+1$. Then the restriction of $\eta_{M}$ to
  $W_{\Psi}$ is injective.
\end{Lemma}  

\subsubsection{ Proof of the main results on partial realization theory of formal power series}
\label{part_pow_pf:pf}
Now we are ready to present the proof of Theorem \ref{part_real_pow:theo} - \ref{sect:svd:theo1}.
\begin{pf}[Proof of Theorem \ref{part_real_pow:theo}]

\textbf{Proof of Part \ref{part_theo:part1} }
 We will define a $2N+1$ partial representation $\hat{R}_{N}$ of $\Psi$ on
 the linear space $W_{\Psi,N,N}$. The representation $\hat{R}_{N}$
 will be reachable and observable and will satisfy a number of properties.

 To this end, define the map 
 \begin{equation*}
    \eta_{N+1,N}: \mathbb{R}^{p}
       \ll X^{\le N+1} \gg \rightarrow \mathbb{R}^{p} \ll X^{\le N}\gg 
  \end{equation*}
 by requiring that for all $S \in \mathbb{R}^{p}\ll X^{\le N+1} \gg$ and for each word
 $w \in X^{\le N}$,
 \begin{equation} 
 \label{part_pow:pf6}
  \eta_{N+1,N}(S)(w)=S(w) 
  \end{equation}
It is easy to see $\eta_{N+1,N}$ is a surjective linear map. Moreover,
$\eta_{N+1,N}$ maps the sub-space $W_{\Psi,N+1,N}$ onto $W_{\Psi,N,N}$.
 Using Corollary \ref{part_pow:pf3:col3}, (\ref{part_theo:eq1})
 can be rewritten as 
 \begin{equation}
 \label{part_pow:pf4}
\dim W_{\Psi,N,N}=\dim W_{\Psi,N+1,N}=\dim W_{\Psi,N,N+1}
 \end{equation}
 But (\ref{part_pow:pf4}) implies that
 \( W_{\Psi,N,N+1}=W_{\Psi,N,N} \)
 and that the restriction of $\eta_{N+1,N}$ to 
 $W_{\Psi,N,N}$ is injective. The latter means that the inverse map 
 \begin{equation}
 \eta_{N+1,N}^{-1}: W_{\Psi,N,N} \rightarrow W_{\Psi,N+1,N}
 \end{equation}
 exists. 

 Define now the representation $\hat{R}_{N}$ as follows.
 \begin{equation}
 \label{part_pow:pf5}
  \hat{R}_{N}=(W_{\Psi,N,N+1}, \{ \hat{A}_{\sigma} \}_{\sigma \in X}, \hat{B},\hat{C})
 \end{equation}
 where
 \begin{itemize}
 \item
  For each letter $\sigma \in X$, the linear map
  \[ \hat{A}_{\sigma}:W_{\Psi,N,N}=W_{\Psi,N,N+1} \rightarrow W_{\Psi,N,N+1} \]
  is defined as follows.
  Consider the map $T_{\sigma}: W_{\Psi,N+1,N} \rightarrow \mathbb{R}^{p}\ll X^{\le N}\gg$ 
  defined by 
  $T_{\sigma}(Z)(v)=Z(\sigma v)$, for all $v \in X^{\le N}$ and
  $Z \in W_{\Psi,N+1,N}$ . It is easy to see that
  $T_{\sigma}$ is a linear map. In addition, if $Z=\eta_{N+1}(w \circ S_j)$ for
  some $w \in X^{\le N}$, $j \in J$, then  $T_{\sigma}(Z)=\eta_N(w\sigma \circ S_j)$.
  Indeed, for any $v \in X^{\le N}$,
  $T_{\sigma}(Z)(v)=w \circ S_{j}(\sigma v)=S_j(w\sigma v)=w \sigma \circ S_j(v)$. 
  Hence, because of linearity of $T_{\sigma}$, the range of $T_{\sigma}$ is 
  $W_{\Psi,N,N+1}$ and thus $T_{\sigma}$ can be viewed as a map 
  $T_{\sigma}:W_{\Psi,N+1,N} \rightarrow W_{\Psi,N,N+1}$.

   Then for all $S \in W_{\Psi,N,N}=W_{\Psi,N,N+1}$, define
   \[ \hat{A}_{\sigma}(S)=T_{\sigma}(\eta_{N+1,N}^{-1}(S)) \]
   It is easy to see that $\hat{A}_{\sigma}$ is a well-defined linear map and it has the 
   property that it commutes with the shift by $\sigma$ of a formal
   power series in $W_{\Psi,.,N}$, i.e. for all $S \in W_{\Psi,.,N}$,
   \begin{equation} 
   \label{part_pow:pf10}
    \hat{A}_{\sigma}(\eta_{N}(S)) =\eta_{N}(\sigma \circ S)
   \end{equation}
    Here, we used the notation of (\ref{part_pow:pf1}).
    Indeed, notice that $\eta_{N+1,N}^{-1}(\eta_{N}(S))=\eta_{N+1}(S)$.
    In addition, from the discussion above it follows that 
    $T_{\sigma}(\eta_{N+1}(S))=\eta_{N}(\sigma \circ S)$. Combining this with
    the definition of $A_{\sigma}$ yields (\ref{part_pow:pf10}).

    In addition, for all $S \in W_{\Psi,N,N}$, and for all $v \in X^{\le N-1}$ and 
    $\sigma \in X$, 
    \begin{equation}
    \label{part_pow:pf10.1}
       \{\hat{A}_{\sigma}(S)\}(v)=\{T_{\sigma}(\eta_{N+1,N}^{-1}(S))\}(v)=
        \{\eta_{N+1,N}^{-1}(S)\}(\sigma v) = S(\sigma v)
    \end{equation}
 \item
   The family $\hat{B}=\{ \hat{B}_{j} \in W_{\Psi,N,N} \mid j \in J\}$ is defined as
   \begin{equation}
   \label{part_pow:pf11}
        \hat{B}_{j}=\eta_{N}(S_{j}) \mbox{ for all } j \in J 
   \end{equation}
   i.e. $\hat{B}_{j}$ is just the restriction of $S_j$ to $X^{\le N}$.
 \item
   The map $\hat{C}$ is defined as
   \begin{equation} 
   \label{part_pow:pf12}
    \hat{C}: W_{\Psi,N,N+1} \ni S \mapsto S(\epsilon) \in \mathbb{R}^{p} 
   \end{equation}
   i.e. $\hat{C}$ is just the evaluation of the elements of $W_{\Psi,N+1,N}$ at
   the empty word $\epsilon$.
 \end{itemize}
  It is easy to see that $\hat{R}_{N}$ is a well-defined rational representation.
  Next, we will show that $\hat{R}_{N}$ is an $2N+1$-partial representation of $\Psi$.
To this end, by repeated application of (\ref{part_pow:pf10}) and using the equality
 $W_{\Psi,N,N}=W_{\Psi,N,N+1}$  we get that
 for any word $w \in X^{\le N+1}$ it holds that
\begin{equation}
\label{part_pow:pf8}
\hat{A}_{w}\hat{B}_{j}=\hat{A}_{w}(\eta_{N}(S_{j}))=\eta_{N}(w \circ S_{j})
\end{equation}
 Using (\ref{part_pow:pf8}) and repeatedly applying 
 (\ref{part_pow:pf10.1}) yields that for each $v \in X^{\le N}$,
\begin{equation}
\label{part_pow:pf13}
     \hat{C}\hat{A}_{v}\hat{A}_{w}\hat{B}_{j}=
     \{\hat{A}_{v}\eta_{N}(w \circ S_{j})\}(\epsilon)=\{w \circ S_{j}\}(v)=S_{j}(wv) 
\end{equation}
By noticing that $\hat{A}_{v}\hat{A}_{w}=\hat{A}_{wv}$ and that any word $\hat{w}$ in
$X^{\le 2N+1}$ can be represented as the concatenation $\hat{w}=wv$ 
of a word $w \in X^{\le N+1}$ with
a word $v \in X^{\le N}$, we get that for all $\hat{w} \in X^{\le 2N+1}$, 
\[ CA_{\hat{w}}\hat{B}_j = S_{j}(\hat{w}) \]
That is, $\hat{R}_{N}$ is an $2N+1$-partial representation of $\Psi$.
From (\ref{part_pow:pf8}) it is easy to deduce that $\hat{R}_{N}$ is reachable.
Observability of $\hat{R}_{N}$ can be derived as follows. 
From (\ref{part_pow:pf10.1}) it follows that for any 
$T \in W_{\Psi,N,N+1}=W_{\Psi,N,N}$, $\hat{C}\hat{A}_{w}T=T(w)$ for all words
$w \in X^{\le N}$.
Hence, if $\hat{C}\hat{A}_{w}(T)=0$ for all $w \in X^{*}$, then
for all $w \in X^{\le N}$, $T(w)=0$, i.e. $T=0$. That is, $O_{\hat{R}_{N}}=\{0\}$.

Define now the representation $R_{N}$ as the isomorphic copy of $\hat{R}_{N}$ with
the isomorphism $\psi_{N}$ defined in (\ref{part_pow:pf3.1}).
That is, using the notation of the theorem, define the parameters $A_{\sigma}$,
$\sigma \in X$, $C$ and $\{B_{j} \mid j \in J\}$ of $R_{N}$ as follows.
\begin{equation}
\label{part_pow:pf14}
\begin{array}{rcl}
 A_{\sigma}=\psi_{N} \hat{A}_{\sigma} \psi_{N}^{-1}  \mbox{ for all } \sigma \in X,  & 
 C=\hat{C}\psi_{N}^{-1},  &   B_{j}=\psi_{N}(\hat{B}_{j}) \mbox{ for all } j \in J.
\end{array}
\end{equation}
Since $\psi_{N}$ defines a representation isomorphism 
$\psi_{N}:\hat{R}_{N} \rightarrow R_{N}$,
it follows that $R_{N}$  is a $2N+1$-rational 
representation of $\Psi$ as well, and it is reachable and 
observable. In addition, combining (\ref{part_pow:pf14}) with
(\ref{part_pow:pf10},\ref{part_pow:pf11},\ref{part_pow:pf12}) yields
(\ref{part_theo:eq4},\ref{part_theo:eq5},\ref{part_theo:eq6}).

\textbf{Proof of Part \ref{part_theo:part2} }\\
 Recall that $\Rank H_{\Psi,M,K}=\dim W_{\Psi,M,K}$ for all $K,M \in \{N,N+1\}$
 and $\dim W_{\Psi}=\Rank H_{\Psi}$.
 It is easy to see that 
 \begin{equation} 
 \label{part_pow:pf-1}
   \dim W_{\Psi,N,N} \le \dim W_{\Psi,K,M} \le \dim W_{\Psi} 
    \mbox{ for all } K,M \in \{N,N+1\}
 \end{equation}
 Hence, (\ref{part_theo:eq2}) implies that  
 \begin{equation*} 
   \dim W_{\Psi,N,N}=\dim W_{\Psi,M,K}=\dim W_{\Psi} \mbox{ for all } M,K \in \{N,N+1\}
 \end{equation*}
  Hence, (\ref{part_theo:eq1}) holds and by Part \ref{part_theo:part1} the representation
  $R_N$ is well-defined and it is a $2N+1$-partial representation of $\Psi$.
  Recall that (\ref{part_theo:eq1}) implies that $W_{\Psi,N,N+1}=W_{\Psi,N,N}$.
  
  It is left to show that $R_{N}$ is a representation of $\Psi$. To this end, recall
  from (\ref{project:def}) the definition of the map $\eta_{N}$ and recall from 
  (\ref{part_pow:pf1})
  the definition of the set $W_{\Psi,.,N}$.
 Since $\eta_{N}(W_{\Psi,.,N})=W_{\Psi,N,N}$, we get that 
 \[ \dim W_{\Psi} \ge \dim W_{\Psi,.,N} \ge \dim W_{\Psi,N,N}=\dim W_{\Psi}, \] 
 hence $W_{\Psi,.,N}=W_{\Psi}$.
 Notice that $\eta_{N}$ maps $W_{\Psi,.,K}$ onto $W_{\Psi,N,K}$ for $K \in \{N,N+1\}$. 
 Hence,  
 the map $\eta_{N}$ maps $W_{\Psi}$ onto
 $W_{\Psi,N,N}=W_{\Psi,N,N+1}$.
 That is, the restriction $\eta_{N}$ to $W_{\Psi}$ yields a linear isomorphism
 \( \eta_{N}|_{W_{\Psi}}:W_{\Psi} \rightarrow W_{\Psi,N,N+1} \).
Recall from Lemma \ref{sect:form:theo1:lemma1} 
the definition of the free representation $R_{\Psi}$
of $\Psi$, defined on the space $W_{\Psi}$. 
Using (\ref{part_pow:pf10}),(\ref{part_pow:pf11}) and (\ref{part_pow:pf12}),
it is easy to see that the
restriction of $\eta_{N}$ to $W_{\Psi}$ in fact yields a representation isomorphism
$\eta_{N}|_{W_{\Psi}}:R_{\Psi} \rightarrow \hat{R}_{N}$, where
$\hat{R}_{N}$ is as defined in (\ref{part_pow:pf5}). Since, $R_{\Psi}$ is a minimal
representation of $\Psi$, then so is its isomorphic copy $\hat{R}_{N}$. Since
$R_{N}$ is merely an isomorphic copy of $\hat{R}_{N}$, the same conclusion holds
for $R_{N}$.

\textbf{Proof of Part \ref{part_theo:part3}} \\
  Recall from (\ref{part_pow:pf1}) and (\ref{part_pow:def2}) 
  the definition of the spaces
  $W_{\Psi,.,N}$ and $W_{\Psi,M,N}$ for $M \in \mathbb{N}$, $M > 0$.
  From Lemma \ref{sect:form:theo1:lemma2} it follows that if $\dim R \le N+1$,
  then $\dim W_{\Psi}=\Rank H_{\Psi} \le \dim R \le N+1$.
  Then by Corollary \ref{part_pow:pf1:col1},
  $W_{\Psi,.,K}=W_{\Psi}$ for $K=N,N+1$.
  Moreover, notice that the image of $W_{\Psi,.,N}$ by
  $\eta_{M}$ equals 
  $W_{\Psi,M,N}$.
  Applying the above observations for $M,K \in \{N,N+1\}$ and using 
  Lemma \ref{part_pow:lemma3} 
  we get
  that $\eta_{N}$ and $\eta_{N+1}$ are injective and hence
  \begin{equation}
  \label{part_pow:pf2} 
   \dim W_{\Psi}=\dim W_{\Psi,.,K}=\dim W_{\Psi,M,K}
    \mbox{ where } K,M=N,N+1
  \end{equation}
  Combining (\ref{part_pow:pf2}) with Corollary \ref{part_pow:pf3:col3} and 
  Remark \ref{sect:form:def1:rem1} 
  we get that (\ref{part_theo:eq2}) holds. Combining this and
  Part \ref{part_theo:part2} we get that the statement of Part \ref{part_theo:part3}
  holds.
  
\end{pf}
Next, we will present the proof of Theorem \ref{part_real_pow:theo2}.
\begin{pf}[Proof of Theorem \ref{part_real_pow:theo2}]
  We will prove the statements of the theorem one by one. However, before proceeding to
  the actual proof, we need to introduce some notation. Let $R$ be $p$-$J$ rational
  representation, i.e. $R$ is of the form (\ref{repr:def0}). Define the
  family of formal power series \emph{$\Psi_{R}$ associated with the representation $R$}
  as follows.
  The family $\Psi_{R}$ is indexed by elements of $J$, i.e. it can be written as
  $\Psi_{R}=\{S_{j}^{R} \in \mathbb{R}^{p}\ll X^{*} \gg \mid j \in J\}$. For each index
  $j \in J$, and for each word $w \in X^{*}$, the value of the formal power series
  $S_j^{R}$ at $w$ is defined as \( S^{R}_{j}(w)=CA_wB_{j} \).
  It is easy to see that the rational representation $R$ is a representation of $\Psi_{R}$,
  i.e. $\Psi_{R}$ is rational. In fact,  $R$ is a representation of the  family of formal
  power series $\Psi$, if and only if $\Psi_R=\Psi$, i.e. $S_{j}=S^{R}_{j}$ for all $j \in J$.
  Furthermore,
  $R$ is a $2N+1$-partial representation of $\Psi$, if and only if 
  $S_j^{R}(w)=S_j(w)$ for all $j \in J$ and $w \in X^{\le 2N+1}$.
  In addition, if $R$ is a $2N+1$-partial representation of $\Psi$, then 
  \begin{equation}
  \label{part_real_pow:theo2:pfeq1}
   H_{\Psi,K,L}=H_{\Psi_{R},K,L} \mbox{ for all } K+L \le 2N+1 
  \end{equation}
  In particular,
  $H_{\Psi,N,N}=H_{\Psi_R,N,N}$.
  Now we are ready to proceed to the actual proof of the theorem.

 \textbf{Proof of Part \ref{part_real_pow:theo2:res1}} \\ 
   Assume that $R$ is a $2N+1$-partial representation of $\Psi$. From the discussion above
   it follows that $\Rank H_{\Psi,N,N}=\Rank H_{\Psi_R,N,N}$. It is easy to see
   that $\Rank H_{\Psi_R,N,N} \le \Rank H_{\Psi_R}$. Moreover, from Theorem \ref{sect:form:theo3} it follows
   that $\Rank H_{\Psi_R} \le \dim R$. That is, the rank of $H_{\Psi,N,N}$ is not greater than
   the dimension of the representation $R$.
   But from Theorem \ref{part_real_pow:theo} it follows the the dimension of the
   representation $R_N$ is precisely $\Rank H_{\Psi,N,N}$ and that $R_N$ is an
   $2N+1$-partial  representation of $\Psi$.  Hence, we get that 
   $\dim R_N \le \dim R$. In other words,  $R_N$ is a minimal $2N+1$ representation of $\Psi$.

 \textbf{Proof of Part \ref{part_real_pow:theo2:res2}} \\
   We have just shown that the rational representation $R_N$ of Theorem
   \ref{part_real_pow:theo} is a minimal $2N+1$-partial representation of $\Psi$.
   If $R$ is another minimal $2N+1$-partial representation of $\Psi$, then
   the dimension of $R$ must be equal to that of $R_N$, which in turn equals 
   $\Rank H_{\Psi,N,N}$. Hence, the second part of the statement is proved.
   Assume now that $R$ is a minimal $2N+1$-partial representation of $\Psi$. Then
   it means that $S_{j}^{R}(w)=S_j(w)$ for each $j \in J$ and $w \in X^{\le 2N+1}$,
   where $S_j^{R}$ denotes the element of $\Psi_{R}$ indexed by the index $j \in J$.
   If $R$ is not reachable and observable, then it is not a minimal representation of 
   $\Psi_{R}$. From Theorem \ref{sect:form:theo3} it follows then that there exists a reachable
   and observable representation $R_m$ of $\Psi_R$ such that $\dim R_m < \dim R$.
   But if $R_m=(\mathcal{X}^{m},\{A^{m}_{\sigma}\}_{\sigma \in X}, B^{m},C^{m})$ 
   is a representation of $\Psi_R$, then it is a $2N+1$-partial representation of $\Psi$,
   since, $S_j(w)=S_j^{R}(w)=C^{m}A_{w}^{m}B^{m}$ for all $j \in J$ and $w \in X^{\le 2N+1}$.
   Hence, $R$ is not a minimal $2N+1$-partial representation of $\Psi$, which is a contradiction.

 \textbf{Proof of Part \ref{part_real_pow:theo2:res3}} \\
    Suppose $R$ and $\hat{R}$ are two minimal $2N+1$-partial representations of $\Psi$.
    Recall the definition of the families of formal power series
    $\Psi_{R}$ and $\Psi_{\hat{R}}$ associated with $R$. Since $R$ and $\hat{R}$ have
    to be reachable and observable by Part \ref{part_real_pow:theo2:res2} of the
    current theorem, by Theorem \ref{sect:form:theo3}
   the representation $R$ (respectively $\hat{R}$) is a minimal
    representation of $\Psi_R$ (respectively $\Psi_{\hat{R}}$).
    We will show that $\Psi_R$ and $\Psi_{\hat{R}}$ are equal, i.e.
    $S_j^{R}=S^{\hat{R}}_j$ for all $j \in J$. This in turn implies that $R$ and $\hat{R}$
    are both minimal representations of the same family of formal power series, and 
    hence by Theorem \ref{sect:form:theo3} they are isomorphic.
    
    The proof that $\Psi_R$ and $\Psi_{\hat{R}}$ are identical proceeds as follows.
    Recall that $R$ and $\hat{R}$ are both of dimension $\Rank H_{\Psi,N,N}$, and
    recall (\ref{part_real_pow:theo2:pfeq1}).
    Hence, it follows that $\Rank H_{\Psi,N,N}=\Rank H_{\Psi_{R},N,N}=H_{\Psi_R}=\dim R$ and 
    $\Rank H_{\Psi,N,N}=\Rank H_{\Psi_{\hat{R}},N,N}=\Rank H_{\Psi_{\hat{R}}}=\dim \hat{R}$. 
    It also follows that 
    $\Psi_R$ and $\Psi_{\hat{R}}$ satisfy  (\ref{part_theo:eq2}) with $\Psi$ replaced
    by $\Psi_R$ and $\Psi_{\hat{R}}$ respectively. 
    Denote by $R^{1}_{N}$ and $R^{2}_{N}$ the minimal $2N+1$-partial
    representations of $\Psi_R$ and respectively $\Psi_{\hat{R}}$, obtained
    by applying Part \ref{part_theo:part1} of
    Theorem \ref{part_real_pow:theo} to $\Psi_{R}$ and $\Psi_{\hat{R}}$
    respectively.
    By Part \ref{part_theo:part2}
    of Theorem \ref{part_real_pow:theo} the representations $R^{1}_{N}$ and $R^{2}_N$
    are minimal representations of $\Psi_R$ and $\Psi_{\hat{R}}$ respectively.
    However the construction of $R^{1}_{N}$ and $R^{2}_N$ depends only on
    $H_{\Psi_R,N,N+1}$ and $H_{\Psi_{\hat{R}},N,N+1}$. The latter two matrices  both coincide
    with $H_{\Psi,N,N+1}$. Hence, $R^{1}_{N}$ and $R^{2}_N$ are both equal to the
    partial representation $R_N$ obtained by applying Part \ref{part_theo:part1} of
    Theorem \ref{part_real_pow:theo}
    to $\Psi$.
    That is,  rational representation $R_{N}$, obtained by applying 
    Part \ref{part_theo:part1} of Theorem \ref{part_real_pow:theo} to $\Psi$, is a 
    minimal representation of both $\Psi_R$ and $\Psi_{\hat{R}}$. This implies that
    $\Psi_R$ and $\Psi_{\hat{R}}$ coincide.
\end{pf}

We conclude with the proof of Theorem \ref{sect:svd:theo1}.
\begin{pf}[Proof of Theorem \ref{sect:svd:theo1}]

\textbf{Proof of Part \ref{sect:svd:theo1:part1}}\\   
 We have to show that Algorithm \ref{PartRepr} returns a representation $\widetilde{R}_{N}$, 
 and $\widetilde{R}_{N}$ is isomorphic to $R_{N}$. 
 By Theorem \ref{part_real_pow:theo}, if (\ref{sect:svd:theo1:eq1}) 
 holds then $R_{N}$ is an $2N+1$-representation of $\Psi$. From Theorem \ref{part_real_pow:theo2}
 it then follows that $R_N$ is a minimal $2N+1$-partial representation. 
 Hence, if $\widetilde{R}_{N}$ is an isomorphic copy of $R_N$, then we get that
 $\widetilde{R}_{N}$ is a minimal $2N+1$-partial representation of $\Psi$.

 Recall from (\ref{part_pow:pf5}) in the proof of Theorem \ref{part_real_pow:theo} the definition
 of the representation $\hat{R}_{N}$.
 Using the remark in Subsection \ref{sect:prelim:not_matrix}, we will interpret
 the matrix $O \in \mathbb{R}^{I_{N+1} \times r}$ as a linear map
 $O:\mathbb{R}^{r} \rightarrow \mathbb{R}^{I_{N+1}}$ defined by
 $(Ox)_{i}=\sum_{k=1}^{r} O_{i,k}x_k$ for all $x \in \mathbb{R}^{r}$.
 Using this interpretation, we define the map 
 \begin{equation} 
 \label{sect:svd:theo1:pf1}
   \hat{\xi}=\eta_{N,N+1} \circ \psi_{N+1}^{-1} \circ O: \mathbb{R}^{r} \rightarrow W_{\Psi,N,N}=W_{\Psi,N,N+1}
 \end{equation}
 Since $\psi_{N}$ defines a representation isomorphism 
 $\psi_{N}:\hat{R}_{N} \rightarrow R_{N}$ 
 it is enough to show that $\hat{\xi}$ defines a 
 representation isomorphism $\hat{\xi}: \widetilde{R}_{N} \rightarrow \hat{R}_{N}$
 and that Algorithm \ref{PartRepr} returns a rational representation if
 (\ref{sect:svd:theo1:eq1}) holds.
 Indeed, in this case $\xi=\psi_{N} \circ \hat{\xi}$ defines a representation
 isomorphism $\xi:\widetilde{R}_{N} \rightarrow R_{N}$.

 It is clear that $\hat{\xi}$ is well defined. Indeed, 
 $\IM O=\IM H_{\Psi,N+1,N}$ by definition of
 matrix factorization. Moreover,  $O: \mathbb{R}^{r} \rightarrow \mathbb{R}^{I_{N+1}}$ 
 is an injective linear map and
 \[ \psi_{N+1}: W_{\Psi,N+1,N} \rightarrow \IM H_{\Psi,N+1,N} \]
 is 
 a linear isomorphism. Furthermore, since 
 \[ \Rank H_{\Psi,N+1,N}=\dim W_{\Psi,N+1,N}=\Rank H_{\Psi,N,N}=\dim W_{\Psi,N,N} \]
 we get that $\eta_{N+1,N}: W_{\Psi,N+1,N} \rightarrow W_{\Psi,N,N}$
 is a linear isomorphism. Thus, 
 $\hat{\xi}=\eta_{N+1,N} \circ \psi_{N+1}^{-1} \circ O$ 
  is a well defined linear isomorphism.
 It is left to show that $\hat{\xi}$ is a representation morphism from $\widetilde{R}_{N}$
 to $\hat{R}_{N}$.

 It is easy to see that for all $w \in X^{\le N}$, and $x \in \mathbb{R}^{r}$.
 \begin{equation}
 \label{sect:svd:theo1:pf:eq2}
 \hat{\xi}(x)(w)= \begin{bmatrix} (Ox)_{(w,1)}, & (Ox)_{(w,2)}, & \cdots & (Ox)_{(w,p)}
                  \end{bmatrix}^{T}
 \end{equation}    
  Recall from Subsection \ref{sect:prelim:not_matrix}
  that $(Ox)_{(w,i)}$ stands for the value (entry) of $Ox \in \mathbb{R}^{I_{N+1}}$
  corresponding to the index $(w,i) \in I_{N+1}$.
  In particular, for any $(v,j) \in J_{N}$,
  \begin{equation}
  \label{sect:svd:theo1:eq2.0}
   \hat{\xi}(R_{.,(v,j)})=\eta_{N}(v \circ S_{j}) 
  \end{equation}
  where $R_{.,(v,j)}$ stands for the column of $R$ indexed by $(v,j)$.
  Indeed, using (\ref{sect:svd:theo1:pf:eq2}) we get that for all $w \in X^{\le N}$,
  \begin{equation} 
  \label{sect:svd:theo1:eq2.0.0}
    \hat{\xi}(R_{.,(v,j)})(w)=
                  \begin{bmatrix} (OR)_{(w,1),(v,j)}, & (OR)_{(w,2),(v,j)}, & \cdots & 
                (OR)_{(w,p),(v,j)}
                  \end{bmatrix}^{T}
  \end{equation}
  But $H_{\Psi,N+1,N}=OR$, and hence the right-hand side of
  (\ref{sect:svd:theo1:eq2.0.0}) equals the vector
  formed by the entries $(H_{\Psi,N+1,N})_{(w,i),(v,j)}$ for $i=1,2\ldots,p$.
  But the entry $(H_{\Psi,N+1,N})_{(w,i),(v,j)}$ equals the $i$th entry of $S_{j}(vw)$, hence
  \[ \hat{\xi}(R_{.,(v,j)})(w)=S_{j}(vw)=v \circ S_j(w) \]
 From (\ref{sect:svd:theo1:pf:eq2}) it follows that
 \begin{equation}
  \label{sect:svd:theo1:eq2.1}
 \begin{split}
   & \hat{C}(\hat{\xi}(x))=\hat{\xi}(x)(\epsilon)=\widetilde{C}x	
\end{split}
\end{equation}
 From (\ref{sect:svd:theo1:eq2.0}), 
 it follows that 
 \begin{equation}
 \label{sect:svd:theo1:eq2.2}
 \begin{split}
    \hat{\xi}(\widetilde{B}_{j})=\eta_{N}(S_j)=\hat{B}_{j}
 \end{split}
 \end{equation}

It is left to show that a unique solution to equation (\ref{ComputePartialRepr})
exists and 
\begin{equation} 
\label{sect:svd:theo1:eq5}
   \hat{A}_{\sigma}\hat{\xi}=\hat{\xi} \widetilde{A}_{\sigma}
\end{equation}
 for all $\sigma \in X$. 
First, notice that $\bar{\Gamma}R=H_{\Psi,N,N}$. Thus,
$\Rank \bar{\Gamma}R=\Rank H_{\Psi,N,N}=r$, i.e. $\Rank \bar{\Gamma}=r$.
Hence, if a solution solution (\ref{ComputePartialRepr})
exists, then this solution is unique. Therefore, if we show that
$\hat{\xi}^{-1} \hat{A}_{\sigma} \hat{\xi}$ is a solution to 
(\ref{ComputePartialRepr}) then 
(\ref{sect:svd:theo1:eq5}) automatically holds.
Notice that $\hat{\xi}^{-1}\hat{A}_{\sigma}\hat{\xi}$ is a linear map from
$\mathbb{R}^{r}$ to $\mathbb{R}^{r}$ and hence it can be identified with its 
matrix representation.
By identifying $\hat{\xi}^{-1}\hat{A}_{\sigma}\hat{\xi}$ with its matrix, and
using (\ref{sect:svd:theo1:eq2.0}) and (\ref{part_pow:pf10}),
it follows that for any $(v,j) \in J_{N}$ and $(w,i) \in I_{N}$,
\begin{equation}
\label{sect:svd:theo1:eq6}
\begin{split}
 & (\bar{\Gamma}\hat{\xi}^{-1}\hat{A}_{\sigma} \hat{\xi}R)_{(w,i),(v,j)}=   
   \bar{\Gamma}_{(w,i),.} (\hat{\xi}^{-1}\hat{A}_{\sigma}\hat{\xi}R)_{.,(v,j)}= \\
   & \bar{\Gamma}_{(w,i),.}(\hat{\xi}^{-1}\hat{A}_{\sigma}\hat{\xi}R_{.,(v,j)})=
    \bar{\Gamma}_{(w,i),.} \hat{\xi}^{-1}\hat{A}_{\sigma}\eta_{N}(v \circ S_{j})=  \\
 & \bar{\Gamma}_{(w,i),.} \hat{\xi}^{-1}\eta_{N}(v\sigma\circ S_j)
\end{split}
\end{equation}
 Here we used the notation of Subsection \ref{sect:prelim:not_matrix} to denote rows
 and columns of matrices. In particular, $M_{.,(v,j)}$ denotes the column of
 a matrix $M$ indexed by $(v,j)$ and $M_{(w,i),.}$ denotes the row of a matrix $M$
 indexed by $(w,i)$.
 Notice that if $w \in X^{\le N}$, then the row of $\bar{\Gamma}$ indexed by $(w,i)$
 equals the row of $O$ indexed by $(w,i)$. Moreover, if
  $x=\hat{\xi}^{-1}(\eta_N(v\sigma \circ S_j))$, then from (\ref{sect:svd:theo1:pf:eq2}) 
  it follows that the row of $O(x)$ indexed
  by $(w,i)$ equals the $i$th entry of $\eta_N(v\sigma \circ S_j)(w)=S_{j}(v\sigma w)$.
 Hence, the last expression of (\ref{sect:svd:theo1:eq6}) can be rewritten as
 \begin{equation}
 \label{sect:svd:theo1:eq6.1}
 \bar{\Gamma}_{(w,i),.} \hat{\xi}^{-1}\eta_{N}(v\sigma\circ S_j) =
  O(\hat{\xi}^{-1}\eta_{N}(v\sigma\circ S_j))_{(w,i)}=(S_{j}(v\sigma w))_{i}
 \end{equation}
 where $(S_{j}(v\sigma w))_i$ denotes the $i$th entry of $S_{j}(v\sigma w) \in \mathbb{R}^{p}$.
On the other hand, the row of $\bar{\Gamma}_{\sigma}$ indexed by $(w,i)$
 equals the row of $O$ indexed by $(\sigma w,i)$, if $w \in X^{\le N}$. Hence, 
\begin{equation}
\label{sect:svd:theo1:eq7}
\begin{split}
  (\bar{\Gamma}_{\sigma}R)_{(w,i),(v,j)}=
  (O)_{(\sigma w,i),.}R_{.,(v,j)}=
 (H_{\Psi,N+1,N})_{(\sigma w,i),(v,j)} =(S_{j}(v\sigma w))_{i}
 \end{split}
\end{equation}
 Combining (\ref{sect:svd:theo1:eq6}), (\ref{sect:svd:theo1:eq6.1}) and 
 (\ref{sect:svd:theo1:eq7}) we get that 
 \( \bar{\Gamma}\hat{\xi}^{-1}\hat{A}_{\sigma}\hat{\xi}R=\bar{\Gamma}_{\sigma}R \).
  Since $\Rank R=r$, i.e the columns of $R$ span the whole space $\mathbb{R}^{r}$,
  the last equality implies that
  $\hat{\xi}^{-1}\hat{A}_{\sigma}\hat{\xi}$ is a solution to
  (\ref{ComputePartialRepr}).

\textbf{Proof of Part \ref{sect:svd:theo1:part2}}\\
 From Theorem \ref{part_real_pow:theo} it follows that 
 if (\ref{sect:svd:theo1:eq2}) holds, then  (\ref{sect:svd:theo1:eq1}) holds
 and $R_{N}$ is a minimal representation of $\Psi$.
 By Part \ref{sect:svd:theo1:part1} 
 of this theorem, if (\ref{sect:svd:theo1:eq1}) holds, then Algorithm \ref{PartRepr} returns a 
 representation $\widetilde{R}_{N}$ and $\widetilde{R}_{N}$ is isomorphic to $R_{N}$.
 Hence $\widetilde{R}_{N}$ is a minimal representation too.

\textbf{Proof of Part \ref{sect:svd:theo1:part3}}\\
  Again, from Theorem \ref{part_real_pow:theo} it follows
  that (\ref{sect:svd:theo1:eq2}) holds in this case. The statement follows now from
  (\ref{sect:svd:theo1:eq2}) and Part \ref{sect:svd:theo1:part2} of this theorem.
\end{pf}

\subsubsection{ Proof of the auxiliary results}
\label{part_pow_pf:aux}
We conclude the section with presenting the
proof of Lemma \ref{part_pow:lemma1} and Lemma \ref{part_pow:lemma3}
used in the proof of Theorem \ref{part_real_pow:theo}. 
The proof relies on the following chain of results, which are interesting
on their own right.
\begin{Lemma}
\label{part_pow:lemma4}
 Let $\mathcal{X}$ be finite-dimensional vector space, 
 $\dim \mathcal{X} \le N$. Let 
$A_{\sigma}: \mathcal{X} \rightarrow \mathcal{X}$, $\sigma \in X$ 
be a family of linear maps. Then for each
 $y \in \mathcal{X}$, for each $w \in X^{*}$, the vector $A_{w}y$
 is a linear combination of the vectors $A_{v}y$, for finitely many words 
$v \in X^{\le N-1}$.
\end{Lemma}
\begin{pf}
If $|w| < N$, then the statement of the lemma is trivially true.
 First we prove the lemma for $|w|=N$. Assume that 
 $w=\sigma_{1}\sigma_2\cdots \sigma_{N}$, $\sigma_{1},\sigma_2,\ldots, \sigma_{N} \in X$. 
Consider the elements $A_{\sigma_{1}\sigma_2\cdots \sigma_{i}}y$, $i=0,\ldots, N$.
Since every $N+1$ elements of $\mathcal{X}$ are linearly
dependent, we get that there exist 
$i=0,\ldots,N$ such that
$A_{\sigma_1\cdots \sigma_i}y$ is a linear combination of 
$A_{\sigma_1\cdots \sigma_j}y$, $j=0,2,\ldots,i-1$.
 Then we get that
\( A_{w}y=A_{\sigma_{i+1}\cdots \sigma_{N}}(A_{\sigma_{1}\cdots \sigma_{i}}y) \) is a linear
 combination of vectors of the form $A_{\sigma_1\cdots \sigma_j\sigma_{i+1}\cdots \sigma_{N}}y$,
$j=0,\ldots,i-1$. Notice that each word
$\sigma_1\cdots \sigma_j\sigma_{i+1}\cdots \sigma_N$ is of length at most $N-1$.
 To prove the lemma for arbitrary $|w| \ge N$ we proceed by
 induction on $|w|-N$. The case of $|w|=N$ we proved above. 
 Assume that the statement of the lemma holds for $|w|\le n+N$. 
 Assume now that $w$ is of the form $w=s\sigma$, where $s \in X^{*}$,
 $|s| \le n+N$ and $\sigma \in X$ is the last letter of $w$.
 From the induction hypothesis it follows that $A_{s}y$ is a linear
 combination of vectors of the form $A_{v}y$, $v \in X^{\le N-1}$.
 But then $A_{w}y=A_{\sigma}(A_{s}y)$ is a linear combination of 
 vectors of the form $A_{\sigma}A_{v}y=A_{v\sigma}y$. For each $v \in X^{\le N-1}$,
either $v\sigma \in X^{\le N-1}$ or $|v\sigma|=N$. In the latter case, by
induction hypothesis
 $A_{v\sigma}y$ is again a linear combination of vectors of the form
 $A_{\hat{v}}y$, $\hat{v} \in X^{\le N-1}$. Altogether, we get that 
 $A_{w}y$ is a linear combination of the vectors $A_{v}y$, $v \in X^{\le N-1}$.
\end{pf}
 The lemma above yields the following characterization of
 the reachability and observability
 subspaces (see Definition \ref{reach:rep:def} and \ref{obs:rep:def}) 
 of a rational representation.
\begin{Corollary}
\label{part_pow:col1}
Consider a $p-J$ representation $R=(\mathcal{X},\{A_{\sigma}\}_{\sigma \in X},B,C)$.
Assume that $\dim R \le N$. With the notation of
(\ref{sect:pow:reachobs:eq1.1}--\ref{sect:pow:reachobs:eq2.1}), the following holds.
\begin{equation*}
 \begin{split}
 O_{R}&=\bigcap_{v \in X^{\le N-1}} \ker CA_{v} \\
 W_{R}&=\SPAN\{ A_{v}B_{j} \mid j \in J, v \in X^{\le N-1}\} \\
 \end{split}
\end{equation*}
\end{Corollary}
 If $\mathcal{X}=\mathbb{R}^{n}$ and $J$ is finite, then
 the above corollary states that the observability subspace $O_R$ is the
 kernel of a finite matrix, and the reachability subspace $W_R$ is the image of
 a finite matrix.
\begin{pf}[Proof of Corollary \ref{part_pow:col1}]
 It is clear that $O_{R} \subseteq \bigcap_{v \in X^{\le N-1}} \ker CA_{v}$.
 We will show that the reverse inclusion holds as well.
 Assume that $x \in \mathcal{X}$ is such that
 $CA_{v}x=0$ for all $v \in X^{\le N-1}$.
 By Lemma \ref{part_pow:lemma4}, for each $w \in X^{*}$,
 $A_{w}x$ is a linear combination of $A_{v}x$ for $v \in X^{\le N-1}$. Hence,
 $CA_{w}x$ is a linear combination of $CA_{v}x$, $v \in X^{\le N-1}$ 
 and the latter implies that $CA_{w}x=0$. That is, 
 $\bigcap_{v \in X^{\le N-1}} \ker CA_{v} \subseteq O_{R}$.

 Similarly, it is clear from the definition that $\SPAN\{ A_{v}B_{j} \mid v \in X^{\le N-1}, j \in J\} \subseteq W_{R}$. The reverse inclusion follows from 
 Lemma \ref{part_pow:lemma4}, according to which for each word $w \in X^{*}$,
 $A_{w}B_{j}$ is a linear combination of the vectors $A_{v}B_{j}$, $v \in X^{\le N-1}$.
 Hence, we get that $W_{R} \subseteq \SPAN\{ A_{v}B_{j} \mid v \in X^{\le N-1}, j \in J\}$.
\end{pf}
 Finally, we will present the proof of Lemma \ref{part_pow:lemma1} and
 Lemma \ref{part_pow:lemma3}.
\begin{pf}[Proof of Lemma \ref{part_pow:lemma1}]
  We will use the fact that $\dim W_{\Psi}=\Rank H_{\Psi} \le N+1$.
  Consider the free representation 
  $R_{\Psi}=(W_{\Psi}, \{A_{\sigma}\}_{\sigma \in X}, B,C)$ of $\Psi$ defined in
  Theorem \ref{sect:form:theo1}. 
  Notice that for any word $X^{*}$, 
  $w \circ S_{j}=A_{w}B_{j}$.
  Apply Lemma \ref{part_pow:lemma4} to $W_{\Psi}$,
  $A_{\sigma}: W_{\Psi} \rightarrow W_{\Psi}$, $\sigma \in X$ and $y=B_{j}$.
  Then we get that for any word $w \in X^{*}$, the formal power series
   $A_{w}B_j=w \circ S_{j}$ 
  is a linear combination of the formal power series
   $A_{v}B_j=v \circ S_{j}$ with $v \in X^{\le N}$.
\end{pf}

\begin{pf}[Proof of Lemma \ref{part_pow:lemma3}]
 It is easy to see that $\eta_{M}$ is a surjective linear map. 
 Consider the 
 free realization $R_{\Psi}=(W_{\Psi}, \{ A_{\sigma} \}_{\sigma \in X}, B,C)$ of
 $\Psi$. 
 From Theorem \ref{sect:form:theo3} we know that $R_{\Psi}$ is minimal
 and therefore it is reachable and observable, i.e. 
 $O_{R_{\Psi}}=\{0\}$.
 From Corollary \ref{part_pow:col1} we also know that
 if $\dim R_{\Psi}=\dim W_{\Psi}=\Rank H_{\Psi} \le M+1$, then
 $O_{R_{\Psi}}=\bigcap_{v \in X^{\le M}} \ker CA_{v}$.
 Consider the kernel of $\eta_{M}$. For any formal power series $S \in W_{\Psi}$, 
 $\eta_{M}(S)=0$ if and only if
 $S(w)=0$ for all $w \in X^{\le M}$. Hence, $CA_{w}S=S(w)=0$ for each $w \in X^{\le M}$,
 i.e. $S \in O_{R_{\Psi}}=\{0\}$. Thus, $\eta_{M}$ is injective.
\end{pf}
\end{document}